\documentclass[11pt]{article}
\usepackage{txfonts}
\usepackage{color,epsfig}
\usepackage{mathrsfs}
\usepackage{amssymb}
\usepackage{amsfonts,graphicx,psfrag}

\newtheorem{lemma}{Lemma}[section]
\newtheorem{theorem}[lemma]{Theorem}
\newtheorem{corollary}[lemma]{Corollary}

\newtheorem{proposition}[lemma]{Proposition}
\newtheorem{definition}[lemma]{Definition}

\newtheorem{remark}[lemma]{Remark}
\let\lutzremark=\remark
\def\remark{\lutzremark\normalfont}

\newtheorem{notation}[lemma]{Notation}

\def\be{\begin{equation}}
\def\ee{\end{equation}}
\def\bea{\begin{eqnarray}}
\def\eea{\end{eqnarray}}
\def\bes{\begin{eqnarray*}}
\def\ees{\end{eqnarray*}}

\def\nn{\nonumber}
\def\<{\langle}
\def\>{\rangle}
\def\lb{\label}
\def\bs{\setminus}
\def\pt{\partial}

\def\R{{\bf R}}
\def\C{{\bf C}}
\def\Z{{\bf Z}}
\def\N{{\bf N}}
\def\U{{\bf U}}

\def\aa{{\alpha}}
\def\bb{{\beta}}
\def\ga{{\gamma}}
\def\Ga{{\Gamma}}

\def\th{{\theta}}
\def\om{{\omega}}
\def\Om{{\Omega}}
\def\ep{{\epsilon}}
\def\lm{{\lambda}}

\def\sg{{\sigma}}

\def\P{{\cal P}}

\def\Im{{\rm Im}}
\def\Re{{\rm Re}}

\def\Sp{{\rm Sp}}

\def\dm{{\rm \diamond}}

\def\hb{\vrule height0.18cm width0.14cm $\,$}
\def\ol#1{\overline{#1}}
\def\td#1{\tilde{#1}}

\title{Linear stability of elliptic Lagrangian solutions of\\
the planar three-body problem via index theory}

\author{Xijun Hu$^{1}$\thanks{Partially supported by NSFC (No.10801127, 11131004). E-mail: xjhu@sdu.edu.cn},
\quad Yiming Long$^{2}$\thanks{Partially supported by NSFC (No.11131004), MCME, RFDP, LPMC of MOE of China,
Nankai University, and the Beijing Center for Mathematics and Information Interdisciplinary Sciences at CNU.
E-mail: longym@nankai.edu.cn},
\quad Shanzhong Sun$^{3,4}$\thanks{Partially supported by NSFC (No.10731080, 11131004), PHR201106118, PCSIRT,
the Institute of Mathematics and Interdisciplinary Science at CNU.  E-mail: sunsz@mail.cnu.edu.cn }\\ \\
$^{1}$ Department of Mathematics, Shandong University\\
Jinan, Shandong 250100, The People's Republic of China\\
$^{2}$ Chern Institute of Mathematics and LPMC, Nankai University \\
Tianjin 300071, The People's Republic of China\\
$^{3}$ Department of Mathematics, Capital Normal University\\
Beijing 100048, The People's Republic of China\\
$^{4}$ Beijing Center for Mathematics and Information
Interdisciplinary Sciences\\
Beijing 100048, The People's Republic of China}
\date{}

\begin{document}

\maketitle

\begin{abstract}
{It is well known that the linear stability of Lagrangian elliptic equilateral triangle homographic
solutions in the classical planar three-body problem depends on the mass parameter
$\bb=27(m_1m_2+m_2m_3+m_3m_1)/(m_1+m_2+m_3)^2\in [0,9]$ and the eccentricity $e\in [0,1)$. We are not
aware of any existing analytical method which relates the linear stability of these solutions to the
two parameters directly in the full rectangle $[0,9]\times [0,1)$, aside from perturbation methods for $e>0$
small enough, blow-up techniques for $e$ sufficiently close to $1$, and numerical studies. In this paper,
we introduce a new rigorous analytical method to study the linear stability of these solutions in terms
of the two parameters in the full $(\bb,e)$ range $[0,9]\times [0,1)$ via the $\om$-index theory of
symplectic paths for $\om$ belonging to the unit circle of the complex plane, and the theory of linear
operators. After establishing the $\om$-index decreasing property of the solutions in $\bb$ for fixed
$e\in [0,1)$, we prove the existence of three curves located from left to right in the rectangle
$[0,9]\times [0,1)$, among which two are $-1$ degeneracy curves and the third one is the right envelop
curve of the $\om$-degeneracy curves, and show that the linear stability pattern of
such elliptic Lagrangian solutions changes if and only if the parameter $(\bb,e)$ passes through each of
these three curves. Interesting symmetries of these curves are also observed. The linear stability of the
singular case when the eccentricity $e$ approaches $1$ is also analyzed in detail. }
\end{abstract}

{\bf Keywords:} planar three-body problem, Lagrangian solution, linear stability, Maslov-type $\om$-index,
perturbations of linear operators.

{\bf AMS Subject Classification}: 58E05, 37J45, 34C25

\renewcommand{\theequation}{\thesection.\arabic{equation}}

\setcounter{equation}{0}
\setcounter{figure}{0}
\section{Introduction and main results}
\label{sec:1}

We consider the classical planar three-body problem in celestial mechanics. Denote by $q_1,q_2,q_3\in \R^2$
the position vectors of three particles with masses $m_1,m_2,m_3>0$ respectively. By Newton's second law
and the law of universal gravitation, the system of equations for this problem is
\be   m_i\ddot{q}_i=\frac{\partial U}{\partial q_i}, \qquad {\rm for}\quad i=1, 2, 3, \lb{1.1}\ee
where $U(q)=U(q_1,q_2,q_3)=\sum_{1\leq i<j\leq 3}\frac{m_im_j}{\|q_i-q_j\|}$ is the
potential or force function by using the standard norm $\|\cdot\|$ of vector in $\R^2$.
For periodic solutions with period $2\pi$, the system is the Euler-Lagrange equation
of the action functional
$$ \mathcal{A}(q)=\int_{0}^{2\pi}\left[\sum_{i=1}^3\frac{m_i\|\dot{q}_i(t)\|^2}{2}+U(q(t))\right]dt $$
defined on the loop space $W^{1,2}(\R/2\pi\Z,\hat{\mathcal {X}})$, where
$$  \hat{\mathcal {X}}:=\left\{q=(q_1,q_2,q_3)\in (\R^2)^3\,\,\left|\,\,
       \sum_{i=1}^3 m_iq_i=0,\,\,q_i\neq q_j,\,\,\forall i\neq j \right. \right\}  $$
is the configuration space of the planar three-body problem. The periodic solutions of (\ref{1.1}) correspond
to critical points of the action functional.

It is a well-known fact that (\ref{1.1}) can be reformulated as a Hamiltonian system. Let
$p_1, p_2, p_3\in \R^2$ be the momentum vectors of the particles respectively.
The Hamiltonian system corresponding to (\ref{1.1}) is
\be \dot{p}_i=-\frac{\partial H}{\partial q_i},\,\,\dot{q}_i
  = \frac{\partial H}{\partial p_i},\qquad {\rm for}\quad i=1,2,3,  \lb{1.2}\ee
with Hamiltonian function
\be H(p,q)=H(p_1,p_2,p_3, q_1,q_2,q_3)=\sum_{i=1}^3\frac{\|p_i\|^2}{2m_i}-U(q_1,q_2,q_3).  \lb{1.3}\ee

In 1772, Lagrange (\cite{Lag}) discovered some celebrated homographic periodic solutions, now named after
him, to the planar three-body problem, namely the three bodies form an equilateral triangle at any instant
of the motion and at the same time each body travels along a specific Keplerian elliptic orbit about the
center of masses of the system.

When $0\le e<1$, the Keplerian orbit is elliptic, following Meyer and Schmidt (\cite{MS}), we call such
elliptic Lagrangian solutions {\it elliptic relative equilibria}. Specially when $e=0$, the Keplerian
elliptic motion becomes circular motion and then all the three bodies move around the center of masses
along circular orbits with the same frequency, which are called {\it relative equilibria} traditionally.

Our main concern in this paper is the linear stability of these homographic solutions.
For the planar three-body problem with masses $m_1, m_2, m_3>0$, it turns out that the stability
of elliptic Lagrangian solutions depends on two parameters, namely the mass parameter $\beta\in [0,9]$
defined below and the eccentricity $e\in [0,1)$,
\be  \beta=\frac{27(m_1m_2+m_1m_3+m_2m_3)}{(m_1+m_2+m_3)^2}.   \lb{1.4}\ee
Note that besides local perturbation method or blow up technique which study only the case for small
enough $e>0$ or $e<1$ sufficiently close to $1$, we are not aware of any rigorous analytical method
dealing with this problem for the major part of the full range of the $(\bb,e)$ rectangle
$[0,9]\times [0,1)$, except the recent paper \cite{HS2} of the first and the third named authors.
Continuing with \cite{HS1} and \cite{HS2}, the current paper is devoted to introduce a new rigorous
analytical method to study the linear stability of the elliptic Lagrangian solutions in the full
range of the $(\bb,e)$ rectangle $[0,9]\times [0,1)$ via the index theory of symplectic paths and
the perturbation theory of linear operators.

The linear stability of relative equilibria was known more than a century ago and it is due
to Gascheau (\cite{Ga}, 1843) and Routh (\cite{R2}, 1875) independently. In this case, using
the Floquet theory one can work out all the details explicitly by hands.

After initial considerations of Danby (\cite{Dan}, 1964), Roberts (\cite{R1}, 2002) reduced all
the symmetries of the problem and their first integrals and studied the case of sufficiently
small $e\ge 0$ by perturbation techniques. He also got a partial stability bifurcation diagram
in this case, where the stability patterns are clearly presented.

In 2005, Meyer and Schmidt (cf. \cite{MS}) used heavily the central configuration nature
of the elliptic Lagrangian orbits and decomposed the fundamental solution of the elliptic Lagrangian
orbit into two parts symplectically, one of which is the same as that of the Keplerian solution and
the other is the essential part for the stability. In the current paper, the fundamental
solution of the linearized Hamiltonian system of the essential part of the elliptic Lagrangian orbit
is denoted by $\gamma_{\beta,e}(t)$ for $t\in [0,2\pi]$, which is a path of $4\times 4$ symplectic
matrices starting from the identity. They also did the stability analysis by normal form theory for
small enough $e\ge 0$.

In 2004-2006, Mart\'{\i}nez, Sam\`{a} and Sim\'{o} (\cite{MSS},\cite{MSS1},\cite{MSS2})
studied the stability problem when $e>0$ is small enough by using normal form theory, and $e<1$
and close to $1$ enough by using blow-up technique in general homogeneous potential. They
further gave a much more complete bifurcation diagram numerically and a beautiful figure was
drawn there for the full $(\bb,e)$ range, which we repeat here as Figure 1. It is one of our
primary motivations to understand this diagram globally and analytically in the present work.

\begin{figure}
\begin{center}
\resizebox{11cm}{5cm}{\psfrag{A}[][][1.5]{$\frac34$}\includegraphics*[0cm,0cm][23cm,14cm]{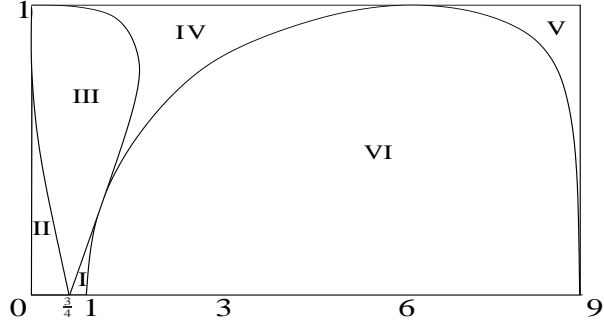}}
\caption{Stability bifurcation diagram of Lagrangian equilateral triangular
homographic orbits in the $(\beta,e)$ rectangle $[0,9]\times [0,1)$.
It is the Fig. 5 in \cite{MSS1}. Here the regions I, II, III, IV, V and VI
are EE, EE, EH, HH, HH and CS respectively.}
\end{center}
\end{figure}
\vspace{2mm}


Let $\U$ denote the unit circle in the complex plane $\C$. As in \cite{MSS1}, the following notations for
the different parameter regions are used in Figure 1:
\begin{itemize}
\item[$\bullet$] EE (elliptic-elliptic), if $\ga_{\bb,e}(2\pi)$ possesses two pairs of eigenvalues in $\U\bs\R$;

\item[$\bullet$] EH (elliptic-hyperbolic), if $\ga_{\bb,e}(2\pi)$ possesses a pair of eigenvalues in $\U\bs\R$ and
a pair of eigenvalues in $\R\bs\{0,\pm 1\}$;

\item[$\bullet$] HH (hyperbolic-hyperbolic), if $\sg(\ga_{\bb,e}(2\pi))\;\subset\;\R\bs\{0,\pm 1\}$;

\item[$\bullet$] CS (complex-saddle), if $\sg(\ga_{\bb,e}(2\pi))\;\subset\;\C\bs (\U\cup\R)$.
\end{itemize}

In summary, after these authors, the following results are rigorously proved:

(i) the relative equilibrium, i.e., the case of $e=0$, is linearly stable if
and only if $\beta<1$; it is only spectrally stable and not linearly stable when $\beta=1$,
and linearly unstable (in fact, CS) when $\beta>1$. More precisely, when $e=0$
and $\beta$ goes from $0$ to $1$, the two pairs of elliptic characteristics $\om_1$, $\ol{\om}_1$, and
$\om_2$, $\ol{\om}_2$, starting from two pairs of characteristics $+1$, and without loss of generality,
we can assume that both $\om_1$ and $\om_2$ move
clock-wisely around the unit circle with different speeds. One, say $\om_1$, moves faster and arrives
at $-1$ of the unit circle when $\beta=3/4$ and then continues to move forward around the unit circle.
At the same time, $\om_2$ moves slower along the unit circle. Then $\om_1$ and $\ol{\om}_2$ as well as $\om_2$
and $\ol{\om}_1$ collide respectively on somewhere which is not $\pm 1$ in the unit circle when $\beta=1$,
and then become CS when $\beta>1$. When $\beta=9$, they become a pair of positive double eigenvalues. So
when $e=0$ the only possible bifurcation points in the $(\beta,e)$ plane are $(3/4,0)$ and $(1,0)$.
We refer readers to \S 3.3 and Figure 3 for more detailed discussions.

(ii) It turns out that if $\beta=3/4$ is fixed and $e$ increases from $0$ slightly, the pair of $-1$
characteristics switches to a real hyperbolic pair, and two period-doubling bifurcation curves born out
from the point $(\beta, e)=(3/4, 0)$. When $\beta=1$, if $e$ increases from $0$ slightly, the two pairs of
corresponding elliptic characteristics collide and a Krein collision curve bifurcates from $(\beta, e)=(1,0)$
for such small enough $e>0$. An interesting phenomenon occurs here, namely, when $\beta$ is slightly larger
than $1$, some of the elliptic relative equilibrium with $e>0$ small can be linearly stable even though the
relative equilibrium with $e=0$ is not.

(iii) When $e<1$ and $e$ is close to $1$ enough, the relative equilibria are all HH under some non-degenerate
conditions, which is not satisfied at $\beta=6$ by numerical computations.

(iv) But the major part of the intermediate region in the $(\beta,e)$ rectangle $[0,9]\times [0,1)$ is totally
unknown theoretically, besides numerical results.

Inspired by the second named author's works on the index iteration theory of periodic orbits of Hamiltonian systems
(cf. \cite{Lon4}), the first and the third named authors initiated the program of applying the ideas of index theory and
its iteration theory in calculus of variations (cf. \cite{Lon4}) to the stability problem of periodic orbits in
celestial mechanics, especially the elliptical Lagrangian solutions (\cite{HS2}) as well as the celebrated Figure-eight
periodic orbits due to Chenciner and Montgomery in the planar three-body problem (\cite{HS1}). In \cite{HS2}, the
stability of elliptical Lagrangian solutions is studied and related to the Morse indies of their iterations,
i.e., Theorems 2.4 and 2.5 below. But in \cite{HS2}, specially the $1$-non-degeneracy of elliptic Lagrangian solutions is not
proved, and the separation curves of different index regions and thus the stability regions in $[0,9]\times [0,1)$ are
not identified.

In the current paper, we develop a new method using the $\om$-index theory for symplectic paths introduced by C. Conley,
E. Zehnder and Y. Long when $\om=1$ (cf. \cite{Lon4}) and by Y. Long when $\om\in \U\bs\{1\}$ in \cite{Lon2} and
linear differential operator theory to understand the linear stability diagram of elliptic Lagrangian solutions
theoretically in the full range of $(\bb,e)$. Specially the main purpose here is to relate such a linear stability
directly to the two major parameters of the motion: the mass parameter $\beta$ and the eccentricity $e$. For each
fixed $e\in [0,1)$, we prove that the $\om$-index $i_{\om}(\ga_{\bb,e})$ of the essential part of the elliptic Lagrangian
solutions is non-decreasing in $\bb$ for all $\om\in\U$. Then we use this important property to prove all these solutions
are $1$-non-degenerate, find the two $-1$ degeneracy curves and right envelop curve of all $\om$-degeneracy curves for
$\om\in\U\bs\{1\}$, and determine the linear stability of all sub-regions separated by these three curves. This
establishes rigorously most parts of the linear stability properties observed numerically in Figure 1, and find more
interesting properties. Note that the symplectic coordinate decomposition of Meyer and Schmidt fits quite well with
index theory, and our study will concentrate on the fundamental solution $\gamma_{\beta,e}(t)$ of the linearized
Hamiltonian system of the essential part of the elliptic Lagrangian orbit for $(\bb,e)\in [0,9]\times [0,1)$.

Denote by $\Sp(2n)$ the symplectic group of real $2n\times 2n$ matrices. For any
$\omega\in\U$ and $M\in\Sp(2n)$, let $\nu_{\om}(M)=\dim_{\C}\ker_{\C}(M-\om I_{2n})$, and $M$ is called
$\omega$-{\it degenerate} ($\om$-{\it non-degenerate} respectively) if $\nu_{\om}(M)> 0$
($\nu_{\om}(M)=0$ respectively). When $\om=1$ and if there is no confusion, we shall simply omit the
subindex $1$ and say just {\it degenerate} or {\it non-degenerate}. Let $e(M)$ be the total algebraic
multiplicity of all eigenvalues of $M$ on $\U$. We call $M\in\Sp(2n)$ {\it spectrally stable} if
$e(M)=2n$, and {\it linearly stable} if $M$ is spectrally stable and semi-simple. A symplectic matrix
$M$ is called {\it strongly linearly stable} if there is some $\ep>0$ such that all symplectic matrices
$N$ satisfying $\|M-N\|<\ep$ are linearly stable. And $M$ is {\it hyperbolic}, if $e(M)=0$.

The following is the first part of our main results in this paper.

\begin{theorem}\label{T1.1}
In the planar three-body problem with masses $m_1, m_2$, and $m_3>0$, for the elliptic Lagrangian
solution $q=(q_1(t), q_2(t), q_3(t))$ with eccentricity $e$ and mass parameter $\beta$ as given
in (\ref{1.4}), the essential part $\ga_{\bb,e}(2\pi)\in\Sp(4)$ of the monodromy matrix of the
fundamental solution along this orbit is non-degenerate for all $(\bb,e)\in (0,9]\times [0,1)$;
and when $\beta=0$, it is degenerate. Note that the Maslov-type index satisfies $i_1(\ga_{\bb,e})=0$
for all $(\bb,e)\in [0,9]\times [0,1)$.
\end{theorem}

In the proof of this theorem, we consider the second order differential operators $A(\bb,e)$
(see (\ref{2.25})) corresponding to the linear variation equation to the essential part
$\ga_{\bb,e}(t)$ of its fundamental solution along the orbit. The main ingredient of the proof
is the non-decreasing property of $\om$-index proved in Lemma \ref{L4.4} and Corollary \ref{C4.5} below
for all $\om\in\U$, by which we further prove that the operator $A(\bb,e)$ is positive definite, and
thus $1$-non-degenerate.

The rest part of this paper, specially Theorems \ref{T1.2} and \ref{T1.8} below, is devoted to rigorous
analytical studies on the existences and properties of three distinct curves $\Ga_s$, $\Ga_m$ and
$\Ga_k$ locating from left to right in the parameter $(\bb,e)$ rectangle $[0,9)\times [0,1)$. We prove
that the linear stability of the essential part $\ga_{\bb,e}(2\pi)$ of the monodromy matrix and thus
that of the elliptic Lagrange solution change precisely when $(\bb,e)$ passes through each of these
three curves, which yields a complete and rigorous understanding of the linear stability of the elliptic
Lagrange solutions. Note that here $\ga_{\bb,e}(2\pi)$ is always linearly unstable on its hyperbolic
subregion in the $(\bb,e)$ rectangle $[0,9]\times [0,1)$, and our Theorems \ref{T1.2} and \ref{T1.8}
do not distinguish the regions IV, V, and VI in Figure 1.

The main idea in the proofs of Theorems \ref{T1.2} and \ref{T1.8} is the following: By Theorem \ref{T1.1},
when $(\bb,e)$ changes, eigenvalues of $\ga_{\bb,e}(2\pi)$ can leave from the unit circle $\U$ only at
$-1$ or some Krein collision eigenvalues in $\U\bs\{\pm 1\}$. Thus such $-1$ and Krein collision eigenvalues
should correspond to $(\bb,e)$ points which form the above mentioned three curves $\Ga_s$, $\Ga_m$ and $\Ga_k$.
In order to find those $(\bb,e)$ such that $-1\in\sg(\ga_{\bb,e}(2\pi))$, we prove that the $-1$ index
$i_{-1}(\ga_{\bb,e})$ is non-increasing in $\bb\in [0,9]$ for fixed $e\in [0,1)$, and takes values $2$ at
$\bb=0$ and $0$ at $\bb=9$, thus there must exist two $-1$ index strictly decreasing curves $\Ga_s$ and
$\Ga_m$, each of which intersects every horizontal line $e={\rm constant}$ only once for $e\in [0,1)$,
and which then yield precisely the two $-1$ degeneracy curves. Next we prove that the hyperbolic region
of $\ga_{\bb,e}(2\pi)$ in the $(\bb,e)$ rectangle $[0,9]\times [0,1)$ is connected, and its boundary curve
$\Ga_k$ is continuous and thus well defined, which is the third curve for determining the linear stability.
Here the part of $\Ga_k$ which is different from the curve $\Ga_m$ is also the curve of Krein collision
eigenvalues of $\ga_{\bb,e}(2\pi)$. We prove also that the two curves $\Ga_s$ and $\Ga_m$ come from
two real analytic curves and bifurcate out from $(3/4,0)$, the curve $\Ga_k$ starts from $(1,0)$, and all
of them goes up and tends to the point $(0,1)$ as $e$ increases from $0$ to $1$. These three curves were
observed numerically in \cite{MSS1} as shown in the above Figure 1.

In this paper for any $M$ and $N\in \Sp(2n)$, we write $M\approx N$ if $M=P^{-1}NP$ holds for some
$P\in\Sp(2n)$, i.e., $N$ can be obtained from $M$ by a symplectic coordinate change. Recall that as defined
in Chapter 1 of \cite{Lon4}, the normal form of an $M\in\Sp(2n)$ is the simplest matrix $N\in \Sp(2n)$
satisfying $N\approx M$ (cf. Theorem 1.7.3 on p.36 of \cite{Lon4}). Recall also that as introduced in
Definition 1.8.9 and Theorem 1.8.10 on p.41 of \cite{Lon4} (cf. Definition \ref{D2.1} below), the basic
normal form of an $M\in\Sp(2n)$ is the simplest matrix $N\in\Sp(2n)$ such that
$\dim_{\C}\ker_{\C}(N-\om I) = \dim_{\C}\ker_{\C}(M-\om I)$ for all $\om\in\U$. It yields the homotopically
simplest form of $M$ based on the normal form of $M$ for eigenvalues in $\U$. Note that studies at the level
of basic normal forms of $\ga_{\bb,e}(2\pi)$ are easier and already powerful enough for determining the
linear stability. But the results at the level of normal forms of $\ga_{\bb,e}(2\pi)$ are stronger than
basic normal forms and involve more demonstrations. Here we describe our main results in normal forms
in Theorems \ref{T1.2} and \ref{T1.8} below. Note that here the symplectic direct sum $\dm$ is given
in (\ref{2.1}) and the normal form matrices $D(\lambda)$, $R(\theta)$, $N_1(\lm,a)$, $N_2(\om,b)$ and
$M_2(\lm,c)$ used in Theorems \ref{T1.2} and \ref{T1.8} can be found in \S 2.1 below.

\begin{theorem}\label{T1.2}
Using notations in the last theorem, for every $e\in [0,1)$, the $-1$ index $i_{-1}(\ga_{\bb,e})$
is non-increasing, and strictly decreasing only on two values of $\bb=\beta_1(e)$ and
$\bb=\beta_2(e)\in (0,9)$. Define $\Ga_i = \{(\bb_i(e),e)\;|\;e\in [0,1)\}$ for $i=1$ and $2$,
$$ \beta_s(e)=\min\{\beta_1(e),\beta_2(e)\}\quad {\it and}\quad \beta_m(e)=\max\{\beta_1(e),\beta_2(e)\}
            \qquad {\it for}\;\;e\in [0,1),  $$
and
$$ \Ga_s=\{(\bb_s(e),e)\;|\;e\in [0,1)\} \quad {\it and}\quad \Ga_m=\{(\bb_m(e),e)\;|\;e\in [0,1)\}.  $$

For every $e\in [0,1)$ we define
\be  \bb_k(e) = \sup\{\bb'\in [0,9]\;|\;\sg(\ga_{\bb,e}(2\pi))\cap \U \not= \emptyset,
         \;\;\forall\;\bb\in [0,\bb']\}, \lb{1.5}\ee
and
\be  \Ga_k=\{(\bb_k(e),e)\in [0,9]\times [0,1)\;|\;e\in [0,1)\}.   \lb{1.6}\ee

Then $\Ga_s$, $\Ga_m$ and $\Ga_k$ form three curves which possess the following properties.

(i) $0<\beta_i(e)<9$ for $i=1, 2$, and both $\beta=\beta_1(e)$ and $\beta=\beta_2(e)$ are real
analytic in $e\in [0,1)$;

(ii) $\beta_1(0)=\beta_2(0)=3/4$ and $\lim_{e\to 1}\beta_1(e)=\lim_{e\to 1}\beta_2(e)=0$. The two
curves $\Ga_1$ and $\Ga_2$ are real analytic in $e$, and bifurcate out from $(3/4,0)$ with tangents
$-\sqrt{33}/4$ and $\sqrt{33}/4$ respectively, thus they are different and their intersection
points must be isolated if there exist when $e\in (0,1)$; Consequently, $\Ga_s$ and $\Ga_m$ are
different piecewise real analytic curves;

(iii) We have
\be  i_{-1}(\ga_{\bb,e}) = \left\{\matrix{2, & {\it if}\;\;0\le \bb < \bb_s(e), \cr
                                         1, & {\it if}\;\;\bb_s(e) \le \bb < \bb_m(e), \cr
                                         0, & {\it if}\;\;\bb_m(e) \le \bb \le 9, \cr}\right. \lb{1.7}\ee
and $\Ga_s$ and $\Ga_m$ are precisely the $-1$ degeneracy curves of the matrix $\ga_{\bb,e}(2\pi)$ in
the $(\bb,e)$ rectangle $[0,9]\times [0,1)$;

(iv) There holds $\bb_s(e) \le \bb_m(e)\le \bb_k(e) < 9$ for all $e\in [0,1)$;

(v) Every matrix $\ga_{\bb,e}(2\pi)$ is hyperbolic when $\bb \in (\bb_k(e),9]$ and $e\in [0,1)$, and
there holds
\be  \bb_k(e) = \inf\{\bb\in [0,9]\;|\;\sg(\ga_{\bb,e}(2\pi))\cap \U = \emptyset\}, \qquad
          \forall\;e\in [0,1). \lb{1.8}\ee
Consequently $\Ga_k$ is the boundary curve of the hyperbolic region of $\ga_{\bb,e}(2\pi)$ in
the $(\bb,e)$ rectangle $[0,9]\times [0,1)$;

(vi) $\Ga_k$ is continuous in $e\in [0,1)$;

(vii) $\lim_{e\to 1}\bb_k(e) = 0$;

(viii) There exists a point $\tilde{e}\in (0,1]$ such that $\bb_m(e) < \bb_k(e)$
holds for all $e\in [0,\tilde{e})$. Therefore the curve $\Ga_k$ is different from the curve $\Ga_m$
at least when $e\in [0,\tilde{e})$.

(ix) We have $\ga_{\bb,e}(2\pi)\approx R(\th_1)\dm R(\th_2)$ for some $\th_1$ and
$\th_2\in (\pi, 2\pi)$, and thus it is strongly linearly stable on the segment $0<\beta<\beta_s(e)$;

(x) We have $\ga_{\bb,e}(2\pi)\approx D(\lm)\dm R(\th))$ for some $0>\lm\not= -1$ and $\th\in (\pi,2\pi)$,
and it is elliptic-hyperbolic, and thus linearly unstable on the segment $\beta_s(e)<\beta<\beta_m(e)$.

(xi) We have $\ga_{\bb,e}(2\pi)\approx R(\th_1)\dm R(\th_2)$ for some $\th_1\in (0,\pi)$ and
$\th_2\in (\pi,2\pi)$ with $2\pi-\th_2<\th_1$, and thus it is strongly linearly stable on the segment
$\beta_m(e)<\beta<\beta_k(e)$.

Here and below in Theorem \ref{T1.8}, we write $\lm=\lm_{\bb,e}$ and $\th=\th_{\bb,e}$ for short,
all of which depend on the parameters $\beta$ and $e$.
\end{theorem}

By the Bott-type iteration formulas of Maslov-type indices, we can decompose
$W^{2,2}([0,2\pi],\R^2)$ into two subspaces $E_1$ and $E_2$ (see (\ref{7.3}) and (\ref{7.4})
below) according to the boundary conditions. Then using the operator $A(\bb,e)$ (see (\ref{7.9}))
corresponding to the variational equation, we carry out the computations of Morse indices of
$A(\bb,e)|_{E_i}$ with $i=1$ and $2$ via those of $A(0,0)|_{E_i}$ with $i=1$ and $2$.

As a corollary, we have immediately

\begin{corollary}\label{C1.3}
For every $e\in [0,1)$, the Lagrangian orbit is strongly linearly stable if $\beta>0$ is small enough.
\end{corollary}

Furthermore, we can strengthen the conclusion (v) of Theorem \ref{T1.2} to

\begin{proposition}\label{P1.4}
For the equal mass case, i.e., $\beta=9$, the matrix $\ga_{9,e}(2\pi)$ is always hyperbolic and possesses
a pair of positive double eigenvalues $\lambda(e)$ and $\lm(e)^{-1}\neq 1$ for every $0\leq e<1$.
Consequently, the matrix $\gamma_{\bb,e}(2\pi)$ is hyperbolic whenever $\beta<9$ is sufficiently
close to $9$.
\end{proposition}

We establish Proposition \ref{P1.4} in \S 4.1 by the Maslov-type index theory and the theory of linear
differential operators. Then we further have the following

\begin{theorem}\label{T1.5}
(i) For every $\om\in\U\bs\{1\}$ and $e\in [0,1)$, the $\omega$-index $i_{\om}(\ga_{\bb,e})$
is decreasing for $\beta\in [0,9]$.

(ii) There exist precisely two curves in the $(\bb,e)$ rectangle $[0,9]\times [0,1)$, on which the
Maslov-type index $i_{\om}(\ga_{\bb,e})$ decreases strictly. These two curves are given by
$\bb=\beta_1(e,\omega)$ and $\bb=\beta_2(e,\omega)$ for $0\le e< 1$ respectively, where both
$\beta_1(e,\omega)$ and $\beta_2(e,\omega)$ are real analytic functions in $e\in [0,1)$ and satisfy
$\lim_{e\to 1}\beta_1(e,\omega)=\lim_{e\to 1}\beta_2(e,\omega)=0$.
\end{theorem}

This is proved in Theorem \ref{T6.3} below. Similar to the idea in the proof of Theorem \ref{T1.2},
we know that the $\omega$-Morse indices of $A(0,e)$ and $A(9,e)$ are $2$ and $0$ respectively. The
existence of the two $\om$-index strictly decreasing curves follows from the monotonicity of the
 operators involved. Note that $\om$-index strictly decreasing is equivalent to the $\om$-degeneracy
of the operator $A(\bb,e)$ by our Proposition \ref{P6.1} below. With the aid of Dunford-Taylor integral,
the $\omega$-degeneracy of $A(\bb,e)$ is related to the spectral problem of another compact operator
$B(e,\omega)$ (see (\ref{6.5}) below), and the real analyticity of the two index degeneracy curves
follows from the theory of operators.

Although $e<0$ does not have physical meaning, we can extend the fundamental solution to the case
$e\in (-1,1)$ mathematically and some interesting property of the two degeneracy curves follows.
Namely we get

\begin{theorem}\label{T1.6}
(i) The identity $\beta_1(e, -1)=\beta_2(-e, -1)$ holds for all $e\in (-1,1)$. For fixed
$\omega\in \U\bs\{-1\}$ and $i=1$ and $2$, the function $\bb_i(e,\om)$ is also even in
$e\in (-1,1)$ .

(ii) For $e\in (-1,1)$, the function $\bb_k(e)$ is also even in $e$, i.e., $\bb_k(-e)=\bb_k(e)$
for all $e\in (-1,1)$. Consequently $\Ga_k$ can be continuously extended to the region
$[0,9]\times (-1,1)$ as a curve symmetric to the segment $[0,9]\times \{0\}$.
\end{theorem}

Theorem \ref{T1.6} follows from the fact that $A(\bb, e)$ is conjugate to $A(\bb, -e)$ by a
unitary operator. This is proved in Theorems \ref{T6.4} and \ref{T7.2} below.

\begin{figure}
\begin{center}
\resizebox{6cm}{6cm}{\psfrag{Gas}[][][2.0]{$\Gamma_s$}  \psfrag{Gam}[][][2.0]{$\Gamma_m$}
\psfrag{Gak}[][][2.0]{$\Gamma_k$}\psfrag{3/4}[][][1.8]{$\frac{3}{4}$}\includegraphics*[0cm,0cm][20cm,30cm]{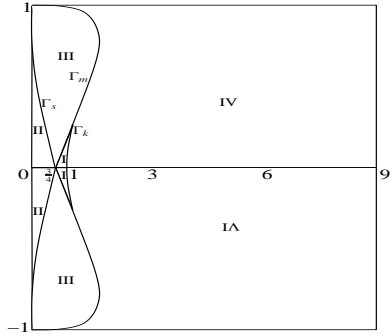}}
\caption{Stability bifurcation diagram of Lagrangian triangular homographic orbits in the
$(\beta,e)$ rectangle $[0,9]\times (-1,1)$. It is symmetric with respect ro the $\beta$-axis. The region IV is hyperbolic (HH or CS).}
\end{center}
\end{figure}
\vspace{2mm}

Figure 2 represents the case of $\omega=-1$ in Theorem \ref{T1.6}. Following Theorem \ref{T1.2},
in Figure 1, the curve separating the regions II and III is $\Ga_s$ and the curve separating the
regions III and the union of I and IV is $\Ga_m$. The curves $\Ga_1$ and the mirror of the curve
$\Ga_2$ in Theorem \ref{T1.2} together give one of the analytic curves in Theorem \ref{T1.6}, and the
another one in the Theorem \ref{T1.6} is derived from the curves $\Ga_2$ and the mirror of $\Ga_1$ in
Theorem \ref{T1.2} as indicated in the Figure 2. So we see that the two seemingly unrelated index
degeneracy curves in Theorem \ref{T1.2} are in fact coming from one degeneracy curve in Theorem
\ref{T1.6}. Note that the curve $\Ga_k$ separates the regions I, IV and VI in Figure 1. Numerical
studies on $\omega\in \U$ are also given in Figure 4 in \S 10.

When $e=1$, the operator $A(\beta,1)$ is singular and its domain is different from that of $A(\beta,e)$
with $e< 1$, and not convenient to be used. Thus we use the corresponding sesquilinear forms to
study the limiting case when $e\to 1$.

\begin{theorem}\label{T1.7}
For any fixed $0<\beta\leq 9$, the matrix $\ga_{\bb, e}(2\pi)$ is hyperbolic when $1-|e|$ is small enough.
\end{theorem}

Note that by definition, at least one pair of eigenvalues of the matrix $\ga_{\bb,e}(2\pi)$
locates on the unit circle $\U$ when $\bb\in [0,\bb_k(e)]$ and $e\in [0,1)$.

For $(\bb,e)$ located on these three special curves, we have the following

\begin{theorem}\label{T1.8}
For the normal forms of $\ga_{\bb,e}(2\pi)$ when $(\bb,e)\in \Ga_s$, $\Ga_m$ or $\Ga_k$, we
have the following results.

(i) If $\bb_s(e)<\bb_m(e)$, we have $\ga_{\bb_s(e),e}(2\pi)\approx N_1(-1,1)\dm R(\th)$ for some
$\th\in (\pi,2\pi)$, and thus it is spectrally stable and linearly unstable;

(ii) If $\bb_s(e)=\bb_m(e)<\bb_k(e)$, we have $\ga_{\bb_s(e),e}(2\pi)\approx -I_2 \dm R(\th)$ for
some $\th\in (\pi,2\pi)$, and thus it is linearly stable, but not strongly linearly stable;

(iii) If $\bb_s(e)<\bb_m(e)<\bb_k(e)$, we have $\ga_{\bb_m(e),e}(2\pi)\approx N_1(-1,-1)\dm R(\th)$
for some $\th\in (\pi,2\pi)$, and thus it is spectrally stable and linearly unstable;

(iv) If $\bb_s(e)\le\bb_m(e)<\bb_k(e)$, we have $\ga_{\bb_k(e),e}(2\pi)\approx N_2(e^{\sqrt{-1}\th},b)$
for some $\th\in (0,\pi)$ and $b=\left(\matrix{b_1 & b_2\cr
                                               b_3 & b_4\cr}\right)$
satisfying $(b_2-b_3)\sin\th>0$, i.e., $N_2(e^{\sqrt{-1}\th},b)$ is trivial in the sense of Definition
1.8.11 on p.41 of \cite{Lon4} (cf. \S 2.1 below). Consequently the matrix $\ga_{\bb_k(e),e}(2\pi)$ is
spectrally stable and linearly unstable;

(v) If $\bb_s(e)<\bb_m(e)=\bb_k(e)$, we have either $\ga_{\bb_k(e),e}(2\pi)\approx N_1(-1,1)\dm D(\lm)$
for some $-1\not=\lm < 0$ and is linearly unstable; or $\ga_{\bb_k(e),e}(2\pi)\approx M_2(-1,c)$ with
$c_1, c_2\in\R$ and $c_2\not=0$, and it is spectrally stable and linearly unstable;

(vi) If $\bb_s(e)=\bb_m(e)=\bb_k(e)$, either $\ga_{\bb_k(e),e}(2\pi)\approx M_2(-1,c)$ with $c_1, c_2\in\R$
and $c_2=0$ which possesses basic normal form $N_1(-1,1)\dm N_1(-1,1)$, or
$\ga_{\bb_k(e),e}(2\pi)\approx N_1(-1,1)\dm N_1(-1,1)$. Thus $\ga_{\bb_k(e),e}(2\pi)$ is spectrally
stable and linearly unstable.
\end{theorem}

\begin{theorem}\label{T1.9}
For any fixed $e\in[0,1)$, the set
$$  I_e=\{\beta\in (0,9]\,\left|\, {\it the\,spectrum\,of}\,\ga_{\bb,e}(2\pi)\,{\it is\,CS}\}\right.  $$
is a non-empty open set.
\end{theorem}

\begin{remark}\label{R1.10}
Note that our above results yield that the two curves $\Ga_s$ and $\Ga_m$ can have only isolated
intersection points, but it is not clear if they do have any when $e\in (0,1)$. It is not clear so
far whether $\td{e}<1$ in Theorem \ref{T1.2} and whether $\Ga_m$ and $\Ga_k$ coincide completely when
$e\in (\td{e},1)$. It is also not clear whether the $(\bb,e)$ sub-region in which $\sg(\ga_{\bb,e}(2\pi))$
is CS is connected or not.
\end{remark}

This paper is organized as follows. In \S 2, we give the definitions of $\omega$-Maslov-type index
to fix notations and its relation to the $\omega$-Morse index. Some basic variational facts on the elliptic
Lagrangian solutions are also recalled. In \S 3, we study the linear stability along the three
boundary segments of the $(\bb,e)$ rectangle $[0,9]\times [0,1)$. In \S 4, we prove the hyperbolicity of
the elliptic Lagrangian solutions in the case of equal masses (Proposition \ref{P1.4}) and the non-degeneracy
stated in Theorem \ref{T1.1}. In \S 5, the stability behavior in the limit case $e\rightarrow 1$ is
considered by the sesquilinear forms of linear operators, and Theorem \ref{T1.7} is proved.  In \S 6, we
investigate the $\omega$ degeneracy curves for general $\om\in\U\backslash\{1\}$ in the unit circle and
establish Theorem \ref{T1.5}. In \S 7, we concentrate on the $-1$ degeneracy curves.
In \S 8, we study the non-hyperbolic region and prove Theorem \ref{T1.6} and the first half of Theorem
\ref{T1.2} including its items (i)-(iii) and (ix)-(x). Section 9 is on the hyperbolic region, and we prove
the second half of Theorem \ref{T1.2} including its items (iv)-(viii) and (xi), as well as Theorems \ref{T1.8}
and \ref{T1.9}. Finally in the conclusion Section 10, we will give more comparisons for results of
Mart\'{\i}nez, Sam\`{a} and Sim\'{o} and our theorems as well as some possible future considerations.

\setcounter{equation}{0}

\section{Preliminaries}\label{sec:2}

\subsection{$\omega$-Maslov-type indices and $\omega$-Morse indices}

Let $(\R^{2n},\Omega)$ be the standard symplectic vector space with coordinates
$(x_1,...,x_n,y_1,...,y_n)$ and the symplectic form $\Omega=\sum_{i=1}^{n}dx_i \wedge dy_i$.
Let $J=\left(\matrix{0&-I_n\cr
                 I_n&0\cr}\right)$ be the standard symplectic matrix, where $I_n$
is the identity matrix on $\R^n$.

As usual, the symplectic group $\Sp(2n)$ is defined by
$$ \Sp(2n) = \{M\in {\rm GL}(2n,\R)\,|\,M^TJM=J\}, $$
whose topology is induced from that of $\R^{4n^2}$. For $\tau>0$
we are interested in paths in $\Sp(2n)$:
$$ \P_{\tau}(2n) = \{\ga\in C([0,\tau],\Sp(2n))\,|\,\ga(0)=I_{2n}\}, $$
which is equipped with the topology induced from that of $\Sp(2n)$.
For any $\om\in\U$ and $M\in\Sp(2n)$, the following real function was
introduced in \cite{Lon2}:
$$ D_{\om}(M) = (-1)^{n-1}\ol{\om}^n\det(M-\om I_{2n}). $$
Thus for any $\om\in\U$ the following codimension $1$ hypersurface
in $\Sp(2n)$ is defined (\cite{Lon2}):
$$ \Sp(2n)_{\om}^0 = \{M\in\Sp(2n)\,|\, D_{\om}(M)=0\}.  $$
For any $M\in \Sp(2n)_{\om}^0$, we define a co-orientation of
$\Sp(2n)_{\om}^0$ at $M$ by the positive direction
$\frac{d}{dt}Me^{t J}|_{t=0}$ of the path $Me^{t J}$ with $0\le t\le
\varepsilon$ and $\varepsilon$ being a small enough positive number. Let
\bea
\Sp(2n)_{\om}^{\ast} &=& \Sp(2n)\bs \Sp(2n)_{\om}^0,   \nn\\
\P_{\tau,\om}^{\ast}(2n) &=&
      \{\ga\in\P_{\tau}(2n)\,|\,\ga(\tau)\in\Sp(2n)_{\om}^{\ast}\}, \nn\\
\P_{\tau,\om}^0(2n) &=& \P_{\tau}(2n)\bs \P_{\tau,\om}^{\ast}(2n). \nn\eea
For any two continuous paths $\xi$ and $\eta:[0,\tau]\to\Sp(2n)$ with
$\xi(\tau)=\eta(0)$, we define their concatenation by:
$$ \eta\ast\xi(t) = \left\{\matrix{
            \xi(2t), & \quad {\rm if}\;0\le t\le \tau/2, \cr
            \eta(2t-\tau), & \quad {\rm if}\; \tau/2\le t\le \tau. \cr}\right. $$
Given any two $2m_k\times 2m_k$ matrices of square block form
$M_k=\left(\matrix{A_k&B_k\cr
                                C_k&D_k\cr}\right)$ with $k=1, 2$,
the symplectic sum of $M_1$ and $M_2$ is defined (cf. \cite{Lon2} and \cite{Lon4}) by
the following $2(m_1+m_2)\times 2(m_1+m_2)$ matrix $M_1\dm M_2$:
\be M_1\dm M_2=\left(\matrix{A_1 &   0 & B_1 &   0\cr
                             0   & A_2 &   0 & B_2\cr
                             C_1 &   0 & D_1 &   0\cr
                             0   & C_2 &   0 & D_2\cr}\right),   \lb{2.1} \ee
and $M^{\dm k}$ denotes the $k$ copy $\dm$-sum of $M$. For any two paths $\ga_j\in\P_{\tau}(2n_j)$
with $j=0$ and $1$, let $\ga_0\dm\ga_1(t)= \ga_0(t)\dm\ga_1(t)$ for all $t\in [0,\tau]$.

As in \cite{Lon4}, for $\lm\in\R\bs\{0\}$, $a\in\R$, $\th\in (0,\pi)\cup (\pi,2\pi)$,
$b=\left(\matrix{b_1 & b_2\cr
                 b_3 & b_4\cr}\right)$ with $b_i\in\R$ for $i=1, \ldots, 4$, and $c_j\in\R$
for $j=1, 2$, we denote respectively some normal forms by
\bea
&& D(\lm)=\left(\matrix{\lm & 0\cr
                         0  & \lm^{-1}\cr}\right), \qquad
   R(\th)=\left(\matrix{\cos\th & -\sin\th\cr
                        \sin\th  & \cos\th\cr}\right),  \nn\\
&& N_1(\lm, a)=\left(\matrix{\lm & a\cr
                             0   & \lm\cr}\right), \qquad
   N_2(e^{\sqrt{-1}\th},b) = \left(\matrix{R(\th) & b\cr
                                           0      & R(\th)\cr}\right),  \nn\\
&& M_2(\lm,c)=\left(\matrix{\lm &   1 &       c_1 &         0 \cr
                              0 & \lm &       c_2 & (-\lm)c_2 \cr
                              0 &   0 &  \lm^{-1} &         0 \cr
                              0 &   0 & -\lm^{-2} &  \lm^{-1} \cr}\right). \nn\eea
Here $N_2(e^{\sqrt{-1}\th},b)$ is {\bf trivial} if $(b_2-b_3)\sin\th>0$, or {\bf non-trivial}
if $(b_2-b_3)\sin\th<0$, in the sense of Definition 1.8.11 on p.41 of \cite{Lon4}. Note that
by Theorem 1.5.1 on pp.24-25 and (1.4.7)-(1.4.8) on p.18 of \cite{Lon4}, when $\lm=-1$ there hold
\bea
c_2 \not= 0 &{\rm if\;and\;only\;if}\;& \dim\ker(M_2(-1,c)+I)=1, \nn\\
c_2 = 0 &{\rm if\;and\;only\;if}\;& \dim\ker(M_2(-1,c)+I)=2. \nn\eea
Note that we have $N_1(\lm,a)\approx N_1(\lm, a/|a|)$ for $a\in\R\bs\{0\}$ by symplectic coordinate
change, because
$$ \left(\matrix{1/\sqrt{|a|} & 0\cr
                           0  & \sqrt{|a|}\cr}\right)
   \left(\matrix{\lm & a\cr
                  0  & \lm\cr}\right)
   \left(\matrix{\sqrt{|a|} & 0\cr
                           0  & 1/\sqrt{|a|}\cr}\right) = \left(\matrix{\lm & a/|a|\cr
                                                                         0  & \lm\cr}\right). $$

\begin{definition}\lb{D2.1} (\cite{Lon2}, \cite{Lon4})
For any $\om\in\U$ and $M\in \Sp(2n)$, define
\be \nu_{\om}(M)=\dim_{\C}\ker_{\C}(M - \om I_{2n}).  \lb{2.2}\ee

For every $M\in \Sp(2n)$ and $\om\in\U$, as in Definition 1.8.5 on p.38 of \cite{Lon4}, we define the
{\bf $\om$-homotopy set} $\Om_{\om}(M)$ of $M$ in $\Sp(2n)$ by
$$  \Om_{\om}(M)=\{N\in\Sp(2n)\,|\, \nu_{\om}(N)=\nu_{\om}(M)\},  $$
and the {\bf homotopy set} $\Om(M)$ of $M$ in $\Sp(2n)$ by
\bea  \Om(M)=\{N\in\Sp(2n)\,&|&\,\sg(N)\cap\U=\sg(M)\cap\U,\,{\it and}\; \nn\\
         &&\qquad \nu_{\lm}(N)=\nu_{\lm}(M)\qquad\forall\,\lm\in\sg(M)\cap\U\}.  \nn\eea
We denote by $\Om^0(M)$ (or $\Om^0_{\om}(M)$) the path connected component of $\Om(M)$ ($\Om_{\om}(M)$)
which contains $M$, and call it the {\bf homotopy component} (or $\om$-{\bf homtopy component}) of $M$ in
$\Sp(2n)$. Following Definition 5.0.1 on p.111 of \cite{Lon4}, for $\om\in \U$ and $\ga_i\in \P_{\tau}(2n)$
with $i=0, 1$, we write $\ga_0\sim_{\om}\ga_1$ if $\ga_0$ is homotopic to $\ga_1$ via
a homotopy map $h\in C([0,1]\times [0,\tau], \Sp(2n))$ such that $h(0)=\ga_0$, $h(1)=\ga_1$, $h(s)(0)=I$,
and $h(s)(\tau)\in \Om_{\om}^0(\ga_0(\tau))$ for all $s\in [0,1]$. We write also $\ga_0\sim \ga_1$, if
$h(s)(\tau)\in \Om^0(\ga_0(\tau))$ for all $s\in [0,1]$ is further satisfied.
\end{definition}

Following Definition 1.8.9 on p.41 of \cite{Lon4}, we call the above matrices $D(\lm)$, $R(\th)$, $N_1(\lm,a)$
and $N_2(\om,b)$ basic normal forms of symplectic matrices. As proved in \cite{Lon2} and \cite{Lon3} (cf.
Theorem 1.9.3 on p.46 of \cite{Lon4}), every $M\in\Sp(2n)$ has its basic normal form decomposition in $\Om^0(M)$
as a $\dm$-sum of these basic normal forms. This is very important when we derive basic normal forms
for $\ga_{\bb,e}(2\pi)$ to compute the $\om$-index $i_{\om}(\ga_{\bb,e})$ of the path $\ga_{\bb,e}$ later
in this paper.

We define a special continuous symplectic path $\xi_n\subset {\Sp}(2n)$ by
\be \xi_n(t) = \left(\matrix{2-\frac{t}{\tau} & 0 \cr
                             0 &  (2-\frac{t}{\tau})^{-1}\cr}
               \right)^{\dm n} \qquad {\rm for}\;0\le t\le \tau.  \lb{2.3}\ee

\begin{definition} (\cite{Lon2}, \cite{Lon4})\lb{D2.2}
{For any $\tau>0$ and $\ga\in \P_{\tau}(2n)$, define
\be \nu_{\om}(\ga)= \nu_{\om}(\ga(\tau)).  \lb{2.4}\ee

If $\ga\in\P_{\tau,\om}^{\ast}(2n)$, define
\be i_{\om}(\ga) = [\Sp(2n)_{\om}^0: \ga\ast\xi_n],  \lb{2.5}\ee
where the right hand side of (\ref{2.5}) is the usual homotopy intersection number, and
the orientation of $\ga\ast\xi_n$ is its positive time direction under homotopy with
fixed end points.

If $\ga\in\P_{\tau,\om}^0(2n)$, we let $\mathcal{F}(\ga)$ be the set of all open
neighborhoods of $\ga$ in $\P_{\tau}(2n)$, and define
\be i_{\om}(\ga) = \sup_{U\in\mathcal{F}(\ga)}\inf\{i_{\om}(\beta)\,|\,
                       \beta\in U\cap\P_{\tau,\om}^{\ast}(2n)\}.      \lb{2.6}\ee
Then
$$ (i_{\om}(\ga), \nu_{\om}(\ga)) \in \Z\times \{0,1,\ldots,2n\}, $$
is called the index function of $\ga$ at $\om$. }
\end{definition}

We refer to \cite{Lon4} for more details on this index theory of symplectic matrix paths
and periodic solutions of Hamiltonian system.

For $T>0$, suppose $x$ is a critical point of the functional
$$ F(x)=\int_0^TL(t,x,\dot{x})dt,  \qquad \forall\,\, x\in W^{1,2}(\R/T\Z,\R^n), $$
where $L\in C^2((\R/T\Z)\times \R^{2n},\R)$ and satisfies the
Legendrian convexity condition $L_{p,p}(t,x,p)>0$. It is well known
that $x$ satisfies the corresponding Euler-Lagrangian
equation:
\bea
&& \frac{d}{dt}L_p(t,x,\dot{x})-L_x(t,x,\dot{x})=0,    \label{2.7}\\
&& x(0)=x(T),  \qquad \dot{x}(0)=\dot{x}(T).    \label{2.8}\eea

For such an extremal loop, define
\bea
P(t) &=& L_{p,p}(t,x(t),\dot{x}(t)),  \nn\\
Q(t) &=& L_{x,p}(t,x(t),\dot{x}(t)),  \nn\\
R(t) &=& L_{x,x}(t,x(t),\dot{x}(t)).  \nn\eea
Note that
\be F\,''(x)=-\frac{d}{dt}(P\frac{d}{dt}+Q)+Q^T\frac{d}{dt}+R. \lb{2.9}\ee

For $\omega\in\U$, set
\be  D(\omega,T)=\{y\in W^{1,2}([0,T],\C^n)\,|\, y(T)=\omega y(0) \}.   \lb{2.10}\ee
We define the $\omega$-Morse index $\phi_\omega(x)$ of $x$ to be the dimension of the
largest negative definite subspace of
$$ \langle F\,''(x)y_1,y_2 \rangle, \qquad \forall\;y_1,y_2\in D(\omega,T), $$
where $\langle\cdot,\cdot\rangle$ is the inner product in $L^2$. For $\omega\in\U$, we
also set
\be  \ol{D}(\omega,T)= \{y\in W^{2,2}([0,T],\C^n)\,|\, y(T)=\omega y(0), \dot{y}(T)=\om\dot{y}(0) \}.
                     \lb{2.11}\ee
Then $F''(x)$ is a self-adjoint operator on $L^2([0,T],\R^n)$ with domain $\ol{D}(\omega,T)$.
We also define
\[\nu_\omega(x)=\dim\ker(F''(x)).\]

In general, for a self-adjoint operator $A$ on the Hilbert space $\mathscr{H}$, we set
$\nu(A)=\dim\ker(A)$ and denote by $\phi(A)$ its Morse index which is the maximum dimension
of the negative definite subspace of the symmetric form $\langle A\cdot,\cdot\rangle$. Note
that the Morse index of $A$  is equal to the total multiplicity of the negative eigenvalues
of $A$.

On the other hand, $\td{x}(t)=(\partial L/\partial\dot{x}(t),x(t))^T$ is the solution of the
corresponding Hamiltonian system of (\ref{2.7})-(\ref{2.8}), and its fundamental solution
$\gamma(t)$ is given by
\bea \dot{\gamma}(t) &=& JB(t)\gamma(t),  \lb{2.12}\\
     \gamma(0) &=& I_{2n},  \lb{2.13}\eea
with
\be B(t)=\left(\matrix{P^{-1}(t)& -P^{-1}(t)Q(t)\cr
                       -Q(t)^TP^{-1}(t)& Q(t)^TP^{-1}(t)Q(t)-R(t)\cr}\right). \lb{2.14}\ee

\begin{lemma}(Long, \cite{Lon4}, p.172)\lb{L2.3}  
For the $\omega$-Morse index $\phi_\omega(x)$ and nullity $\nu_\omega(x)$ of the solution $x=x(t)$
and the $\omega$-Maslov-type index $i_\omega(\gamma)$ and nullity $\nu_\omega(\gamma)$ of the symplectic
path $\ga$ corresponding to $\td{x}$, for any $\omega\in\U$ we have
\be \phi_\omega(x) = i_\omega(\gamma), \qquad \nu_\omega(x) = \nu_\omega(\gamma).  \lb{2.15}\ee
\end{lemma}

A generalization of the above lemma to arbitrary  boundary conditions is given in \cite{HS1}.
For more information on these topics, we refer to \cite{Lon4}.

\subsection{Stability criteria via Morse indices}

In this subsection we recall some results of \cite{HS2}.

Gordon's classical theorem (cf. \cite{Go}) says that the Keplerian solution is a minimizer in the loop
space under some topological constraint. To the essential part, a theorem of Venturelli in
\cite{V1} as well as Zhang and Zhou in \cite{ZZ1} tells us that the Lagrangian solution is
the minimizer among the loops in its some homology class. Note that, up to now, these are the
only known variational facts under topological constraints on the loop spaces.

By these theorems, we got criteria for the stability in terms of Morse indices. Let
$\phi^k=\phi_1(q^k)$ be the Morse index of the $k$-th iteration $q^k$ of the Lagrangian
solution $q$ in the variational problem, and according to \cite{V1} and \cite{ZZ1},
$\phi^1=0$ holds.

The Lagrangian solution is called linearly stable (spectrally stable) if $\gamma(2\pi)$ is
linearly stable(spectrally stable). The first and the third named authors proved the following

\begin{theorem} (Hu-Sun, \cite{HS2})\label{T2.4}
For the monodromy matrix $M$ corresponding to the elliptic Lagrangian solution
$q=(q_1(t),q_2(t),q_3(t))$, $2\leq\phi^2\leq 4$ and,
\be  e(M)/2\geq \phi^2. \lb{2.16}\ee
Moreover

(i) If $\phi^2=4$, then the Lagrangian solution is spectrally stable;

(ii) If $\phi^2=3$, then the Lagrangian solution is linearly unstable;

(iii) If $\phi^2=2$, then the Lagrangian solution is spectrally stable if there exists some
integer $k\geq 3$, such that $\phi^k>2(k-1)$.

(iv) If $\phi^k=2(k-1)$, for all $k\in \N$, then the Lagrangian solution is linearly unstable.
\end{theorem}

Moreover, if the essential part $\ga=\ga_{\bb,e}(t)$ (cf. Subsection 2.3 below) of the monodromy
matrix at $t=4\pi$ is non-degenerate, we can get the normal forms of $\ga_{\bb,e}(2\pi)$.

\begin{theorem} (Hu-Sun, \cite{HS2})\lb{T2.5}
In the same setting of the above theorem, suppose $\ga_{\bb,e}(4\pi)=\ga_{\bb,e}(2\pi)^2$ is
non-degenerate.

(i) If $\phi^2=4$, then $\ga_{\bb,e}(2\pi)\approx R(2\pi-\theta_1)\dm R(2\pi-\theta_2)$ holds for
some $\theta_1$ and $\theta_2\in (0,\pi)$, and is linearly stable;

(ii) If $\phi^2=3$, then $\ga_{\bb,e}(2\pi)\approx D(\lambda)\dm R(2\pi-\theta)$ for some $\lambda<0$
and $\theta\in(0,\pi)$, and is linearly unstable;

(iii) If $\phi^2=2$ and there exists some integer $k\geq 3$ such that $\phi^k>2(k-1)$, then
$\ga_{\bb,e}(2\pi)\approx R(2\pi-\th_1)\dm R(\th_2)$ holds with $0<\th_1<\th_2<\pi$,
and is linearly stable;

(iv) If $\phi^k=2(k-1)$ for all $k\in \N$,  then $\ga_{\bb,e}(2\pi)$ is hyperbolic or spectrally
stable and linearly unstable.
\end{theorem}

\subsection{The essential part of the fundamental solution of the elliptic Lagrangian orbit}

Following Meyer and Schmidt (cf. p.275 of \cite{MS}), the essential part $\ga=\ga_{\bb,e}(t)$ of the
fundamental solution of the Lagrangian orbit satisfies
\bea
\dot{\gamma}(t) &=& JB(t)\gamma(t),   \lb{2.17}\\
\gamma(0) &=& I_{4},    \lb{2.18}\eea
with
\be B(t)=\left(\matrix{1 & 0 & 0 & 1\cr
                       0 & 1 & -1 & 0 \cr
                       0 & -1 &\frac{2e\cos(t)-1-\sqrt{9-\beta}}{2(1+e\cos(t))} & 0 \cr
                       1 & 0 & 0 & \frac{2e\cos(t)-1+\sqrt{9-\beta}}{2(1+e\cos(t))} \cr}\right), \lb{2.19}\ee
where $e$ is the eccentricity, and $t$ is the truly anomaly.

Let
\be  J_2=\left(\matrix{ 0 & -1 \cr 1 & 0 \cr}\right), \qquad
    K_{\bb,e}(t)=\left(\matrix{\frac{3+\sqrt{9-\beta}}{2(1+e\cos(t))} & 0 \cr
                                     0 & \frac{3-\sqrt{9-\beta}}{2(1+e\cos(t))} \cr}\right),  \lb{2.20}\ee
and set
\be L(t,x,\dot{x})=\frac{1}{2}\|\dot{x}\|^2 + J_2x(t)\cdot\dot{x}(t) + \frac{1}{2}K_{\bb,e}(t)x(t)\cdot x(t),
       \qquad\quad  \forall\;x\in W^{1,2}(\R/2\pi\Z,\R^2),  \lb{2.21}\ee
where $a\cdot b$ denotes the inner product in $\R^2$. Obviously the origin in the configuration space is a
solution of the corresponding Euler-Lagrange system. By Legendrian transformation, the corresponding
Hamiltonian function is
$$   H(t,z)=\frac{1}{2}B(t)z\cdot z,\qquad \forall\; z\in\R^4.  $$

Note first that the elliptical Lagrangian solution is a local minimizer of the action functional $\mathcal{A}$
on the homology class of curves with winding number $(1,1,1)$ or $(-1,-1,-1)$ in $H_1(\hat{\mathcal{X}})$
(cf. \cite{V1}, \cite{ZZ1}, and Lemma 4.1 of \cite{HS2}). Note secondly that curves with winding number not
equal to $\pm 1$ can not approximate curves with winding number $1$ or $-1$ in $C(\R/\Z, \R^2)$. Therefore
the Morse index of $\mathcal{A}$ at the elliptical Lagrangian solution in the whole space $\hat{\mathcal{X}}$
and that restricted in the homology class of curves with winding number $(1,1,1)$ or $(-1,-1,-1)$ in
$H_1(\hat{\mathcal{X}})$ coincide and take the value zero.

Here note that let $\ga_1$ be the fundamental solution of the Kepler orbit. Then $\ga_1\dm\ga_{\bb,e}$
is the fundamental solution of the Lagrangian orbit. By Theorem 7.3.1 on p.168 of \cite{Lon4} and the
additivity of the index theory (cf. (ii) of Theorem 6.2.7 on p.147 of \cite{Lon4}), we obtain
\be   \phi^k = i_1(\ga_1^k)+i_1(\ga_{\bb,e}^k),  \qquad \forall\;k\in\N,  \lb{2.22}\ee
where $\phi^k$ is defined at the beginning of \S 2.2. Note that by Proposition 3.6 in p.110 of \cite{HS2},
we have
\be  i_1(\ga_1^k)=2(k-1), \qquad \forall\;k\in\N.  \lb{2.23}\ee
Thus by our above discussions and (\ref{2.22})-(\ref{2.23}) with $k=1$, we obtain
\be i_1(\ga_{\bb,e})=\phi^1=0, \qquad \forall \, (\bb,e)\in [0,9]\times [0,1). \lb{2.24}\ee

\subsection{A modification on the path $\ga_{\bb,e}(t)$}

In order to transform the Lagrangian system (\ref{2.19}) to a simpler linear operator corresponding to
a second order Hamiltonian system with the same linear stability as $\ga_{\bb,e}(2\pi)$, using $R(t)$
and $R_4(t)\equiv N_2(e^{t\sqrt{-1}},0)$ defined in \S 2.1, we let
\be  \xi_{\bb,e}(t) = R_4(t)\ga_{\bb,e}(t), \qquad \forall\; t\in [0,2\pi], (\bb,e)\in [0,9]\times [0,1). \lb{2.25}\ee
One can show by direct computation that
\be  \frac{d}{dt}\xi_{\bb,e}(t)
  = J \left(\matrix{I_2 & 0 \cr
                    0 & R(t)(I_2-K_{\bb,e}(t))R(t)^T \cr}\right)\xi_{\bb,e}(t). \lb{2.26}\ee
Note that $R_4(0)=R_4(2\pi)=I_4$, so $\ga_{\bb,e}(2\pi)=\xi_{\bb,e}(2\pi)$ holds and the linear stabilities
of the systems (\ref{2.18}) and (\ref{2.26}) are precisely the same.

By (\ref{2.25}) the symplectic paths $\ga_{\bb,e}$ and $\xi_{\bb,e}$ are homotopic to each other via the
homotopy $h(s,t)=R_4(st)\ga_{\bb,e}(t)$ for $(s,t)\in [0,1]\times [0,2\pi]$. Because $R_4(s)\ga_{\bb,e}(2\pi)$
for $s\in [0,1]$ is a loop in $\Sp(4)$ which is homotopic to the constant loop $\ga_{\bb,e}(2\pi)$, we have
$\ga_{\bb,e} \sim_1 \xi_{\bb,e}$ by the homotopy $h$. Then by Lemma 5.2.2 on p.117 of \cite{Lon4}, the
homotopy between $\ga_{\bb,e}$ and $\xi_{\bb,e}$ can be realized by a homotopy which fixes the end point
$\ga_{\bb,e}(2\pi)$ all the time. Therefore by the homotopy invariance of the Maslov-type index (cf. (i) of
Theorem 6.2.7 on p.147 of \cite{Lon4}) we obtain
\be  i_{\om}(\xi_{\bb,e}) = i_{\om}(\ga_{\bb,e}), \quad \nu_{\om}(\xi_{\bb,e}) = \nu_{\om}(\ga_{\bb,e}),
       \qquad \forall \,\omega\in\U, \; (\bb,e)\in [0,9]\times [0,1). \lb{2.27}\ee
On the other hand, the first order linear Hamiltonian system (\ref{2.26}) corresponds to the following second order
linear Hamiltonian system
\be  \ddot{x}(t)=-x(t)+R(t)K_{\bb,e}(t)R(t)^Tx(t). \lb{2.28}\ee

For $(\bb,e)\in [0,9)\times [0,1)$, the second order differential operator corresponding to (\ref{2.28}) is given by
\bea  A(\bb,e)
&=& -\frac{d^2}{dt^2}I_2-I_2+R(t)K_{\bb,e}(t)R(t)^T  \nn\\
&=& -\frac{d^2}{dt^2}I_2-I_2+\frac{1}{2(1+e\cos t)}(3I_2+\sqrt{9-\beta}S(t)),  \lb{2.29}\eea
where $S(t)=\left(\matrix{ \cos 2t & \sin 2t \cr
                           \sin 2t & -\cos 2t \cr}\right)$, defined on the domain $\ol{D}(\omega,2\pi)$
in (\ref{2.11}). Then it is self-adjoint and depends on the parameters $\bb$ and $e$. By Lemma
\ref{L2.3}, we have for any $\bb$ and $e$, the Morse index $\phi_{\om}(A(\bb,e))$ and nullity $\nu_{\om}(A(\bb,e))$
of the operator $A(\bb,e)$ on the domain $\ol{D}(\omega,2\pi)$ satisfy
\be  \phi_{\om}(A(\bb,e)) = i_{\om}(\xi_{\bb,e}), \quad \nu_{\om}(A(\bb,e)) = \nu_{\om}(\xi_{\bb,e}), \qquad
           \forall \,\om\in\U. \lb{2.30}\ee
Specially by Lemma 4.1, (55) and (58) of \cite{HS2} and the above (\ref{2.24}), we obtain
\be  i_1(\xi_{\bb,e})=\phi_1(A(\bb,e))=i_1(\ga_{\bb,e})=0, \qquad \forall \,(\bb,e)\in [0,9]\times [0,1), \lb{2.31}\ee

In the rest part of this paper, we shall use both of the paths $\ga_{\bb,e}$ and $\xi_{\bb,e}$ to study
the linear stability of $\ga_{\bb,e}(2\pi)=\xi_{\bb,e}(2\pi)$. Because of (\ref{2.27}), in many cases and
proofs below, we shall not distinguish these two paths.

\setcounter{equation}{0}
\section{Stability on the three boundary segments of the rectangle $[0,9]\times [0,1)$}
\label{sec:3}

We need more precise information on stabilities and indices of the three boundary segments of the
$(\bb,e)$ rectangle $[0,9]\times [0,1)$.

\subsection{The boundary segment $\{0\}\times [0,1)$}

When $\beta=0$, this is the case with two zero masses, and the essential part of the fundamental solution of
Lagrangian orbit is also the fundamental solution of the Keplerian orbits. In fact, when $\bb=0$, without
loss of generality, by the definition (\ref{1.4}) of $\bb$ we may assume $m_2=m_3=0$ and $m_1>0$, and then
every elliptic Lagrangian solution becomes the motion along two Keplerian solutions of the two points $q_2$ and
$q_3$ with zero masses going along their elliptic orbits around fixed point $q_1(t)\equiv 0$ with mass
$m_1>0$. When $\bb=0$, the matrix $B(t)$ in (\ref{2.19}) becomes
\be B(t)=\left(\matrix{1 & 0  & 0  & 1\cr
                       0 & 1  & -1 & 0 \cr
                       0 & -1 & -\frac{2-e\cos(t)}{1+e\cos(t)} & 0 \cr
                       1 & 0  & 0  & 1 \cr}\right), \lb{3.1}\ee
which coincides with the coefficient matrix $\td{B}(t)$ in (17) on p.275 of \cite{MS}. Note that by the
different sign choice of the standard symplectic matrices on p.259 of \cite{MS}, the order of the independent
variables there are different from that in our system (\ref{2.17})-(\ref{2.19}).

{\bf (A)} {\it The case of $e=0$.}

In this case, the matrix $B(t)$ becomes independent of $t$:
\be B\equiv B(t) = \left(\matrix{1 & 0  & 0  & 1\cr
                                 0 & 1  & -1 & 0 \cr
                                 0 & -1 & -2 & 0 \cr
                                 1 & 0  & 0  & 1 \cr}\right).    \lb{3.2}\ee
Then one can find the fundamental solution $\ga_{0,0}(t)$ of the corresponding system (\ref{2.1})
with constant coefficient $JB$ explicitly:
\be  \ga_{0,0}(t)
= \left(\matrix{2-\cos t   & 3t-2\sin t  & 3t-\sin t   & 1-\cos t   \cr
                -\sin t    & 2\cos t -1  & \cos t -1   & -\sin t    \cr
                \sin t     & 2-2\cos t   & 2-\cos t    & \sin t     \cr
                2\cos t -2 & 4\sin t -3t & 2\sin t -3t & 2\cos t -1 \cr}\right). \lb{3.3}\ee
Letting
$$ P =  \left(\matrix{1 & 0         & 0  & 6\pi  \cr
                      0 & -1/(6\pi) & -1 & 0     \cr
                      0 & 0         & 1  & 0     \cr
                      0 & 0         & 0  & -6\pi \cr}\right), $$
we then obtain
$$  \check{\ga}(t) \equiv P^{-1}\ga_{0,0}(t)P
= \left(\matrix{\cos t                 & -\frac{2\sin t}{6\pi}       & -\sin t              & 0 \cr
                0                      & 1                           & 0                    & 0 \cr
                \sin t                 & \frac{2\cos t -2}{6\pi}     & \cos t               & 0 \cr
                \frac{2-2\cos t}{6\pi} & \frac{4\sin t -3t}{36\pi^2} & \frac{2\sin t}{6\pi} & 1 \cr}\right). $$
Next for $\ep\in [0,1]$ we consider the following homotopy path $\check{\ga}_{\ep}(t)$:
$$  \check{\ga}_{\ep}(t)
= \left(\matrix{\cos t                 & -\ep\frac{2\sin t}{6\pi}       & -\sin t              & 0 \cr
                0                      & 1                           & 0                    & 0 \cr
                \sin t                 & \ep\frac{2\cos t -2}{6\pi}     & \cos t               & 0 \cr
                \ep\frac{2-2\cos t}{6\pi} & \frac{4\sin t -3t}{36\pi^2} & \ep\frac{2\sin t}{6\pi} & 1 \cr}\right). $$
Then $\check{\ga}_{\ep}(t)\in \Sp(4)$ and $\check{\ga}_{\ep}(0)=I_4$ hold for all $t\in\R$
and $\ep\in [0,1]$. We have $\check{\ga}_{1}(t)=\check{\ga}(t)$ and
$$  \check{\ga}_{0}(t)
\;=\; \left(\matrix{\cos t & 0                           & -\sin t & 0 \cr
                0      & 1                           & 0       & 0 \cr
                \sin t & 0                           & \cos t  & 0 \cr
                0      & \frac{4\sin t -3t}{36\pi^2} & 0       & 1 \cr}\right)
\;=\; R(t)\dm\left(\matrix{1 & 0 \cr
                       \frac{4\sin t -3t}{36\pi^2} & 1 \cr}\right)
\;\sim\;  R(t)\dm\left(\matrix{1 & \frac{t}{2\pi} \cr
                              0 & 1              \cr}\right),  $$
and specially we obtain
\be  \ga_{0,0}(2\pi)
= \left(\matrix{1 & 6\pi  & 6\pi  & 0 \cr
                0 & 1     & 0     & 0 \cr
                0 & 0     & 1     & 0 \cr
                0 & -6\pi & -6\pi & 1 \cr}\right)
= P\left(\matrix{1 & 0               & 0 & 0 \cr
                 0 & 1               & 0 & 0 \cr
                 0 & 0               & 1 & 0 \cr
                 0 & -\frac{1}{6\pi} & 0 & 1 \cr}\right)P^{-1}
\approx I_2\dm N_1(1,1).  \lb{3.4}\ee

{\bf (B)} {\it The case of $e\in [0,1)$.}

Note that it is well known that the matrix $\Phi(t)$ with $t=f$ in Lemma 3.1 on p.271 of \cite{MS} consists
of two parts, one of which corresponds to the Keplerian elliptic solution of the two-body problem and the other
of which corresponds to the coefficient matrix $B(t)$ of the essential part $\ga_{\bb,e}(t)$ in our notations.
Specially when $\bb=0$ (i.e., $\sg=0$ in \cite{MS}), both the two parts of the matrix $\Phi(t)$ coincide to
each other, and yields precisely $\Phi(t)=\td{B}(t)\dm \td{B}(t)$ there. Therefore when $\bb=0$, the
essential part $\ga_{\bb,e}(t)$ corresponds to the Keplerian elliptic solution of the two-body problem.

Thus by Lemma 3.3 and the discussion on (46) of \cite{HS2}, the matrix $\ga_{0,e}(2\pi)$ satisfies:
\be   \ga_{0,e}(2\pi) \approx I_2\dm N_1(1,1), \qquad \forall\,e\in [0,1),   \lb{3.5}\ee
and it includes (\ref{3.4}) as a special case.

{\bf (C)} {\it The indices $i_{\om}(\ga_{0,e})$ for $\om\in\U$.}

By (\ref{2.31}) we obtain
\be   \phi_1(A(0,e)) = i_1(\xi_{0,e}) = i_1(\ga_{0,e}) = 0.  \lb{3.6}\ee
By (\ref{3.5}) and properties of splitting numbers in Chapter 9 of \cite{Lon4}, as in (56) of
\cite{HS2}, for $\om\in\U\bs\{1\}$ and $M=\ga_{0,e}(2\pi)$ we obtain
\bea  i_{\om}(\ga_{0,e})
&=& i_{\om}(\xi_{0,e})  \nn\\
&=& i_1(\xi_{0,e}) + S_{M}^+(1) - S_{M}^-(\om)  \nn\\
&=& 0 + S_{I_2}^+(1) + S_{N_1(1,1)}^+(1) - 0  \nn\\
&=& 2.    \lb{3.7}\eea
For every $e\in [0,1)$, note that (\ref{3.5}) yields also
\be  \nu_{\om}(\ga_{0,e}) = \nu_{\om}(\xi_{0,e})
= \left\{\matrix{3, &\quad {\rm if}\;\;\om = 1, \cr
                 0, &\quad {\rm if}\;\;\om \in \U\bs\{1\}. \cr}\right.  \lb{3.8}\ee

\subsection{The boundary segment $\{9\}\times [0,1)$}

This is studied in our Proposition \ref{P1.4} and more precisely in \S 4 below.
Specially we have
$$   \sg(\ga_{9,e}(2\pi)) \subset \R_+\bs\{1\}, \qquad \forall\;e\in [0,1)   $$
and it possesses a pair of double positive real hyperbolic eigenvalues. By
(\ref{2.31}),  we have $i_1(\ga_{9,e}) = 0$ for all
$e\in [0,1)$. By our studies in \S 4.1, the matrix $\ga_{9,e}(2\pi)$ possesses always two double
positive eigenvalues not equal to $1$, thus by the definition of the splitting number in
\S 9.1 of \cite{Lon4}, we have $S_M^{\pm}(\om)=0$ for $M=\ga_{9,e}(2\pi)$ with $e\in [0,1)$
and all $\om\in\U$. Then for all $e\in [0,1)$ and $\om\in\U$ this yields
\be i_{\om}(\ga_{9,e}) = i_{\om}(\ga_{9,0}) = 0, \qquad
    \nu_{\om}(\ga_{9,e}) = \nu_{\om}(\ga_{9,0}) = 0.  \lb{3.9}\ee

\subsection{The boundary segment $[0,9]\times \{0\}$}

In this case $e=0$. It is considered in (A) of Subsection 3.1 when $\bb=0$. When $\bb\in (0,9]$, this
is the case of circular orbits with three positive masses. It was studied in Section 4 of \cite{R1}
by Roberts and in pp.275-276 of \cite{MS} by Meyer and Schmidt. Below, we shall first recall the
properties of eigenvalues of $\ga_{\bb,0}(2\pi)$. Then we carry out the computations of normal
forms of $\ga_{\bb,0}(2\pi)$, and $\pm 1$ indices $i_{\pm 1}(\ga_{\bb,0})$ of the path
$\ga_{\bb,0}$ for all $\bb\in [0,9]$, which are new.

In this case, the essential part of the motion (\ref{2.17})-(\ref{2.19}) becomes an ODE
system with constant coefficients:
\be B = B(t) = \left(\matrix{1 & 0 & 0 & 1\cr
                             0 & 1 & -1 & 0 \cr
                             0 & -1 & -\frac{\sqrt{9-\beta}+1}{2} & 0 \cr
                             1 & 0 & 0 & \frac{\sqrt{9-\beta}-1}{2} \cr}\right).  \lb{3.10}\ee
The characteristic polynomial $\det(JB-\lambda I)$ of $JB$ is given by
\be \lambda^4 + \lambda^2 + \frac{\bb}{4} = 0.  \lb{3.11}\ee
Letting $\aa=\lambda^2$, the two roots of the quadratic polynomial $\aa^2 + \aa +\frac{\bb}{4}$
are given by $\aa=\frac{1}{2}(-1\pm \sqrt{1-\bb})$. Therefore the four roots of the polynomial
(\ref{3.11}) are given by
\be
\aa_{1,\pm} = \pm\sqrt{\frac{1}{2}(-1+\sqrt{1-\bb})}, \qquad
\aa_{2,\pm} = \pm\sqrt{\frac{1}{2}(-1-\sqrt{1-\bb})}. \lb{3.12}\ee

{\bf (A)  Eigenvalues of $\ga_{\bb,0}(2\pi)$ for $\bb\in [0,9]$.}

When $0\le \bb\le 1$, by (\ref{3.12}), we get the four well-known characteristic multipliers
of the matrix $\ga_{\bb,0}(2\pi)$
\be   \rho_{i,\pm}(\beta) = e^{2\pi\aa_{i,\pm}} =  e^{\pm 2\pi\sqrt{-1}\th_i(\bb)}, \qquad {\rm for}\;\;i=1, 2,  \lb{3.13}\ee
where
\be  \th_1(\bb) = \sqrt{\frac{1}{2}(1-\sqrt{1-\bb})}, \qquad \th_2(\bb) = \sqrt{\frac{1}{2}(1+\sqrt{1-\bb})}. \lb{3.14}\ee

When $1<\bb\le 9$, from (\ref{3.12}) by direct computation the four characteristic multipliers
of the matrix $\ga_{\bb,0}(2\pi)$ are given by
\be   \rho_{\pm,\pm}  = e^{\pm\pi\sqrt{\sqrt{\bb}-1}}e^{\pm \pi\sqrt{-1}\sqrt{\sqrt{\bb}+1}}.
              \lb{3.15}\ee

Specially, we obtain the following results:

When $\bb=0$, we have $\sg(\ga_{0,0}(2\pi)) = \{1, 1, 1, 1\}$, which coincides with  (\ref{3.5}).

When $0<\bb<3/4$, in (\ref{3.14}) the angle $\th_1(\bb)$ increases strictly from $0$ to $1/2$ as $\bb$ increases
from $0$ to $3/4$. Therefore $\rho_{1,+}(\bb)=e^{2\pi \sqrt{-1}\th_1(\bb)}$ runs from $1$ to $-1$ counterclockwise along
the upper semi-unit circle in the complex plane $\C$ as $\bb$ increases from $0$ to $3/4$. Correspondingly
$\rho_{1,-}(\bb)=e^{-2\pi \sqrt{-1}\th_1(\bb)}$ runs from $1$ to $-1$ clockwise along the lower semi-unit circle in
$\C$ as $\bb$ increases from $0$ to $3/4$. At the same time, because $\th_2(\bb)$ decreases strictly from $1$
to $\sqrt{3}/2$ as $\bb$ increases from $0$ to $3/4$, therefore $\rho_{2,+}(\bb)=e^{2\pi \sqrt{-1}\th_2(\bb)}$ runs from
$1$ to $e^{\sqrt{-1}\sqrt{3}\pi}$ clockwise along the lower semi-unit circle in $\C$ as $\bb$ increases from $0$ to $3/4$.
Correspondingly $\rho_{2,-}(\bb)=e^{-2\pi \sqrt{-1}\th_2(\bb)}$ runs from $1$ to $e^{-\sqrt{-1}\sqrt{3}\pi}$ counterclockwise
along the upper semi-unit circle in $\C$ as $\bb$ increases from $0$ to $3/4$. Thus specially we obtain
$\sg(\ga_{\bb,0}(2\pi)) \subset \U\bs\R$ for all $\bb\in (0,3/4)$.

When $\bb=3/4$, we have $\th_1(3/4)=1/2$ and $\th_2(3/4)=\sqrt{3}/2$. Therefore we obtain
$\rho_{1,\pm}(3/4)=e^{\pm \sqrt{-1} \pi} = -1$ and $\rho_{2,\pm}(3/4)=e^{\pm \sqrt{-1} \sqrt{3}\pi}$.

When $3/4<\bb<1$, the angle $\th_1(\bb)$ increases strictly from $1/2$ to $\sqrt{2}/2$ as $\bb$ increase
from $3/4$ to $1$. Thus $\rho_{1,+}(\bb)=e^{2\pi \sqrt{-1}\th_1(\bb)}$ runs from $-1$ to $e^{\sqrt{-1} \sqrt{2}\pi}$
counterclockwise along the lower semi-unit circle in $\C$ as $\bb$ increases from $3/4$ to $1$. Correspondingly
$\rho_{1,-}(\bb)=e^{-2\pi \sqrt{-1}\th_1(\bb)}$ runs from $-1$ to $e^{-\sqrt{-1} \sqrt{2}\pi}$ clockwise along the
upper semi-unit circle in $\C$ as $\bb$ increases from $3/4$ to $1$. Because $\th_2(\bb)$ decreases strictly
from $\sqrt{3}/2$ to $\sqrt{2}/2$ as $\bb$ increases from $3/4$ to $1$, we obtain that
$\rho_{2,+}(\bb)=e^{2\pi \sqrt{-1}\th_2(\bb)}$ runs from $e^{\sqrt{-1}\sqrt{3}\pi}$ to $e^{\sqrt{-1}\sqrt{2}\pi}$ clockwise along the
lower semi unit circle in $\C$ as $\bb$ increases from $3/4$ to $1$. Correspondingly
$\rho_{2,-}(\bb)=e^{-2\pi \sqrt{-1}\th_2(\bb)}$ runs from $e^{-\sqrt{-1}\sqrt{3}\pi}$ to $e^{-\sqrt{-1}\sqrt{2}\pi}$ counterclockwise
along the upper semi unit circle in $\C$ as $\bb$ increases from $3/4$ to $1$. Thus we obtain
$\sg(\ga_{\bb,0}(2\pi)) \subset \U\bs\R$ for all $\bb\in (3/4,1)$.

When $\bb=1$, we obtain $\th_1(1)=\th_2(1)=\sqrt{2}/2$, and then we have double eigenvalues
$\rho_{1,\pm}(1) = \rho_{2,\pm}(1) = e^{\pm \sqrt{-1}\sqrt{2}\pi}$.

When $1<\bb<9$, using notations defined in (\ref{3.15}), the four characteristic multipliers of
$\ga_{\bb,0}(2\pi)$ satisfy $\sg(\ga_{\bb,0}(2\pi)) \subset \C\bs (\U\cup\R)$ for all $\bb\in (1,9)$.

When $\bb=9$, $\sqrt{\sqrt{9}+1}\pi=2\pi$. By (\ref{3.15}), we get the two positive double
characteristic multipliers of $\ga_{9,0}(2\pi)$ given by
$\rho_{\pm,\pm} = e^{\pm\sqrt{2}\pi}e^{\pm \sqrt{-1}2\pi} = e^{\pm\sqrt{2}\pi}\,\in \R_+\bs\{1\}$, where we denote
by $\R_+=\{r\in\R\,|\,r>0\}$.

{\bf (B)  Normal forms of $\ga_{\bb,0}(2\pi)$ for $\bb\in [0,9]$.}

By our above analysis, the matrix $\ga_{\bb,0}(2\pi)$ possesses no eigenvalues $\pm 1$ for $\bb\in (3/4,9]$.
Note also that the matrix $\ga_{\bb,0}(2\pi)$ is homotopic to $N_2(\om_0,b)$ with $\om_0=e^{\sqrt{-1}\pi\sqrt{2}}$
and $b_2-b_3\not=0$ when $\bb=1$, and is hyperbolic when $\bb\in (1,9]$. Therefore by
Definition \ref{D2.1} we have
\be  \ga_{\bb,0} \sim_{\pm 1} \ga_{\bb_0,0}, \qquad
      \forall\;\bb\in [1,9]\;\;{\rm and}\;\;\bb_0\in (3/4,1). \lb{3.16}\ee
Thus next we are especially interested in the normal forms of $\ga_{\bb,0}(2\pi)$ with $\bb\in [0,1)$.

For that purpose we construct a family of continuous curves $f(\bb)(t)$ in $\Sp(4)$ for $\bb\in [0,1)$ and
$t\in [0,2\pi]$ satisfying $f(\bb)(0)=I$ such that $\ga_{\bb,0}(2\pi)\sim_{\pm 1} f(\bb)$ for all $\bb\in [0,1)$,
which implies $\ga_{\bb,0}(2\pi)\approx f(\bb)$ for all $\bb\in [0,1)$. After that we shall study the normal
forms of $\ga_{\bb,0}(2\pi)$ for $\bb\in [1,9]$.

The curve $f$ is defined separately according to $\bb$ as follows.

(i) {\it Normal forms of $\ga_{\bb,0}(2\pi)$ when $\bb\in [0,3/4]$.}

Here, we splits $f$ into a symplectic sum of two $\Sp(2)$-paths:
\be  f(\bb) = f_+(\bb)\dm f_-(\bb), \qquad f_{\pm}(\bb)\,\in\,\Sp(2). \lb{3.17}\ee

Using the cylindrical coordinate representation (which is denoted by CCR for short below) of $\Sp(2)$
introduced in \cite{Lon1} (cf. pp.48-50 of \cite{Lon4}, specially Figure 1 and 2 on pp.49-50 there),
we can describe the matrix curves in $\Sp(2)$ more precisely, which is shown below in Figure 3.

Using notations in \S 2 we let
\be  M_0=D(2)R(2\pi-\arcsin(3/5))=\left(\matrix{8/5 & 6/5\cr
                                  -3/10 & 2/5\cr}\right).  \lb{3.18}\ee
Then $M_0\approx N_1(1,1)$ and thus
$M_0 \in \Sp(2)^0_{1,-}\cap \{(r,\th,z)\in\R^3\bs\{z{\rm -axis}\}\;|\;z=0\}$ in Figure 2.1.2 on p.50
of \cite{Lon4}. By (\ref{3.3}) and our discussions in part (A), we define
\be  f_+(0) = I_2 \qquad {\rm and}\qquad f_-(0) = M_0.    \lb{3.19}\ee

By our analysis in part (A), when $\bb$ increases from $0$ to $3/4$ in $(0,3/4)$, we define
\be f_+(\bb) = R(2\pi - \frac{4}{3}\pi\bb),   \qquad {\rm for}\quad \bb\in [0,3/4].  \lb{3.20}\ee
Then the matrix curve $f_+(\bb)$ runs from $I_2$ to $R(\pi)=-I_2$ along the left semi-circle clockwise
in CCR in the left diagram of Figure 3 below.

Then we can choose some $\bb_1\in (0,3/4)$ and define
\be  f_-(\bb) = D(2-\frac{\bb}{\bb_1})R((2\pi-\arcsin(3/5))(1-\frac{\bb}{\bb_1})+\bb\frac{9\pi}{5}),
                   \qquad {\rm for}\quad \bb\in [0,\bb_1],  \lb{3.21}\ee
and
\be  f_-(\bb) = R(\sqrt{3}\pi - \frac{3/4-\bb}{3/4-\bb_1}(\sqrt{3}-\frac{9}{5})\pi),
             \qquad {\rm for}\quad \bb\in [\bb_1,3/4].  \lb{3.22}\ee
Thus the matrix curve $f_-(\bb)$ runs from $f_-(0)=M_0$ to $f_-(\bb_1)=R(9\pi/5)$ when $\bb$ runs from
$0$ to $\bb_1$, and then runs from $f_-(\bb_1)=R(9\pi/5)$ to $f_-(3/4)=R(\sqrt{3}\pi)$ when $\bb$ runs
from $\bb_1$ to $3/4$. The image of $f_-(\bb)$ is shown along the left semi-circle clockwise in CCR in
the right diagram of Figure 3 below.

Then by adjusting the running speeds of $f_+(\bb)$ and $f_-(\bb)$ according to those of
the corresponding eigenvalues in (\ref{3.13}) respectively, we can have
\be  f(\bb) = f_+(\bb)\dm f_-(\bb) \approx \ga_{\bb,0}(2\pi), \qquad {\rm for}\quad \bb\in [0,3/4].
                     \lb{3.23}\ee
Specially when $\bb=3/4$, we obtain
\be  f(3/4) = -I_2\dm R(\sqrt{3}\pi).  \lb{3.24}\ee

(ii) {\it Normal forms of $\ga_{\bb,0}(2\pi)$ when $\bb\in [3/4,1)$.}

For $\bb\in [3/4,1)$, following part (A), we define
\bea
f_+(\bb) &=& R(4(1-\sqrt{2})\pi\bb + (3\sqrt{2}-2)\pi),    \lb{3.25}\\
f_-(\bb) &=& R(4(\sqrt{2}-\sqrt{3})\pi\bb + (4\sqrt{3}-3\sqrt{2})\pi).  \lb{3.26}\eea
By adjusting the two curves $f_+(\bb)$ and $f_-(\bb)$ suitably when $\bb<1$ and
close to $1$ and then adjusting their running speeds in $\bb$ suitably, we can suppose that
when $\bb$ increases from $3/4$ to $1$ in $(3/4,1)$, the matrix curve $f_+(\bb)$ runs from $-I_2$
and tends to $R((2-\sqrt{2})\pi)$ along the right semi-circle clockwise in CCR in the left diagram
of Figure 3 below as $\bb$ runs from $3/4$ and tends to $1$. Simultaneously the matrix curve
$f_-(\bb)$ runs from $R(\sqrt{3}\pi)$ and tends to $R(\sqrt{2}\pi)$ along the left semi-circle
clockwise in CCR in the right diagram of Figure 3 below as $\bb$ runs from $3/4$ and tends to $1$.
Here we shall not explain how $f(\bb)=f_+(\bb)\dm f_-(\bb)$ gets to its limit when
$\bb\to 1$. Note that in this case we have also (\ref{3.23}) holds when $\bb\in [3/4,1)$.

(iii) Note that for $\bb\in [0,1)$, we have $f_+(\bb)=R(\aa(\bb))$, where the value of $\aa(\bb)$
is uniquely given by (\ref{3.20}) and (\ref{3.25}). Then every matrix
$f(\bb)=f_+(\bb)\dm f_-(\bb)\,\in\,\Sp(4)$ in (\ref{3.17}) can be reached by a path starting from
$I_4$ in $\Sp(4)$ as follows:
\bea
f_+(\bb)(t) &=& R(\aa(\bb)\frac{t}{2\pi}), \qquad {\rm for}\;\;0\le t\le 2\pi,      \lb{3.27}\\
f_-(\bb)(t) &=& \left\{\matrix{
      D(2t/\pi)R((2\pi-\arcsin(3/5))t/\pi), & {\rm if}\;\;0\le t\le \pi, \cr
      f_-(\bb\frac{t-\pi}{\pi}), & {\rm if}\;\;\pi <t \le 2\pi,  \cr}\right.  \lb{3.28}\eea
for $\bb\in [0,1)$, where we have used the expression of the matrix $M_0$ in (\ref{3.18}).

(iv) {\it The normal form of $\ga_{\bb,0}(2\pi)$ when $\bb=1$.}

Firstly by Corollary \ref{C4.5} below, for fixed $e\in [0,1)$ the index $i_{\om}(\ga_{\bb,e}(2\pi))$
is non-increasing when $\bb$ increases from $0$ to $9$ for $\om\in \U\bs\{1\}$. By Proposition
\ref{P6.1} below, the index $i_{\om}(\ga_{\bb,e}(2\pi))$ can change only when the matrix path
$\ga_{\bb,e}$ in $\bb$ passes through an $\om$ degeneracy point and there should be either one or two
$\om$-degeneracy points and their total $\om$ degenerate multiplicity is $2$.

By part (A), $\bb=1$ is a Krein collision point of the matrix path $\ga_{\bb,0}(2\pi)$ for
$\bb\in [0,9]$ with
$\sg(\ga_{1,0}(2\pi))=\{\om_0, \ol{\om}_0, \om_0, \ol{\om}_0\}$ for $\om_0=e^{\sqrt{-1}\sqrt{2}\pi}\in \U$.
Note that by our above discussions, when $\bb$ increases in the open interval $(0,3/4)$, the curve
$f_+(\bb)$ passes through each $\om$ singular surface in $\Sp(2)$ precisely once for all
$\om\,\in\,\U\bs\R$, which contributes precisely a $1$ to $\sum_{\bb\in [0,9]}\nu_{\om}(\ga_{\bb,0}(2\pi))$.
Therefore by Proposition \ref{P6.1} below we must have
\be  \nu_{\om_0}(\ga_{1,0}(2\pi)) \equiv \dim_{\C}\ker_{\C}(\ga_{1,0}(2\pi)-\om_0 I) = 1. \lb{3.29}\ee
Then, we must have
\be  \ga_{1,0}(2\pi) \approx N_2(\om_0, b)
  = \left(\matrix{R(\sqrt{2}\pi) & b \cr
                  0              & R(\sqrt{2}\pi) \cr}\right)\quad {\rm with}\quad
    b = \left(\matrix{b_1 & b_2 \cr
                      b_3 & b_4 \cr}\right),  \lb{3.30}\ee
satisfying $b_2-b_3\not=0$.

(v) {\it Normal forms of $\ga_{\bb,0}(2\pi)$ when $\bb\in (1,9]$.}

For $\bb\in (1,9)$, by part (A), we have
\be  \ga_{\bb,0}(2\pi) \approx
\left(\matrix{e^{\pi\sqrt{\sqrt{\bb}-1}}R(\pi\sqrt{\sqrt{\bb}+1}) & b \cr
              0 & e^{-\pi\sqrt{\sqrt{\bb}-1}}R(\pi\sqrt{\sqrt{\bb}+1})\cr}\right),  \lb{3.31}\ee
for some matrix $b=b(\bb)$, which is CS-hyperbolic. When $\bb=9$ we get
\be \ga_{9,0}(2\pi) \approx D(e^{\sqrt{2}\pi})\dm D(e^{\sqrt{2}\pi}).  \lb{3.32}\ee

\begin{remark}\label{R3.1} Because $B(t)$ is a constant matrix depending only on $\bb$ when $e=0$,
similarly to what we did for $\ga_{0,0}(t)$ it is possible to compute out the fundamental matrix path
$\ga_{\bb,0}(t)$ explicitly when $\bb>0$. But the computations on $\ga_{\bb,0}(t)$ when $\bb>0$ are rather
delicate and tedious and thus are omitted here. From this computation, especially when $0\le \bb <1$ we obtain
that $\ga_{\bb,0}(t)\approx R(-2\pi\th_1(\bb)t)\dm R(2\pi\th_2(\bb)t)$ for some $\th_i(\bb)$ with $i=1, 2$,
and the $\bb$-paths $R(-2\pi\th_1(\bb))$ and $R(2\pi\th_2(\bb))$ are homotopic respectively to $f_+(\bb)$ and
$f_-(\bb)$ which we constructed above.
\end{remark}

\begin{figure}
\begin{center}
\resizebox{11cm}{5cm}{\includegraphics*[1cm,1cm][24cm,12cm]{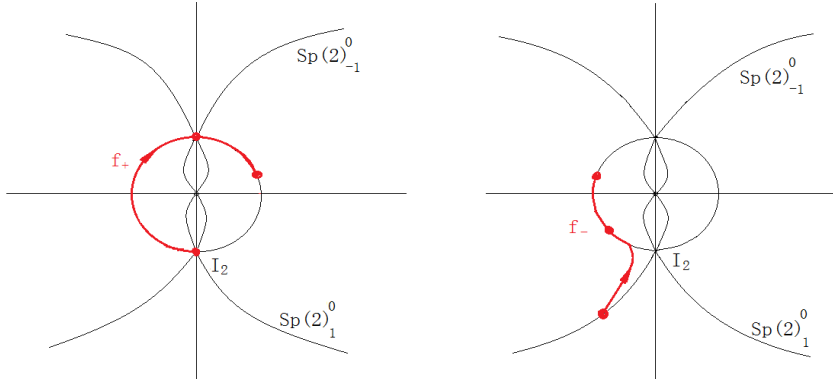}}
\caption{The pictures of the two paths $f_+(\bb)$ and $f_-(\bb)$ in $\Sp(2)$, where
the three dark points are their images at $\bb=0$, $3/4$ and $1$.}
\end{center}
\end{figure}
\vspace{2mm}

{\bf (C)  Indices $i_{\pm 1}(\ga_{\bb,0})$ for $\bb\in [0,9]$.}

From the above discussions as well as Figure 3 we obtain
\bea
&& i_1(f_+(\bb)) = 1, \qquad \forall\,\bb\in [0,1),  \lb{3.33}\\
&& \nu_1(f_+(\bb)) = \left\{\matrix{
                 2, &  {\rm if}\;\;\bb=0, \cr
                 0, &  {\rm if}\;\;\bb\in (0,1). \cr}\right. \lb{3.34}\\
&& i_1(f_-(\bb)) = -1, \qquad \forall\,\bb\in [0,1),  \lb{3.35}\\
&& \nu_1(f_-(\bb)) = \left\{\matrix{
                 1, &  {\rm if}\;\;\bb=0, \cr
                 0, &  {\rm if}\;\;\bb\in (0,1). \cr}\right. \lb{3.36}\eea
Therefore by (\ref{3.17}) we get
\bea
&& i_1(\ga_{\bb,0}) = 0, \qquad \forall\,\bb\in [0,9],  \lb{3.37}\\
&& \nu_1(\ga_{\bb,0}) = \left\{\matrix{
                 3, &  {\rm if}\;\;\bb=0, \cr
                 0, &  {\rm if}\;\;\bb\in (0,9]. \cr}\right. \lb{3.38}\eea

Similarly for the $-1$ index, we obtain
\bea
&& i_{-1}(f_+(\bb)) =  \left\{\matrix{
                 2, &  {\rm if}\;\;\bb\in [0,3/4), \cr
                 0, &  {\rm if}\;\;\bb\in [3/4,1). \cr}\right. \lb{3.39}\\
&& \nu_{-1}(f_+(\bb)) = \left\{\matrix{
                 0, &  {\rm if}\;\;\bb\in [0,1)\bs\{3/4\}, \cr
                 2, &  {\rm if}\;\;\bb=3/4. \cr}\right. \lb{3.40}\\
&& i_{-1}(f_-(\bb)) = 0, \qquad \nu_{-1}(f_-(\bb)) = 0, \qquad {\rm if}\;\;\bb\in [0,1). \lb{3.41}\eea
Therefore by (\ref{3.16}) and Proposition \ref{P6.1} below, we get
\bea
&& i_{-1}(\ga_{\bb,0}) =  \left\{\matrix{
                 2, &  {\rm if}\;\;\bb\in [0,3/4), \cr
                 0, &  {\rm if}\;\;\bb\in [3/4,9]. \cr}\right. \lb{3.42}\\
&& \nu_{-1}(\ga_{\bb,0}) = \left\{\matrix{
                 0, &  {\rm if}\;\;\bb\in [0,9]\bs\{3/4\}, \cr
                 2, &  {\rm if}\;\;\bb=3/4. \cr}\right. \lb{3.43}\eea

The other $\om$-indices of $\ga_{\bb,0}$ can be computed similarly for
$\om\not= e^{\pm \sqrt{-1}\sqrt{2}\pi}\equiv \om_0$, and $i_{\om_0}(\ga_{\bb,0})$
can be computed using the decreasing property of the index proved in Corollary 4.5
and its values at $\bb=0$ and $\bb=9$ as we did in the proof of Theorem \ref{T1.2}
below.

Here we point out specially that $\bb=1$ is the only Krein collision
point on the segment $[0,9]\times \{0\}$.

\setcounter{equation}{0}
\section{Non-degeneracy of elliptic Lagrangian solutions}\label{sec:4}

Note that the complete eigenvalue $1$ non-degeneracy implies that there is no linear
stability change near the positive real line in the complex plane $\C$.

\subsection{Hyperbolicity of elliptic Lagrangian solutions with equal masses}

In the equal mass case, that is $\beta=9$, we have
\be A(9,e)=-\frac{d^2}{dt^2}I_2-I_2+\frac{3}{2(1+e\cos t)}I_2. \lb{4.1}\ee
Let
\be  A_1(e)=-\frac{d^2}{dt^2}-1+\frac{3}{2(1+e\cos t)},  \lb{4.2}\ee
which is a self-adjoint operator with domain $\ol{D}(1,2\pi)\subset W^{2,2}([0,2\pi],\R)$ under the  periodic
boundary conditions $x(t)=x(t+2\pi),\dot{x}(t)=\dot{x}(t+2\pi)$ for all $t\in\R$. Then
\be   A(9,e)=A_1(e)\oplus A_1(e).   \lb{4.3}\ee
By (\ref{2.31}), the Morse index of $A(9,e)$ is zero, so $A_1(e)$ is a
non-negative operator. Moreover, we have
\begin{proposition}\label{P4.1} $A_1(e)>0$ for all $0\leq e<1$ on its domain $\ol{D}(1,2\pi)$.
\end{proposition}

{\bf Proof.} It suffices to show $\ker(A_1(e))=\{0\}$ for all $0\leq e<1$. We argue by
contradiction. Suppose $0\neq x\in\ker(A_1(e))$ is expressed as a Fourier series
$x=x(t)=\sum_{k\in\Z} a_k\exp(\sqrt{-1}kt)$. Then we have
\bea 0
&=& 2(1+e\cos t)A_1(e)x(t)    \nn\\
&=& (2+e\exp(\sqrt{-1}t)+e\exp(-\sqrt{-1}t))\sum_{k\in\Z} a_k(k^2-1)\exp(\sqrt{-1}kt)+\sum_{k\in\Z} 3a_k\exp(\sqrt{-1}kt)   \nn\\
&=& \sum_{k\in\Z}(2a_k(k^2-1)+ea_{k-1}((k-1)^2-1)+ea_{k+1}((k+1)^2-1)+3a_k)\exp(\sqrt{-1}kt). \lb{4.4}\eea
This implies
$$  2a_k(k^2-1)+ea_{k-1}((k-1)^2-1)+ea_{k+1}((k+1)^2-1)+3a_k=0  $$
holds for every $k\in\Z$. Let $k=0$, we have $a_0=0$. So $x$ belongs to the subspace $V$ which
is spanned by $\{\exp{\sqrt{-1}kt},k\neq 0\}$. Note that $(-\frac{d^2}{dt^2}-1)|_V\geq 0$ and
$\frac{3}{2(1+e\cos t)}>0$, so $A_1(e)$ is positive on $V$, which then implies that $x$ must be
zero. This contradiction completes the proof. \hb

\begin{remark}\label{R4.2} For $e=1$, the operator $A_1(e)$ is singular. By the same argument
as above one can show that $A_1(1)x=0$ also implies that $x=0$ in $L^2([0,2\pi],\R)$.
\end{remark}

{\bf Proof of Proposition \ref{P1.4}.} Let $\bar{\xi}_e(t)$ be the fundamental solution
of the first order linear Hamiltonian system corresponding to $A_1(e)$. Then it satisfies
\bea
\frac{d}{dt}\bar{\xi}_e(t)
&=& J_2\left(\matrix{ 1 & 0 \cr
                      0 & 1-\frac{3}{2(1+e\cos t)} \cr}\right)\bar{\xi}_e(t), \lb{4.5}\\
\bar{\xi}_e(0) &=& I_2. \lb{4.6}\eea
Thus we have
\be  \xi_{9,e}(t)=\bar{\xi}_e(t)\diamond \bar{\xi}_e(t), \lb{4.7}\ee
and so
\be  i_\omega(\xi_{9,e}) = 2i_\omega(\bar{\xi}_e), \qquad
\nu_\omega(\xi_{9,e}(2\pi)) = 2\nu_\omega(\bar{\xi}_e(2\pi)).  \lb{4.8}\ee

By Proposition \ref{P4.1}, we have
\be  \dim\ker(\bar{\xi}_e(2\pi)-I_2)=0, \qquad \forall \,e\in [0,1). \lb{4.9}\ee
Since $\bar{\xi}_e(2\pi)$ is a $2\times 2$ symplectic matrix, its eigenvalues are in pair
$\{\lambda,\lambda^{-1}\}$ with $\lambda$ real and $\lambda\neq 1$ or
$\lambda\in\U\backslash\{1\}$. Because $i_1(\bar{\xi}_e)=0$ by Proposition 4.1, the matrix
$\bar{\xi}_e(2\pi)$ must have normal form
$\left(\matrix{ \lambda & 0 \cr
                0 & \lambda^{-1} \cr}\right)$ with $\lambda>0$ or
$\left(\matrix{ 1 & -1 \cr
                0 & 1 \cr}\right)$ by \cite{Lon4} (cf. pp.179-183 there). By (\ref{4.9}), $1$
can not be an eigenvalue of $\bar{\xi}_e(2\pi)$ for any $e\in [0,1)$, so the only possible
case is that $\bar{\xi}_e(2\pi)$ has a pair of positive real eigenvalues not equal to $1$, that
is $\lambda$ and $\lm^{-1}>0$ and $\lambda\neq 1$, and the eigenvalues of $\bar{\xi}_e(2\pi)$
depend continuously on $e$. So the eigenvalues of $\bar{\xi}_e(2\pi)$ are real positive
numbers not equal to $1$ for any $e\in [0,1)$. Specially from (\ref{4.7}) with $t=2\pi$, (\ref{4.8})
and the above discussions, we obtain
\be  i_\omega(\xi_{9,e}) = 2i_\omega(\bar{\xi}_e)=0, \qquad
\nu_\omega(\xi_{9,e}(2\pi)) = 2\nu_\omega(\bar{\xi}_e(2\pi))=0.  \lb{4.10}\ee
Thus by (\ref{2.27}) and (\ref{2.30}), the Proposition \ref{P1.4} is proved. \hb

\begin{corollary}\label{C4.3} $A(9,e)$ is a positive operator for any $\omega$ boundary conditions.
\end{corollary}

{\bf Proof.} When $\bb=9$, since $i_1(\bar{\xi}_e)=0$ and $\bar{\xi}_e(2\pi)$ is
hyperbolic as we proved above, we have $i_{\omega}(\bar{\xi}_e)=0$ for any $\omega\in\U$
and $e\in [0,1)$. Thus the $\omega$ Morse index $i_{\omega}(\xi_{9,e})=\phi_{\omega}(A(9,e))=0$,
that is, $A(9,e)\geq 0$ in any $\omega$ boundary conditions. Then $A(9,e)>0$ follows from the fact
that $A(9,e)$ is non-degenerate at any $\omega$ boundary conditions by (\ref{4.10}). \hb

\subsection{The non-degeneracy of elliptic Lagrangian solutions for $\omega=1$}

For $(\bb,e)\in [0,9)\times [0,1)$, let $\bar{A}(\bb,e)=\frac{A(\bb,e)}{\sqrt{9-\beta}}$. Using
(\ref{2.29}) we can rewrite $A(\bb,e)$ as follows
\be A(\bb,e) = A(9,e)+\frac{\sqrt{9-\beta}}{2(1+e\cos t)}S(t)
= \sqrt{9-\beta}\left(\frac{A(9,e)}{\sqrt{9-\beta}}+\frac{S(t)}{2(1+e\cos t)}\right)
= \sqrt{9-\beta}\,\bar{A}(\bb,e). \lb{4.11}\ee
Then we have
\bea
\phi_\omega(A(\bb,e)) &=& \phi_\omega(\bar{A}(\bb,e)), \lb{4.12}\\
\nu_\omega(A(\bb,e)) &=& \nu_\omega(\bar{A}(\bb,e)).   \lb{4.13}\eea
Following from Corollary \ref{C4.3}, i.e., the fact that $A(9,e)$ is positive definite for any $\omega$
boundary condition, we get the following important lemma:

\begin{lemma}\label{L4.4} (i) For each fixed $e\in [0,1)$, the operator $\bar{A}(\bb,e)$ is increasing
with respect to $\beta\in [0,9)$ for any fixed $\omega\in\U$. Specially
\be  \frac{\pt}{\pt\beta}\bar{A}(\beta,e)|_{\bb=\bb_0} = \frac{1}{2(9-\bb_0)^{3/2}}A(9,e),  \lb{4.14}\ee
for $\bb\in [0,9)$ is a positive definite operator.

(ii) For every eigenvalue $\lm_{\bb_0}=0$ of $\bar{A}(\bb_0,e_0)$ with $\om\in\U$ for some
$(\bb_0,e_0)\in [0,9)\times [0,1)$, there holds
\be \frac{d}{d\bb}\lm_{\bb}|_{\bb=\bb_0} > 0.  \lb{4.15}\ee
\end{lemma}

{\bf Proof.} It suffices to prove (ii). Let $x_0=x_0(t)$ with unit norm such that
\be  \bar{A}(\bb_0,e_0)x_0=0.     \lb{4.16}\ee
Fix $e_0$. Then $\bar{A}(\bb,e_0)$ is an analytic path of strictly increasing self-adjoint operators
with respect to $\beta$. Following Kato (\cite{Ka}, p.120 and p.386), we can choose a smooth path of
unit norm eigenvectors $x_{\bb}$ with $x_{\beta_0}= x_0$ belonging to a smooth path of real eigenvalues
$\lm_{\bb}$ of the self-adjoint operator $\bar{A}(\bb,e_0)$ on $\ol{D}(\om,2\pi)$ such that for small
enough $|\beta-\beta_0|$, we have
\be  \bar{A}(\bb,e_0)x_\beta=\lambda_\beta x_\beta,   \lb{4.17}\ee
where $\lambda_{\beta_0}=0$. Taking inner product with $x_\beta$ on both sides of (\ref{4.17})
and then differentiating it with respect to $\beta$ at $\beta_0$, we get
\bea  \frac{\pt}{\pt\bb}\lambda_{\beta}|_{\bb=\bb_0}
&=& \<\frac{\pt}{\pt\bb}\bar{A}(\bb,e_0)x_{\bb},x_{\bb}\>|_{\bb=\bb_0}
    + 2\<\bar{A}(\bb,e_0)x_{\bb},\frac{\pt}{\pt\bb}x_{\bb}\>|_{\bb=\bb_0}  \nn\\
&=& \<\frac{\pt}{\pt\bb}\bar{A}(\bb_0,e_0)x_0,x_0\>  \nn\\
&=& \frac{1}{2(9-\bb_0)^{3/2}}\<A(9,e_0)x_0,x_0\>   \nn\\
&>& 0, \nn\eea
where the second equality follows from (\ref{4.16}), the last equality follows from the definition
of $\bar{A}(\bb,e)$ and (\ref{4.11}), the last inequality follows from the positive definiteness of
$A(9,e)$ given by Corollary 4.3, and the fact $x_0\not=0$. Thus (\ref{4.15}) is proved. \hb

Consequently we arrive at
\begin{corollary}\label{C4.5} For every fixed $e\in [0,1)$ and $\om\in \U$, the index function
$\phi_{\om}(A(\bb,e))$, and consequently $i_{\om}(\ga_{\bb,e})$, is non-increasing in $\bb\in [0,9]$.
When $\om\in\U\bs\{1\}$, it decreases from $2$ to $0$.
\end{corollary}

{\bf Proof.} For $0\le \bb_1<\bb_2<9$ and fixed $e\in [0,1)$, when $\bb$ increases from $\bb_1$ to
$\bb_2$, it is possible that negative eigenvalues of $\bar{A}(\bb_1,e)$ pass through $0$ to become
positive ones of $\bar{A}(\bb_2,e)$, but it is impossible that positive eigenvalues of
$\bar{A}(\bb_2,e)$ pass through $0$ to become negative by (ii) of Lemma \ref{L4.4}. Therefore the first
claim holds. The second claim follows from (\ref{3.7}), (\ref{3.9}) and Corollary \ref{C4.3}.  \hb

From now on in this section, we will focus on the case of $\omega=1$. Since $\phi_1(A(\bb,e))=0$,
we have $\bar{A}(\bb,e)\geq 0$ for $(\bb,e)\in [0,9)\times [0,1)$. Furthermore, we have

\begin{proposition}\label{P4.6} $A(\bb,e)>0$ for all $(\bb,e)\in (0,9)\times [0,1)$ under the periodic
boundary conditions, i.e., on $\ol{D}(1, 2\pi)$.
\end{proposition}

{\bf Proof.} It suffices to prove $\bar{A}(\bb,e)>0$. This is essentially due to the
fact that $\ker(\bar{A}(\bb,e))=\{0\}$ on $(\bb,e)\in (0,9)\times [0,1)$. In fact, otherwise,
there exists an $x_0=x_0(t)$ with unit norm such that (\ref{4.16}) holds. Then (\ref{4.17})
implies that there exists a negative eigenvalue of $\bar{A}(\bb,e_0)$ when $\beta<\beta_0$ is
sufficiently close to $\bb_0$, which contradicts to $\bar{A}(\bb,e)\geq 0$. \hb

{\bf Proof of Theorem \ref{T1.1}.} Since $\dim\ker(A(\bb,e))=\dim\ker(\ga_{\beta,e}(2\pi)-I_4)$,
we have proved that the elliptic Lagrangian solutions are all non-degenerate on
$(\bb,e)\in (0,9]\times ([0,1)$. It is degenerate when $\beta=0$, by (\ref{3.8}).
The proof is complete. \hb

\setcounter{equation}{0}
\section{The limiting case $e\to 1$}\label{sec:5}

We shall use the sesquilinear form to study the case when $e\rightarrow 1$, and please refer to
Chapter 6 of \cite{Ka} for the details on the sesquilinear form.

In this section we shall deal with the limiting cases $e=1$ and $-1$, then the term $(1+e\cos t)$
in the denominator of the expression of the operator $A(\bb,e)$ and the corresponding functionals
will become zero when $t=(2k+1)\pi$ if $e=1$ or $2k\pi$ if $e=-1$ with $k\in\Z$. Note that the
boundary condition $x(2\pi)=\om x(0)$ is equivalent to the boundary condition $x(t+2\pi)=\om x(t)$
for all $t\in\R$ in the domain $D(\om, 2\pi)$ with $\om\in\U$ of the corresponding functionals.
Therefore in order to move the singular times to the end of the boundary points of the integral
intervals, in this section we use the interval $\Theta(e)=[-\pi,\pi]$ to replace $[0,2\pi]$ when
we study the case $e\to 1$, and keep $\Theta(e)=[0,2\pi]$ when we study the case $e\to -1$ for all
the corresponding integrals.

For any $(\bb,e)\in [0,9]\times (-1,1)$ and $\omega\in\U$, we define the symmetric sesquilinear form
$\Gamma(\bb,e)$ corresponding to the operator $A(\bb,e)$ in (\ref{2.29}) by
\be  \Gamma(\bb,e)(x,y)
=\int_{\Theta(e)}\left[\dot{x}\cdot\dot{y} - x\cdot y
    + \frac{1}{2(1+e\cos t)}((3I_2+\sqrt{9-\beta}S(t))x(t))\cdot y(t)\right]dt,
      \qquad \forall\; x,y\in D(\omega,2\pi). \lb{5.1}\ee
where the domain $D(\omega,2\pi)$ is defined in (\ref{2.10}). Then we denote by
$\Gamma(\bb,e)(x)= \Gamma(\bb,e)(x,x)$ the corresponding quadratic form. Note that here we have extended
the range of $e$ to $(-1,1)$, even to the whole complex plan in the next section. This seems artificial
from the point of view of celestial mechanics, however we can draw some interesting conclusions on the
degeneracy curves in the $(\beta,e)$ rectangle $[0,9]\times (-1,1)$ later. Note that our results proved
for $e\in [0,1)$ in the previous and later sections hold also for $e\in (-1,1)$, which we shall not
specially indicate out explicitly later. Since $\Gamma(\bb,e)+I$ is equivalent to $\Gamma(\bb,e)$ with
respect to the $W^{1,2}$-norm, $\Gamma(\bb,e)$ is closed on domain $D(\omega,2\pi)$ (\cite{Ka}, p.313).

For $e=\pm 1$, we further define $\Gamma(\bb,e)$ as in (\ref{5.1}), but its domain needs to be modified. More
precisely, say, for $e=1$, the domain of $\Gamma(\bb,1)$ is defined by
\be  \hat{D}(\beta) = \left\{x\in W^{1,2}([-\pi,\pi],\R^{2})\;\left|\;
   \int_{-\pi}^{\pi}\left(\frac{1}{2(1+\cos t)}((3I_2+\sqrt{9-\beta}S(t))x(t))\cdot x(t)\right)dt < \infty \right.\right\}.
            \lb{5.2}\ee
Since $3I_2+\sqrt{(9-\beta)}S(t)>0$ for $0<\beta\leq 9$, we have $\hat{D}(\beta)=\hat{D}(9)$ for $0<\beta\leq9$
and $\hat{D}(0)$ is different from them. Note that every $x\in\hat{D}(9)$ must satisfy $x(\pi)=x(-\pi)=0$ which is
the vanishing $\omega$ boundary condition. And $\Gamma(\beta,1)$ is closed since it is the sum of two closed
symmetric forms.

For $e\in (-1,1)$, $A(\beta,e)$ is the Friedrichs extension operator (\cite{Ka}, Theorem 2.1 in p.322) of
$\Gamma(\beta,e)$ under the $\omega$ boundary condition. Let $A(\beta)$ be the Friedrichs extension operator of
$\Gamma(\beta,1)$. Then it has the form
\be  A(\beta)=-\frac{d^2}{dt^2}I_2-I_2+\frac{3}{2(1+\cos t)}I_2+\frac{\sqrt{9-\beta}S(t)}{2(1+\cos t)}, \lb{5.3}\ee
with $\mathrm{dom}(A(\beta))\subset \hat{D}(\beta)$.

\begin{lemma}\label{L5.1} For each $\beta\in (0,9]$, $A(\beta)$ is a self-adjoint operator on $L^2([-\pi,\pi],\R^2)$
with compact resolvent.
\end{lemma}

{\bf Proof.} Since $\Gamma(\beta,1)$ is symmetric, $A(\beta)$ is self-adjoint. It suffices to prove that
$A(\beta)+2I_2$ has a compact inverse. Let $\Gamma'(x)=\Gamma(\beta,1)(x)+2\|x\|^2_{W^{1,2}}$. Then
$A(\beta)+2I_2$ is the Friedrichs extension operator of $\Gamma'$. By the second representation theorem
(\cite{Ka}, p.331),
\be  \Gamma'(x)=\<(A(\beta)+2I_2)^{\frac{1}{2}}x,(A(\beta)+2I_2)^{\frac{1}{2}}x\>. \lb{5.4}\ee
Since $\|x\|^2_{W^{1,2}}\leq \Gamma'(x)$, a $\Gamma'$-bounded set must have a convergent subsequence in
$L^2$. This implies that $(A(\beta)+2I_2)^{-\frac{1}{2}}$ is compact, so $(A(\beta)+2I_2)^{-1}$ is also compact
in $L^2$. \hb

\begin{lemma}\label{L5.2} For $\beta\in (0,9]$ and $x\in \hat{D}(\beta)$, we have
$\Gamma(\beta,e)(x)\rightarrow \Gamma(\beta,1)(x)$ as $e\rightarrow 1$.
\end{lemma}

{\bf Proof.} Let $x\in \hat{D}(\bb)\subset D(\om,2\pi)$ and $e\in [0,1)$. By definitions of $\Ga(\bb,e)$ and
$\Ga(\bb,1)$, when $\cos t>0$, we obtain
\bea
&& \left|(\frac{1}{2(1+e\cos t)}-\frac{1}{2(1+\cos t)})((3I_2+\sqrt{9-\beta}S(t))x(t)\cdot x(t))\right|  \nn\\
&&\qquad\qquad\qquad =\; \frac{1}{2}\left|\frac{(1-e)\cos t}{(1+e\cos t)(1+\cos t)}
                                       ((3I_2+\sqrt{9-\beta}S(t))x(t)\cdot x(t))\right|    \nn\\
&&\qquad\qquad\qquad \le\; (1-e)\left|\frac{1}{2(1+\cos t)}((3I_2+\sqrt{9-\beta}S(t))x(t)\cdot x(t))\right|.    \nn\eea
When $\cos t\le 0$, we have
$$ \left|\frac{1}{2(1+e\cos t)}((3I_2+\sqrt{9-\beta}S(t))x(t)\cdot x(t))\right|
   \le \left|\frac{1}{2(1+\cos t)}((3I_2+\sqrt{9-\beta}S(t))x(t)\cdot x(t))\right|.  $$
Therefore we get
\be \left|\frac{1}{2(1+e\cos t)}((3I_2+\sqrt{9-\beta}S(t))x(t)\cdot x(t))\right|
   \le (2-e)\left|\frac{1}{2(1+\cos t)}((3I_2+\sqrt{9-\beta}S(t))x(t)\cdot x(t))\right|,  \lb{5.5}\ee
for all $t\in [0,2\pi]$. Now the lemma follows from the Lebesgue's dominated convergence theorem. \hb

Since under the periodic boundary condition, $\Gamma(\beta,e)>0$ for $e\in[0,1)$, the above lemma tells
us that $\Gamma(\beta,1)\geq 0$. By Remark \ref{R4.2}, we have $\ker(A(9,1))=\{0\}$, so $\Gamma(9,1)>0$.
By completely the same reasoning as in Proposition \ref{P4.6}, we have
\be \Gamma(\beta,1)>0, \qquad \forall\beta\in (0,9].  \lb{5.6}\ee

For the limiting case $e\to 1$ under general boundary conditions, we consider the sesquilinear form
$\hat{\Gamma}(\beta,e)$, for $\beta\in [0,9]$ and $e\in (-1,1)$, by
\be \hat{\Gamma}(\beta,e)(x,y)=
  \int_{-\pi}^{\pi}\left[\dot{x}\cdot\dot{y} - x\cdot y + (\frac{3I_2+\sqrt{9-\beta}S(t)}{2(1+e\cos t)}x(t)\cdot y(t))\right]dt,
        \qquad \forall x, y\in W^{1,2}(\R/2\pi\Z,\R^2). \lb{5.7}\ee
We have the following

\begin{lemma}\label{L5.3} For $\beta\in (0,9]$, we have $\hat{\Gamma}(\beta,e)>0$ when $1-e$ is small enough.
\end{lemma}

{\bf Proof.} Let $\delta(\beta,e)$ be the largest lower bound of the quadratic form $\hat{\Gamma}(\beta,e)$
for $e\in [0,1]$. Then by (\ref{5.6}) we have
\be   \delta(\beta,1) > 0 \qquad \forall\; \bb\in (0,9].  \lb{5.8}\ee
We need to show
\be   \liminf_{e\to 1}\delta(\beta,e)>0.  \lb{5.9}\ee

For $(\bb,e)\in (0,9]\times [0,1)$, we define
\be f_{\bb,e}(t) = \left\{\matrix{
    \frac{3I_2+\sqrt{9-\beta}S(t)}{2(1+e\cos t)} &\quad {\rm if}\;\; \cos t\leq 0 \cr
    \frac{3I_2+\sqrt{9-\beta}S(t)}{2(1+\cos t)} &\quad {\rm if}\;\; \cos t>0 \cr}\right., \lb{5.10}\ee
and
\be  \tilde{\Gamma}(\beta,e)(x) = \int_{-\pi}^{\pi}\left[\dot{x}\cdot\dot{x} - x\cdot x
             + f_{\bb,e}(t)x(t)\cdot x(t)\right]dt.            \lb{5.11}\ee
Let $\tilde{\delta}(\beta, e)$ be the largest lower bound of $\tilde{\Gamma}(\beta,e)$ on
$W^{1,2}(\R/(2\pi\Z),\R^2)$ for $e\in [0,1]$.

Then by (\ref{5.5}) we obtain
\be  \liminf_{e\to 1}\delta(\beta,e) = \liminf_{e\to 1}\tilde{\delta}(\beta,e).  \lb{5.12}\ee
So it suffices to prove
\be  \liminf_{e\to 1}\tilde{\delta}(\beta, e) > 0.  \lb{5.13}\ee

For $e_2>e_1$, we have
$$  \tilde{\Gamma}(\beta,e_2)(x) - \tilde{\Gamma}(\beta,e_1)(x)
    = \int_{\pi/2}^{\pi}+\int_{-\pi}^{-\pi/2}
    \left(\frac{(e_2-e_1)(-\cos t)}{(1+e_2\cos t)(1+e_1\cos t)}((3I_2+\sqrt{9-\beta}S(t))x(t)\cdot x(t))\right)dt. $$
Note that $3I_2+\sqrt{9-\beta}S(t)$ is positive definite whenever $\bb\in (0,9]$. Thus we obtain
\be  \tilde{\Gamma}(\beta,e_2) \ge \tilde{\Gamma}(\beta,e_1), \qquad {\rm if}\quad e_2>e_1. \lb{5.14}\ee

So $\tilde{\delta}(\beta,e)$ is increasing with respect to $e$. Let
\be  E(\beta,e) = \left\{x\in W^{1,2}(\R/(2\pi\Z), \R^2)\,\left|\, \|x\|_{W^{1,2}}\leq 1, \quad
   \tilde{\Gamma}(\beta,e)(x)\leq \frac{1}{2}\delta(\beta,1)\|x\|_{L^2}^2 \right. \right\}. \lb{5.15}\ee
Then $E(\beta,e)$ is closed in the Hilbert space $W^{1,2}(\R/2\pi\Z,\R^2)$ and $E(\beta,e_2)\subset E(\beta,e_1)$
if $e_2>e_1$. Let $E(\beta,1)=\bigcap_{0<e<1}E(\beta,e)$. We claim
\be  E(\beta,1)=\{0\}. \lb{5.16}\ee

In fact, otherwise, there exists some $x\in E(\beta,1)\bs\{0\}$. We consider two cases depending on
whether $x\in\hat{D}(\beta)$ or not. If $x\in \hat{D}(\beta)$, then $\tilde{\Gamma}(\beta,e)(x)$
converges to $\Gamma(\beta,1)(x)$ when $e\to 1$ by an argument similar to that of Lemma \ref{L5.2}.
Therefore by the definition of $\delta(\bb,1)$ we obtain
$$  \tilde{\Gamma}(\beta,e)(x) > \frac{1}{2}\delta(\beta,1)\|x\|^2,   $$
when $1-e$ is small enough. This contradicts the definition (\ref{5.15}). On the other hand,
$x\in E(\beta,1)\setminus \hat{D}(\beta)$ implies that
$$ \int_{-\pi}^{\pi}\left(\frac{3I_2+\sqrt{9-\beta}S(t)}{2(1+\cos t)}x(t)\cdot x(t)\right)dt $$
is infinite, which contradicts to the definitions of $\Ga(\bb,e)$ and $\Ga(\bb,1)$ as well as
Levi's theorem. Thus (\ref{5.16}) holds.

We then further claim that there exists a constant $e_0\in (0,1)$ such that
\be  E(\beta,e) = \{0\}, \qquad {\rm whenever}\quad e>e_0.  \lb{5.17}\ee

In fact, otherwise, there exists an increasing sequence $e_k\in (0,1)$ with $k\in\N$ such
that $e_k\to 1$ and there exist $x_k\in E(\beta,e_k)$ with $\|x_k\|_{W^{1,2}}=1$. Since
$x_k$ is bounded in $W^{1,2}(\R/(2\pi\Z), \R^2)$, it has a weakly convergent subsequence
$x_{n_k}$, which converges weakly to some $x_0$. Then we have $x_0\in E(\beta,1)$. In
fact, by the weakly lower semi-continuity of norms, we obtain
$$   \|x_0\|_{W^{1,2}} \le \liminf_{k\to \infty} \|x_{n_k}\|_{W^{1,2}}.   $$
On the other hand, by the Sobolev compact embedding theorem, $x_{n_k}$ converges to $x_0$
in the $L^2$ space. By definition, $x\in E(\beta,e)$ is equivalent to
\be   \|x\|_{W^{1,2}}+\int_{-\pi}^{\pi}(f_{\bb,e}(t)x(t)\cdot x(t))dt
           \le (\frac{1}{2}\delta(\beta,1)+2)\|x\|_{L^2}^2. \lb{5.18}\ee
Thus $x_0\in E(\bb,e_{n_k})$ for every $k\in\N$ implies $x_0\in E(\beta,1)$.
On the other hand, (\ref{5.18}) implies that the lower bound of  $\|x_{n_k}\|_{L^2}^2$ is nonzero.
So we have $x_0\neq 0$, which contradicts to (\ref{5.16}), and proves (\ref{5.17}).

Now by (\ref{5.8}), (\ref{5.14}) and (\ref{5.16}), for fixed $\bb\in (0,1]$ and every $e>e_0$
we obtain $\tilde{\delta}(\beta,e)\ge \delta(\beta,1)/2>0$, which completes the proof. \hb

Since $D(\omega,2\pi)\subset W^{1,2}(\R/2\pi\Z,\C^n)$ for any $\omega\in\U$, then
$\hat{\Gamma}(\beta,e)>0$ implies $\Gamma(\beta,e)>0$ for any $\omega\in\U$. Thus we have

\begin{corollary}\label{C5.4} For any fixed $\beta\in (0,9]$, there exists an $e_*\in (0,1)$ such
that for any $\omega\in\U$ there holds $\Gamma(\beta,e)>0$ for all $e\in [e_*,1]$.
\end{corollary}

{\bf Proof of Theorem \ref{T1.7}} The proof of the limiting case $e\to -1$ is similar and thus is
omitted. The above lemmas imply that for any fixed $\bb\in (0,9]$, we have always $A(\beta,e)>0$
for any $\omega\in\U$ whenever $1-|e|$ is small enough. Thus Theorem \ref{T1.7} is proved. \hb

\setcounter{equation}{0}
\section{The $\om$ degeneracy curves of elliptic Lagrangian solutions}\label{sec:6}

For any $\omega$ boundary condition, that is on domain $\ol{D}(\omega,2\pi)$ of (\ref{2.11}),
$A(\beta,e)$ is a closed unbounded operator. If we extend $e$ to the complex plane and denote
the open unit disc by $D=\{e\in\C\,|\,|e|<1\}$, then $A(\beta,e)$ is holomorphic with respect
to $e\in D$ (\cite{Ka}, p.366). It satisfies $A(\beta,\bar{e})=A(\beta,e)^*$. In fact
\be  A(\beta,e) = - \frac{d^2}{dt^2}I_2 - I_2
   + \frac{1}{2}\left(3I_2+\sqrt{9-\beta}S(t))(1-e\cos(t)+e^2\cos^2(t)-e^3\cos^3(t)+\cdots\right),
                      \lb{6.1}\ee
where we have
\bea
3I_2+\sqrt{9-\beta}S(t)\geq 0, &\qquad& {\rm for}\quad 0\leq\beta\leq 9, \nn\\
3I_2+\sqrt{9-\beta}S(t)>0, &\qquad& {\rm for}\quad 0<\beta\leq 9.  \nn\eea
Let $\Omega$ be a small narrow neighborhood of the interval $(-1,1)$ in
the complex plane. For $e\in\Omega$, if its imaginary part $\Im e$ is small enough, $A(9,e)$
is strictly accretive (\cite{Ka}, p.281), i.e., there exists $\delta>0$ such that the real part
$$   \Re(A(9,e)x,x) \ge \delta\|x\|^2   $$
for any $x$ in the domain $\ol{D}(\omega,2\pi)$ of $A(9,e)$. We get $A(9,e)^{-\frac{1}{2}}$ by
the Dunford-Taylor integral (\cite{Ka}, (3.43) in p.282),
\be  A(9,e)^{-\frac{1}{2}} = \frac{1}{\pi}\int_0^{\infty}\mu^{-\frac{1}{2}}(A(9,e)+\mu)^{-1}d\mu,
             \lb{6.2}\ee
which is bounded and holomorphic in $e$ (\cite{Ka}, p.398).

By the idea of the proof of Proposition \ref{P4.6}, we have the following important
results for all $\om\in \U\bs\{1\}$.

\begin{proposition}\label{P6.1} (i) For every $(\bb,e)\in (0,9)\times [0,1)$ and $\om\in \U\bs\{1\}$,
there exists $\ep_0=\ep_0(\bb,e)>0$ sufficiently small such that for all $\ep\in (0,\ep_0]$ there holds
$$  i_{\om}(\ga_{\bb-\ep,e}) - i_{\om}(\ga_{\bb,e}) = \nu_{\om}(\ga_{\bb,e}).   $$

(ii) For every $e\in [0,1)$ and $\om\in\U\bs\{1\}$, the total multiplicity of $\om$ degeneracy of
$\ga_{\bb,e}(2\pi)$ for $\bb\in [0,9]$ is always precisely $2$, i.e.,
$$  \sum_{\bb\in [0,9]}\nu_{\om}(\ga_{\bb,e}(2\pi)) = 2, \qquad \forall\;\om\in\U\bs\{1\}. $$
\end{proposition}

{\bf Proof.} Firstly, for $(\bb,e)\in [0,9)\times [0,1)$, the operator $\bar{A}(\bb,e)$
is a bounded perturbation of the operator $-d^2/dt^2$. Thus as $-d^2/dt^2$ the operator
$\bar{A}(\bb,e)$ possesses only point spectrum, finite Morse index, each of its eigenvalues has
finite multiplicity, and the only accumulation point of its spectrum is $+\infty$. Consequently
its eigenvalues are all isolated.

(i) Fix $\om\in\U\bs\{1\}$ and $e\in [0,1)$. Let $\eta(\bb)$ be a unit norm eigenvector
belonging to an eigenvalue $\lm(\bb)$ of the operator $\bar{A}(\bb,e)$ on $\ol{D}(\om,2\pi)$
for $\bb$ near some $\bb_0\in (0,9)$. Then as in Proposition \ref{P4.6}, we obtain
\bea \left.\frac{d}{d\bb}{\lm}(\bb)\right|_{\bb=\bb_0}
&=& \left.\<\frac{\pt}{\pt\bb}\bar{A}(\bb,e)\eta(\bb),\eta(\bb)\>\right|_{\bb=\bb_0}  \nn\\
&=& \frac{1}{2(9-\bb_0)^{3/2}}\<A(9,e)\eta(\bb_0),\eta(\bb_0)\>   \nn\\
&>& 0.  \lb{6.3}\eea
Therefore eigenvalues of $\bar{A}(\bb,e)$ on $\ol{D}(\om,2\pi)$ for $\bb\in [0,9)$ are strictly
increasing in $\bb$.

Note that by definition the Morse index $k_-=\phi_{\om}(\bar{A}(\bb_0,e))$ is the total multiplicity
of the negative eigenvalues of $\bar{A}(\bb_0,e)$, which is finite.

Suppose $k_0=\dim\ker(\bar{A}(\bb_0,e))>0$ holds on the domain $\ol{D}(\om,2\pi)$. By (\ref{6.3})
there is a smallest positive eigenvalue $\lm_+(\bb_0)$ of the operator $\bar{A}(\bb_0,e)$. Because
$\bar{A}(\bb,e)$ depends analytically on $\bb$, we can choose $\ep>0$ to be small enough so that all
the negative eigenvalues of the operator $\bar{A}(\bb,e)$ with
$\bb\in [\bb_0-2\ep,\bb_0+2\ep]\subset (0,9)$ come only from perturbations of negative and zero
eigenvalues of $\bar{A}(\bb_0,e)$, and are not perturbations from any eigenvalues of $\bar{A}(\bb_0,e)$
larger than or equal to $\lm_+(\bb_0)$. Therefore by (\ref{6.3}) we obtain
\be  \phi_{\om}(\bar{A}(\bb_0-\ep,e)) - \phi_{\om}(\bar{A}(\bb_0,e))
         = \dim\ker(\bar{A}(\bb_0,e)),  \lb{6.4}\ee
together with Lemma \ref{L2.3} and (\ref{2.30}), which yields (i).

(ii) Note first that by Lemma \ref{L2.3}, (\ref{2.30}), (\ref{3.8}) and (\ref{3.9}), the operator
$\bar{A}(\bb,e)$ on $\ol{D}(\om,2\pi)$ for $\bb=0$ or $9$ and $\om\in\U\bs\{1\}$ is non-degenerate.

Following our discussions in (i), at every $\bb_0\in (0,9)$ such that $\bar{A}(\bb_0,e)$ is degenerate,
the $\om$-index must decrease strictly. But by (\ref{3.7}), (\ref{3.9}) and Corollary \ref{C4.5},
there exist at most two values $\bb_1$ and $\bb_2$ at each of which the $\om$-index decreases by $1$
if $\bb_1\not=\bb_2$, or the $\om$-index decreases by $2$ if $\bb_1=\bb_2$. Therefore there exist at
most two $\bb$s in $[0,9]$ at which the operator $\bar{A}(\bb,e)$ degenerates by (i), which we
denote by $\bb_1(e)$ and $\bb_2(e)\in (0,9]$. Thus by Corollary \ref{C4.5} again, we can choose
$\ep>0$ small enough according to $\bb_1(e)$ and $\bb_2(e)$ in the above way so that we have
$\phi_{\om}(\bar{A}(0,e))=\phi_{\om}(\bar{A}(\bb_1(e)-\ep,e))$,
$\phi_{\om}(\bar{A}(\bb_1(e),e))=\phi_{\om}(\bar{A}(\bb_2(e)-\ep,e))$,
$\phi_{\om}(\bar{A}(9,e))=\phi_{\om}(\bar{A}(\bb_2(e),e))$, and (\ref{6.4}) holds for $\bb_0$
replaced by $\bb_1(e)$ and $\bb_2(e)$. Then this yields
\bea 2
&=& \phi_{\om}(\bar{A}(0,e)) - \phi_{\om}(\bar{A}(9,e))  \nn\\
&=& \phi_{\om}(\bar{A}(\bb_1(e)-\ep,e)) - \phi_{\om}(\bar{A}(\bb_1(e),e))
       + \phi_{\om}(\bar{A}(\bb_2(e)-\ep,e)) - \phi_{\om}(\bar{A}(\bb_2(e),e))   \nn\\
&=& \dim\ker(\bar{A}(\bb_1(e),e)) + \dim\ker(\bar{A}(\bb_2(e),e))   \nn\\
&=& \nu_{\om}(\ga_{\bb_1(e),e}(2\pi)) + \nu_{\om}(\ga_{\bb_2(e),e}(2\pi))   \nn\\
&=& \sum_{\bb\in [0,9]}\nu_{\om}(\ga_{\bb,e}(2\pi)),  \nn\eea
which proves the proposition. \hb

Now set
\be  B(e,\omega) = A(9,e)^{-\frac{1}{2}}\frac{1}{2(1+e\cos(t))}S(t)A(9,e)^{-\frac{1}{2}}.
             \lb{6.5}\ee
Be aware that $B(e,\omega)$ depends on $\omega$, since so is $A(9,e)$ on its domain
$\ol{D}(\omega,2\pi)$. Now we have

\begin{lemma}\label{L6.2} For any $\omega$ boundary condition and $e\in \Omega$, $A(\beta,e)$
is $\omega$ degenerate if and only if $\frac{-1}{\sqrt{9-\beta}}$ is an eigenvalue of $B(e,\omega)$.
\end{lemma}

{\bf Proof.} Suppose for $e\in \Omega$, $A(\beta,e)x=0$ holds for some $x\in\ol{D}(\omega,2\pi)$.
Let $y=A(9,e)^{\frac{1}{2}}x$. Then by (\ref{4.11}) we obtain
\bea
&& A(9,e)^{\frac{1}{2}}(\frac{1}{\sqrt{9-\beta}}+B(e,\omega))y(t)    \nn\\
&&\qquad\quad = \left(\frac{1}{\sqrt{9-\bb}}A(9,e) + \frac{1}{2(1+e\cos t)}S(t)\right)x(t)  \nn\\
&&\qquad\quad = \frac{1}{\sqrt{9-\bb}}A(\bb,e)x  \nn\\
&&\qquad\quad = 0.  \lb{6.6}\eea

Conversely, if $(\frac{1}{\sqrt{9-\beta}}+B(e,\omega))y=0$, then $x=A(9,e)^{-\frac{1}{2}}y$ is an
eigenfunction of $A(\beta,e)$ belonging to the eigenvalue $0$ by our computations in
(\ref{6.6}). \hb

\begin{theorem}\label{T6.3} For any $\om\in\U$, there exist two analytic $\om$ degeneracy curves
$(\bb_i(e,\om),e)$ in $e\in (-1,1)$ with $i=1$ and $2$. Specially, each $\beta_i(e,\omega)$ is a real
analytic function in $e\in (-1,1)$, and $0<\bb_i(e,\om)< 9$ and $\ga_{\bb_i(e,\om),e}(2\pi)$ is
$\om$ degenerate for $\om\in \U\bs\{1\}$ and $i=1$ or $2$.
\end{theorem}

{\bf Proof.} For $\omega=1$, we have $\beta_i(e,1)\equiv 0$ for $i=1$ and $2$, by Theorem \ref{T1.1}
and (\ref{2.31}), which is obviously analytic.

For $\omega\in \U\bs\{1\}$, from (\ref{3.7}) we have
$\phi_{\omega}(A(0,e))=\phi_\omega(\bar{A}(0,e))=2$. On the other hand,
$\phi_{\omega}({A}(9,e))=0$ by (\ref{3.9}). Recall that
$\bar{A}(\beta,e)=\frac{A(9,e)}{\sqrt{9-\beta}}$, and it is strictly increasing with
respect to $\beta$ by Lemma \ref{L4.4}. This shows that, for fixed
$e\in (-1,1)$, there are exactly two values $\bb=\beta_1(e,\omega)$ and $\beta_2(e,\omega)$ at which
(\ref{6.6}) is satisfied, and then $\bar{A}(\beta,e)$ at these two $\bb$ values is $\om$ degenerate.
Note that these two $\bb$ values are possibly equal to each other at some $e$ (compare
with the figure in \cite{MSS}), which is not needed in this proof.

Since $\beta=9$ is $\omega$-non-degenerate for any $\omega\in\U$,  we must have $\beta_i(e,\omega)\neq 9$
for $i=1$ and $2$. By Lemma \ref{L6.2}, $\frac{-1}{\sqrt{9-\beta_i(e,\omega)}}$ is an eigenvalue of
$B(e,\omega)$. Note that $B(e,\omega)$ is a compact operator and self-adjoint when $e$ is real.
Moreover it depends analytically on $e$. By \cite{Ka} (Theorem 3.9 in p.392), we know that
$\frac{-1}{\sqrt{9-\beta_i(e,\omega)}}$ with $i=1$ or $2$ is real analytic in $e$. This in turn
implies that both $\beta_1(e,\omega)$ and $\beta_2(e,\omega)$ are real
analytic functions of $e$. \hb

Now we can give

{\bf Proof of Theorem \ref{T1.5}} By Theorem \ref{T6.3}, $\beta_i(e,\omega)$ is real analytic on
$e\in [0,1)$ for $i=1$ or $2$. That $\beta_i(e,\omega)\rightarrow 0$ as $e\rightarrow 1$ for $i=1, 2$
follows by the arguments in the proof of (ii) of Theorem \ref{T1.2} below in Section 8. And Corollary
\ref{C4.5} tells us that $i_\omega(\gamma_{\beta,e})$ is decreasing with respect to
$\beta\in [0,9]$. \hb

Recall that the $\omega$ boundary condition is $x(t)=\omega x(t+2\pi)$,
$\dot{x}(t)=\omega \dot{x}(t+2\pi)$. Let $\psi(x)(t)=x(t+\pi)$, and obviously, $\psi$ preserves
the $\omega$ boundary condition. Also it is a unitary operator and $\psi^*=\psi^{-1}$ is given by
$\psi^*(x)(t)=x(t-\pi)$. One can show that
\be \psi^*A(\beta,e)\psi=A(\beta,-e). \lb{6.8}\ee
In fact,
\bea \psi^*A(\beta,e)\psi x(t)
&=& \psi^*\left(-\frac{d^2}{dt^2}I_2-I_2+\frac{1}{2(1+e\cos t)}
               (3I_2+\sqrt{9-\beta}S(t))\right)x(t+\pi) \nonumber \\
&=& \psi^*\left(-\frac{d^2}{dt^2}I_2-I_2+\frac{1}{2(1-e\cos(t+\pi))}
        (3I_2+\sqrt{9-\beta}S(t+\pi))\right)x(t+\pi) \nonumber \\
&=& A(\beta,-e)x(t). \lb{6.9}\eea
By this property, we know that the $\omega$ degeneracy curve must be symmetric with respect to
$e=0$.

When $e=0$, the eigenvalues of $\ga_{\beta,0}(2\pi)$ has been studied in \S 3.3. Specially
$A(\beta,0)$ has no multiple eigenvalues for $\omega\in \U\bs\{\pm 1\}$ and $0<\beta<9$. So we
have

\begin{theorem}\label{T6.4} For any fixed $\omega\in \U\bs\{\pm 1\}$ and $i=1$ or $2$, the function
$\beta_i(e,\omega)$ is real analytic and even on the interval $(-1,1)$.
\end{theorem}

It then follows that $\frac{\pt}{\pt e}\beta_i(e,\omega)|_{e=0}=0$ when $\omega\in \U\bs\{\pm 1\}$.
But it is not the case when $\omega=-1$, which we now turn to.

\setcounter{equation}{0}
\section{The $-1$ degeneracy curves of elliptic Lagrangian solutions}\label{sec:7}

\subsection{The two $\omega=-1$ degeneracy curves}

For the $\omega=-1$ boundary condition, denote by $g$ the following operator
\be  g(z)(t)=Nz(2\pi-t),   \lb{7.1}\ee
where $N=\left(\matrix{ 1 & 0\cr 0 & -1 \cr}\right)$. Obviously, $g^2=1$ and $g$ is unitary on
$L^2([0,2\pi],\R^2)$. Recall that $E=\ol{D}(-1,2\pi)$ is given by (\ref{2.11}). One can check directly that
\be   A(\beta,e)g=gA(\beta,e).   \lb{7.2}\ee
Let $E_+=\ker(g+I)$ and $E_-=\ker(g-I)$. Following the studies in \S 2.2 and specially
the proof of Theorem 1.1 in \cite{HS1}, the subspaces $E_+$ and $E_-$ are $A(\beta,e)$-orthogonal,
and $E=E_+\oplus E_-$. Note that element $z=(x,y)^T$ in $E_-$ satisfies
$$   x(2\pi-t) = x(t), \quad -y(2\pi-t) = y(t), \qquad \forall\;t\in [0,2\pi].  $$
Thus for all $z=(x,y)^T\in E_-$ we have
$$  x(\pi+t) = x(\pi-t), \quad  -y(\pi+t) = y(\pi-t), \qquad \forall\;t\in [0,2\pi]. $$
By the definition of $z\in \ol{D}(-1,2\pi)$, we have $z(2\pi)=-z(0)$. Thus we have
$x(0)=x(2\pi)=-x(0)$ which implies $x(0)=0$ and $y(\pi)=0$. Similarly for all $z=(x,y)^T\in E_+$
we have
$$  x(\pi+t) = -x(\pi-t), \quad  y(\pi+t) = y(\pi-t), \qquad \forall\;t\in [0,2\pi], $$
and $x(\pi)=0$ and $y(0)=0$.

Therefore by the above discussions, the subspaces $E_-$ and $E_+$ are isomorphic to
the following subspaces $E_1$ and $E_2$ respectively:
\bea
E_1 &=& \{z=(x,y)^T\in W^{2,2}([0,\pi],\R^2)\;|\;x(0)=0,\;y(\pi)=0\}, \lb{7.3}\\
E_2 &=& \{z=(x,y)^T\in W^{2,2}([0,\pi],\R^2)\;|\;x(\pi)=0,\;y(0)=0\}.  \lb{7.4}\eea
For $(\bb,e)\in [0,9]\times [0,1)$, restricting $A(\beta,e)$ to $E_1$ and $E_2$ respectively,
we then obtain
\bea
\phi_{-1}(A(\beta,e)) &=& \phi(A(\beta,e)|_{E_1})+\phi(A(\beta,e)|_{E_2}), \lb{7.5}\\
\nu_{-1}(A(\beta,e)) &=& \nu(A(\beta,e)|_{E_1})+\nu(A(\beta,e)|_{E_2}),  \lb{7.6}\eea
where the left hand sides are the Morse index and nullity of the operator $A(\beta,e)$
on the space $\ol{D}(-1,2\pi)$, i.e., the $-1$ index and nullity of $A(\beta,e)$;
on the right hand sides of (\ref{7.5})-(\ref{7.6}), we denote by
$\phi(A(\beta,e)|_{E_i})$ and $\nu(A(\beta,e)|_{E_i})$ the usual Morse index
and nullity of the operator $A(\beta,e)|_{E_i}$ on the space $E_i$.

By (\ref{4.10}), we have $\phi_{-1}(A(9,e))=0$ and $\nu_{-1}(A(9,e))=0$. Because all
the terms in both sides of (\ref{7.5}) and (\ref{7.6}) are the Morse indices and nullities
which are nonnegative integers, we have
\be \phi(A(9,e)|_{E_1}) = \phi(A(9,e)|_{E_2}) = 0, \qquad
    \nu(A(9,e)|_{E_1}) = \nu(A(9,e)|_{E_2}) = 0.    \lb{7.7}\ee
This shows that $A(9,e)|_{E_i}$ with $i=1$ or $2$ is positive definite.

Since the operator $S(t)$ commutes also with the operator $g$, similarly we have
\be  S(t) = S(t)|_{E_1}\oplus S(t)|_{E_2}. \lb{7.8}\ee
So, for $i=1, 2$, we obtain
\be   A(\beta,e)|_{E_i}
\;=\; A(9,e)|_{E_i}+\left(\frac{\sqrt{9-\beta}}{2(1+e\cos t)}S(t)\right)|_{E_i}
\;=\; \sqrt{9-\beta}\left(\frac{A(9,e)|_{E_i}}{\sqrt{9-\beta}}
             + \frac{S(t)|_{E_i}}{2(1+e\cos t)}\right). \lb{7.9}\ee
Now we want to compute the Morse index of $A(0,e)|_{E_i}$ for $i=1,2$. By (\ref{3.7})
and (\ref{7.5}), we have
\be  \phi(A(0,e)|_{E_1})+\phi(A(0,e)|_{E_2})=2, \qquad \forall\, e\in (-1,1). \lb{7.10}\ee
So the possible value of $\phi(A(0,e)|_{E_i})$ can only be $0,1$ and $2$. By (\ref{3.8})
and (\ref{7.6}), we obtain
\be  \nu(A(0,e)|_{E_1})=\nu(A(0,e)|_{E_2})=0, \qquad\forall\, e\in (-1,1). \lb{7.11}\ee
By the property of Morse index, we have
\be  \phi(A(0,e)|_{E_i})=\phi(A(0,0)|_{E_i}), \qquad \forall\, e\in (-1,1), \;i=1,2. \lb{7.12}\ee
From (\ref{3.43}) and (\ref{7.6}), we obtain
\be  \nu(A(3/4,0)|_{E_1})+\nu(A(3/4,0)|_{E_2})=\nu_{-1}(A(3/4,0))=2.  \lb{7.13}\ee

From the fact that $\xi E_1=E_2$ and $\xi E_2=E_1$ with $\xi$ defined as in (\ref{6.8}), we
have
$$  \nu(A(3/4,0)|_{E_1})=\nu(A(3/4,0)|_{E_2}).  $$
Then we obtain
\be   \nu(A(3/4,0)|_{E_i})=1, \qquad {\rm for}\;\;i=1, 2.   \lb{7.14}\ee

By (\ref{7.9}), for any fixed $e\in (-1,1)$ and $i=1, 2$,
$\frac{A(\beta,e)|_{E_i}}{\sqrt{9-\beta}}$ is increasing with respect to $\beta$ as
proved before. It has the same Morse index and nullity as those of $A(\beta,0)|_{E_i}$.
So we get
\be \phi(A(\beta,0)|_{E_i})=\left\{\matrix{1, &\quad {\rm if}\;\; 0\leq\beta<3/4,\cr
                          0, &\quad {\rm if}\;\; \beta\geq 3/4 ,\cr}\right. \lb{7.15}\ee
This shows that for $-1<e<1$ and $i=1,2$, by (\ref{7.12}) we obtain
\be  \phi(A(0,e)|_{E_i})=1.  \lb{7.16}\ee

By the same idea as in the proof of Theorem \ref{T6.3}, we get

\begin{proposition}\label{P7.1} The $\omega=-1$ degeneracy curve $(\beta_i(e,-1),e)$ is
precisely the degeneracy curve of $A(\beta,e)|_{E_i}$ for $i=1$ or $2$.
\end{proposition}

Following the results of \cite{R1}, \cite{MSS} and \cite{MS}, we know that the curves
$(\beta_1(e,-1),e)$ and $(\beta_2(e,-1),e)$ intersect transversely at the point $(3/4,0)$.
By symmetries of the $\omega$-index gap curves, we have

\begin{theorem}\label{T7.2} $\beta_1(e,-1)=\beta_2(-e,-1)$ holds for all $e\in (-1,1)$.
\end{theorem}

\begin{remark}\label{R7.3} We can also compute $\phi(A(0,0)|_{E_1})$ via the relation
between Maslov-type index and Morse index. Let $V_1=\{(0,x,y,0)\;|\;x,y\in\R\}$ and
$V_2=\{(x,0,0,y)\;|\;x,y\in\R\}$. Then both of them are Lagrangian subspaces of the phase
space $\R^4$ with standard symplectic structure. From Theorem 1.2 of \cite{HS1}, we have
\be \phi(A(0,0)|_{E_1})=\mu(V_2,\gamma_{0,0}(t)V_1), \lb{7.17}\ee
where the right hand side is the Maslov index for paths of Lagrangian subspaces. Note that
for $(\beta,e)=(0,0)$, by (\ref{2.19}), $B(t)\equiv B$ is a constant matrix. Then
$\gamma_{0,0}(t)=\exp(JBt)$ and its Maslov-type index can be computed explicitly as we did
in \S 3.3.
\end{remark}

\subsection{$-1$ degeneracy curve bifurcations from $(3/4,0)$ as $e$ leaves from $0$}

By (\ref{3.22}), $-1$ is a double eigenvalue of the matrix $\ga_{3/4,0}(2\pi)$. As studied by
Roberts (cf. p.212 in \cite{R1}) and Meyer-Schmidt (cf. Section 3 of \cite{MS}), there are two
period doubling curves bifurcating out from $(3/4,0)$ when $e>0$ is sufficiently small. In
\cite{MS}, the tangent directions of these two curves are computed. Note that these two
curves are precisely the $-1$ index gap curves found by our Theorem \ref{T1.2}. For reader's
conveniences, here we give a simple proof on these two tangent directions based on our above
studies.

\begin{proposition}\label{P7.4} The tangent directions of the two curves $\Ga_s$ and $\Ga_m$
bifurcating from $(3/4,0)$ when $e>0$ is small are given by
$$   \beta_s'(e)|_{e=0} = -\frac{\sqrt{33}}{4},
     \qquad \beta_m'(e)|_{e=0} = \frac{\sqrt{33}}{4}.  $$
\end{proposition}

{\bf Proof.} To compute the slope of $-1$ degeneracy curve bifurcating out from $(\beta,e)=(3/4,0)$,
let $(\beta(e),e)$ be one of the curve (say, the $E_1$ degeneracy curve) with $e\in (-\ep, \ep)$
for some small $\ep>0$, and $x_e\in E_1$ be the corresponding eigenvector, that is
\be   A(\beta(e),e)x_e=0.   \lb{7.18}\ee
Here the space $E_1$ is defined in (\ref{7.3}) above. Thus there holds
\be  \<A(\beta(e),e)x_e,x_e\>=0.   \lb{7.19}\ee
Then by direct computations, $\ker(A(3/4,0))\cap E_1$ is generated by $x_0=R(t)z(t)$ with
\be    z(t)=(\frac{7-\sqrt{33}}{4}\sin(t/2),\cos(t/2))^T.    \lb{7.20}\ee
Differentiating both sides of (\ref{7.19}) with respect to $e$ yields
$$ \bb'(e)\<\frac{\pt}{\pt \bb}A(\bb(e),e)x_e,x_e\> + (\<\frac{\pt}{\pt e}A(\bb(e),e)x_e,x_e\>
       + 2\<A(\bb(e),e)x_e,x'_e\> = 0,  $$
where $\bb'(e)$ and $x'_e$ denote the derivatives with respect to $e$. Then evaluating both
sides at $e=0$ yields
\be  \bb'(0)\<\frac{\pt}{\pt \bb}A(3/4,0)x_0,x_0\> + \<\frac{\pt}{\pt e}A(3/4,0)x_0,x_0\> = 0. \lb{7.21}\ee
Then by the definition (\ref{2.29}) of $A(\bb,e)$ we have
\bea
\left.\frac{\pt}{\pt\bb}A(\bb,e)\right|_{(\bb,e)=(3/4,0)}
    &=& \left.R(t)\frac{\pt}{\pt\bb}K_{\bb,e}(t)\right|_{(\bb,e)=(3/4,0)}R(t)^T,  \lb{7.22}\\
\left.\frac{\pt}{\pt e}A(\bb,e)\right|_{(\bb,e)=(3/4,0)}
    &=& \left.R(t)\frac{\pt}{\pt e}K_{\bb,e}(t)\right|_{(\bb,e)=(3/4,0)}R(t)^T,  \lb{7.23}\eea
where $R(t)$ is given in \S 2.1. By direct computations from the definition of $K_{\bb,e}(t)$ in
(\ref{2.20}), we obtain
\bea
&& \frac{\pt}{\pt\bb}K_{\bb,e}(t)\left|_{(\bb,e)=(3/4,0)}
       = \frac{1}{2\sqrt{33}}\left(\matrix{-1 & 0\cr
                                          0 &  1\cr}\right),\right.   \lb{7.24}\\
&& \frac{\pt}{\pt e}K_{\bb,e}(t)\left|_{(\bb,e)=(3/4,0)}
       = \frac{-\cos t}{4}\left(\matrix{6+\sqrt{33} & 0 \cr
                                        0 & 6-\sqrt{33}\cr}\right).\right.   \lb{7.25}\eea
Therefore from (\ref{7.20}) and (\ref{7.22})-(\ref{7.25}) we have
\bea  \<\frac{\pt}{\pt\bb}A(3/4,0)x_0,x_0\>
&=& \<\frac{\pt}{\pt\bb}K_{3/4,0}z,z\>    \nn\\
&=& \int_0^\pi[\frac{1}{2\sqrt{33}}\cos^2(t/2)
                -\frac{1}{2\sqrt{33}}((\frac{7-\sqrt{33}}{4})^2\sin^2(t/2)]dt  \nn\\
&=& \frac{\pi}{4\sqrt{33}}(1-(\frac{7-\sqrt{33}}{4})^2),  \lb{7.26}\eea
and
\bea  \<\frac{\pt}{\pt e}A(3/4,0)x_0,x_0\>
&=& \<\frac{\pt}{\pt e}K_{3/4,0}z,z\>   \nn\\
&=& \int_0^\pi[\frac{1}{4}(6+\sqrt{33})(\frac{7-\sqrt{33}}{4})^2\cos(t)\sin^2(t/2)
                +\frac{1}{4}(6-\sqrt{33})\cos(t)\cos^2(t/2)]dt  \nn\\
&=& \frac{\pi}{16}(6-\sqrt{33}-(6+\sqrt{33})(\frac{7-\sqrt{33}}{4})^2).  \lb{7.27}\eea
Therefore by (\ref{7.21}) and (\ref{7.26})-(\ref{7.27}) we obtain
\be  \beta'(0) = \frac{\sqrt{33}}{4}.  \lb{7.28}\ee
By the above Theorem \ref{T7.2}, the two $-1$ degenerate curves are symmetric with respect to
the $e=0$ axis. Therefore the claim of the Proposition \ref{P7.4} follows from (\ref{7.28}). \hb

\setcounter{equation}{0}
\section{Study on the non-hyperbolic regions}\label{sec:8}

Now we give proofs of the first halves of our main theorems \ref{T1.6} and \ref{T1.2}.

{\bf Proof of (i) of Theorem \ref{T1.6}.} It follows from Theorems
\ref{T6.4} and \ref{T7.2}. \hb

{\bf Proof of the first half of Theorem \ref{T1.2}.} Here we give proofs for items (i)-(iii) and
(ix)-(x) of this theorem.

{\bf (i)} By Theorem \ref{T1.5}, we have got the existence of
two curves defined by $(\beta_i(e),e)$ and $\lim_{e\to 1}\beta_i(e)=0$ with $i=1$ and $2$
such that $\gamma_{\bb,e}(2\pi)$ is degenerate with respect to $\omega=-1$ only on them. Note that
here these two curves may coincide somewhere. Specially we define
\be   0 < \bb_s(e)\equiv \min\{\bb_1(e), \bb_2(e)\} \le \bb_m(e) \equiv \max\{\bb_1(e), \bb_2(e)\}
           < 9,   \qquad {\rm for}\quad e\in [0,1).    \lb{8.1}\ee
Thus (i) is proved.

{\bf (ii)} By the studies in \S 3.3, the only $-1$ degenerate point in the $(\bb,e)$ segment
$[0,9]\times\{0\}$ is $(\bb,e)=(3/4,0)$, which is a $2$-fold $-1$ degenerate point, and there hold
\bea
i_{-1}(\ga_{\bb,0})=2, &\qquad& \nu_{-1}(\ga_{\bb,0})=0,  \qquad {\rm for}\;\;\bb\in [0,3/4),   \lb{8.2}\\
i_{-1}(\ga_{3/4,0})=0, &\qquad& \nu_{-1}(\ga_{3/4,0})=2,    \lb{8.3}\\
i_{-1}(\ga_{\bb,0})=0, &\qquad& \nu_{-1}(\ga_{\bb,0})=0,  \qquad {\rm for}\;\;\bb\in (3/4,9].   \lb{8.4}\eea
Therefore
\be   \bb_i(0)=3/4,  \qquad {\rm for}\quad i=1 \;\;{\rm and}\;\;2. \lb{8.5}\ee

By Meyer and Schmidt in \cite{MS} or our Proposition \ref{P7.4}, the two $-1$ degeneracy curves
bifurcating out from $(\bb,e)=(3/4,0)$ when $e>0$ is sufficiently small must coincide with our curves
$\Ga_s$ and $\Ga_m$ respectively. Because these two curves bifurcate out from $(3/4,0)$ in different angles
with tangents $-\sqrt{33}/4$ and $\sqrt{33}/4$ respectively when $e>0$ is small, they are different from
each other near $(3/4,0)$. By our Theorem \ref{T6.4}, these two curves $\Ga_s$ and $\Ga_m$ are real analytic
with respect to $e$. Therefore they are different curves and their intersection points including the point
$(3/4,0)$ can only be isolated.

By Theorems \ref{T6.3} and \ref{T1.5}, these two curves $\Ga_s$ and $\Ga_m$ must tend to the segment
$[0,9]\times \{1\}$ from $(3/4,0)$ as $e$ increases from $0$ and tends to $1$. By the proof of our Theorem
\ref{T1.8} in \S 8 below, for each $e\in [0,1)$ the function $\bb_k(e)$ defined by (\ref{1.5}) satisfies
$0 < \bb_s(e) \le \bb_m(e) \le \bb_k(e) < 9$, and $\lim_{e\to 1}\bb_k(e) = 0$. Therefore the two curves
$\Ga_s$ and $\Ga_m$ must tend to $(0,1)$ as $e\to 1$.

{\bf (iii)} By our studies on the segments $\{0\}\times [0,1)$ and $\{9\}\times [0,1)$ in \S 3 and the
definitions of $\bb_s(e)$ and $\bb_m(e)$, the index $i_{-1}(\ga_{\bb,e})$ must take the claimed values
$2$, $1$, and $0$ in (\ref{1.7}) respectively when $\bb\in [0,9]\bs \{\bb_s(e), \bb_m(e)\}$ for each
$e\in [0,1)$. Note that when $\bb=\bb_s(e)$ or $\bb_m(e)$, the $-1$ index claim (\ref{1.7}) follows from
(i) of Proposition \ref{P6.1}. The last claim in (iii) follows from Proposition \ref{P6.1}.

{\bf (ix)-(x)} By the Bott-type formula (Theorem 9.2.1 in p.199 of \cite{Lon4}), we obtain
$$ i_1(\ga_{\bb,e}^k) = \sum_{\om^k=1}i_{\om}(\ga_{\bb,e}), \qquad \forall\;k\in\N, $$
and
$$  \nu_1(\ga_{\bb,e}^2) = \nu_1(\ga_{\bb,e}) + \nu_{-1}(\ga_{\bb,e}) = \nu_{-1}(\ga_{\bb,e}),  $$
where the last equality follows from Theorem \ref{T1.1}.

Therefore by (\ref{1.7}), (\ref{3.7}) and (\ref{3.8}), for $e\in [0,1)$ we obtain
\be \phi^2=\left\{\matrix{4, &\quad {\rm if}\;\; 0 \le \beta < \bb_s(e), \cr
                          3, &\quad {\rm if}\;\; \bb_s(e) \le \bb < \bb_m(e), \cr
                          2, &\quad {\rm if}\;\; \bb_m(e) \le \bb \leq 9. \cr}\right. \lb{8.6}\ee
Thus for $(\bb,e)\in (0,9]\times [0,1)$, the matrix $\ga_{\bb,e}(4\pi)=\ga^2_{\bb,e}(2\pi)$ is
non-degenerate with respect to the eigenvalue $1$, whenever $(\bb,e)\not\in \Ga_s\cup \Ga_m$.
Therefore by (\ref{8.6}) we can apply Theorem \ref{T2.5} (i.e., Theorem 1.2 of \cite{HS2}) to
get (ix) and (x). The rest parts of Theorem \ref{T1.2} will be proved in the next section. \hb

\setcounter{equation}{0}
\section{Study on the hyperbolic region}\label{sec:9}

In this section we study the hyperbolic region of $\ga_{\bb,e}(2\pi)$ in the rectangle $[0,9]\times [0,1)$.
By the first halves of Theorems \ref{T1.2} and \ref{T1.6} proved in the \S 8, the function
$\bb_k(e)$ defined by (\ref{1.5}) satisfies
\be  \bb_m(e) \le \bb_k(e), \qquad \forall\;e\in [0,1). \lb{9.1}\ee
We have the following further results.

\begin{lemma}\label{L9.1}
(i) If $0\le\bb_1<\bb_2\le 9$ and $\ga_{\bb_1,e}(2\pi)$ is hyperbolic, so does $\ga_{\bb_2,e}(2\pi)$.
Consequently, the hyperbolic region of $\ga_{\bb,e}(2\pi)$ in $[0,9]\times [0,1)$ is connected.

(ii) For any fixed $e\in [0,1)$, every matrix $\ga_{\bb,e}(2\pi)$ is hyperbolic if $\bb_k(e) < \bb \le 9$
for $\bb_k(e)$ defined by (\ref{1.5}). Thus (\ref{1.8}) holds and $\Ga_k$ is the boundary set of this
hyperbolic region.

(iii) We have
\be  \sum_{\bb\in [0,\bb_k(e)]}\nu_{\om}(\ga_{\bb,e}(2\pi)) = 2, \qquad \forall\;\om\in\U\bs\{1\}. \lb{9.2}\ee
\end{lemma}

{\bf Proof.} (i) By Lemma 4.4, for any fixed $\om\in \U$ and $e\in [0,1)$, the
operator $\bar{A}(\bb,e)=\frac{1}{\sqrt{9-\bb}}A(\beta,e)$ is self-adjoint on $\bar{D}(\om,2\pi)$
and increasing with respect to $\bb$ in the sense that
\be  \bar{A}(\bb_1,e)<\bar{A}(\bb_2,e), \qquad {\rm if}\quad \bb_1 < \bb_2.   \lb{9.3}\ee
Suppose $\ga_{\bb_1,e}(2\pi)$ is hyperbolic. This implies that $A(\bb_1,e)$ is non-degenerate on
$\bar{D}(\om,2\pi)$ for every $\om\in\U$. By (ix)-(x) of Theorem \ref{T1.2}, it also implies $\bb_m(e)<\bb_1$.
Thus by (\ref{2.30}), (\ref{3.7}), Corollary \ref{C4.5}, and Theorem \ref{T1.1}, the $\om$-index
$\phi_{\om}(\bar{A}(\bb_1,e))=0$ for all $\om\in\U$. Then $\bar{A}(\bb_1,e)$ is positive definite on
$\bar{D}(\om,2\pi)$ for every $\om\in\U$. Therefore by (\ref{9.3}) the operator
$\bar{A}(\bb_2,e)$ is positive definite too, and then is non-degenerate on $\bar{D}(\om,2\pi)$ for all
$\om\in\U$. Therefore $\ga_{\bb_2,e}(2\pi)$ must be hyperbolic and so does $\ga_{\bb,e}(2\pi)$ for
all $\bb\in [\bb_1,9)$.

Recall that along the segment $\{9\}\times [0,1)$ the matrix $\ga_{9,e}(2\pi)$ is hyperbolic by our
Proposition \ref{P1.4}. Therefore the hyperbolic region of $\ga_{\bb,e}(2\pi)$ is connected in
$[0,9]\times [0,1)$.

(ii) By the definition of $\bb_k(e)$, there exists a sequence $\{\bb_i\}_{i\in\N}$ satisfying
$\bb_i > \bb_k(e)$, $\bb_i\to \bb_k(e)$, and $\ga_{\bb_i,e}(2\pi)$ is hyperbolic. Therefore
$\ga_{\bb,e}(2\pi)$ is hyperbolic for every $\bb\in (\bb_k(e),9]$ by (i). Then (\ref{1.8}) holds
and $\Ga_k$ is the envelop curve of this hyperbolic region.

Now (iii) follows from (ii) and Proposition \ref{P6.1}, and the proof is complete. \hb

\begin{corollary}\label{C9.2} For every $e\in [0,1)$, we have
\be \sum_{\bb\in (0,\bb_m(e)]}\nu_{-1}(\ga_{\bb,e}(2\pi)) = 2 \quad {\it and}\quad
    \sum_{\bb\in (\bb_m(e),9]}\nu_{-1}(\ga_{\bb,e}(2\pi)) = 0. \lb{9.4}\ee
\end{corollary}

{\bf Proof.} Fix an $e\in [0,1)$. If $\bb_s(e)<\bb_m(e)$, then we obtain
$$ \sum_{\bb\in (0,\bb_m(e)]}\nu_{-1}(\ga_{\bb,e}(2\pi))
    \ge \nu_{-1}(\ga_{\bb_s(e),e}(2\pi)) + \nu_{-1}(\ga_{\bb_m(e),e}(2\pi)) \ge 2. $$
Thus (\ref{9.4}) follows from (ii) of Proposition \ref{P6.1}.

If $\bb_s(e)=\bb_m(e)$, then by (i) of Proposition \ref{P6.1} we obtain
$\nu_{-1}(\ga_{\bb_m(e),e}(2\pi)) = 2$. Therefore we have
$$ \sum_{\bb\in (0,\bb_m(e)]}\nu_{-1}(\ga_{\bb,e}(2\pi))
    \ge \nu_{-1}(\ga_{\bb_m(e),e}(2\pi)) = 2. $$
Thus (\ref{9.4}) follows also from (ii) of Proposition \ref{P6.1}. \hb

Now we can give the

{\bf Proof of the second half of Theorem \ref{T1.2}.} Here we give the proofs for the items
(iv)-(viii) and (xi) of this theorem.

Note that Claims {\bf (iv)} and {\bf (v)} of the theorem follow from (\ref{9.1}) and Lemma
\ref{L9.1}.

{\bf (vi)} In fact, if the function $\bb_k(e)$ is not continuous in $e\in [0,1)$, then there
exist some $\hat{e}\in [0,1)$, a sequence $\{e_i\,|\,i\in\N\}\subset [0,1)\bs\{\hat{e}\}$
and $\bb_0\in [0,9]$ such that
\be   \bb_k(e_i)\to \bb_0 \not= \bb_k(\hat{e}) \quad {\rm and}\quad e_i\to \hat{e}
             \qquad {\rm as}\;\;i\to +\infty .   \lb{9.5}\ee
We continue in two cases according to the sign of the difference $\bb_0-\bb_k(\hat{e})$.

By the definition of $\bb_k(e_i)$ we have $\sg(\ga_{\bb_k(e_i),e_i}(2\pi))\cap\U \not= \emptyset$
for every $e_i$. By the continuity of eigenvalues of $\ga_{\bb_k(e_i),e_i}(2\pi)$ in $i$ and (\ref{9.5}),
we obtain
$$  \sg(\ga_{\bb_0,\hat{e}}(2\pi))\cap\U \not= \emptyset.   $$
Then by Lemma \ref{L9.1}, this would yield a contradiction if $\bb_0>\bb_k(\hat{e})$.

Now suppose $\bb_0 < \bb_k(\hat{e})$. By Lemma \ref{L9.1} for all $i\ge 1$ we have
\be  \sg(\ga_{\bb,e_i}(2\pi))\cap\U \;=\; \emptyset, \qquad \forall\;\bb\in (\bb_k(e_i),9]. \lb{9.6}\ee
Then by the continuity of $\bb_m(e)$ in $e$, (\ref{9.6}) and the definition of $\bb_0$, we obtain
$$  \bb_m(\hat{e}) \le \bb_0 < \bb_k(\hat{e}).  $$
Let $\om_0\in \sg(\ga_{\bb_k(\hat{e}),\hat{e}}(2\pi))\cap \U$, which exists by the definition of
$\bb_k(\hat{e})$.

Let $L=\{(\bb,\hat{e})\;|\;\bb\in (\bb_k(\hat{e}),9]\}$, $V=\{(9,e)\;|\;e\in [0,1)\}$, and
$L_i=\{(\bb,e_i)\;|\;\bb\in (\bb_k(e_i),9]\}$ for all $i\ge 1$. Note that by (\ref{3.9}), (\ref{4.10}),
Corollary \ref{C4.5}, Proposition \ref{P6.1}, Lemma \ref{L9.1}, and the definitions of $\bb_k(e_i)$
and $\bb_k(\hat{e})$, we obtain
\be  i_{\om_0}(\ga_{\bb,e}) = \nu_{\om_0}(\ga_{\bb,e}) = 0, \qquad
       \forall\;(\bb,e)\;\in\; L\cup V\cup \bigcup_{i\,\ge 1}L_i.      \lb{9.7}\ee
Specially we have
$$  i_{\om_0}(\ga_{\bb_k(\hat{e}),\hat{e}}) = 0 \quad {\rm and}\quad
                        \nu_{\om_0}(\ga_{\bb_k(\hat{e}),\hat{e}}) \ge 1. $$
Therefore by Proposition \ref{P6.1} and the definition of $\om_0$, there exists
$\hat{\bb}\in (\bb_0,\bb_k(\hat{e}))$ sufficiently close to $\bb_k(\hat{e})$ such that
\be  i_{\om_0}(\ga_{\hat{\bb},\hat{e}}) = i_{\om_0}(\ga_{\bb_k(\hat{e}),\hat{e}})
        + \nu_{\om_0}(\ga_{\bb_k(\hat{e}),\hat{e}}(2\pi)) \ge 1.      \lb{9.8}\ee
This estimate (\ref{9.8}) in fact holds for all $\bb\in [\hat{\bb},\bb_k(\hat{e}))$ too. Note that
$(\hat{\bb},\hat{e})$ is an accumulation point of $\cup_{i\;\ge 1}L_i$. Consequently for each $i\ge 1$
there exists $(\bb_i,e_i)\in L_i$ such that $\ga_{\bb_i,e_i}\in \P_{2\pi}(4)$ is $\om_0$ non-degenerate,
$\bb_i\to \hat{\bb}$ in $\R$, and $\ga_{\bb_i,e_i}\to \ga_{\hat{\bb},\hat{e}}$ in $\P_{2\pi}(4)$ as
$i\to \infty$. Therefore by (\ref{9.7}), (\ref{9.8}), the Definition 5.4.2 of the $\om_0$-index of
$\om_0$-degenerate path $\ga_{\hat{\bb},\hat{e}}$ on p.129 and Theorem 6.1.8 on p.142 of \cite{Lon4},
we obtain the following contradiction
$$  1 \le i_{\om_0}(\ga_{\hat{\bb},\hat{e}}) \le i_{\om_0}(\ga_{\bb_i,e_i}) = 0,  $$
for $i\ge 1$ large enough. Thus the continuity of $\bb_k(e)$ in $e\in [0,1)$ is proved.

{\bf (vii)} To prove the claim $\lim_{e\to 1}\bb_k(e)=0$, we argue by contradiction, and suppose that
there exist $e_i\to 1$ as $i\to +\infty$, $\bb_0>0$, such that $\lim_{i\to \infty}\bb_k(e_i)=\bb_0$.
Then at least one of the following two cases must occur:

(A) {\it There exists a subsequence $\{\hat{e}_i\}$ of $\{e_i\}$ such that
$\bb_k(\hat{e}_{i+1}) \le \bb_k(\hat{e}_i)$ for all $i\in\N$;}

(B) {\it There exists a subsequence $\{\hat{e}_i\}$ of $\{e_i\}$ such that
$\bb_k(\hat{e}_i) \le \bb_k(\hat{e}_{i+1})$ for all $i\in\N$.}

If Case (A) happens, for this $\bb_0$ by Theorem \ref{T1.7} there exists $e_0>0$ sufficiently close to
$1$ such that $\ga_{\bb_0,e}(2\pi)$ is hyperbolic for all $e\in [e_0,1)$. Then $\ga_{\bb,e}(2\pi)$ is
hyperbolic for all $(\bb,e)$ in the region $[\bb_0,9]\times [e_0,1)$ by Lemma \ref{L9.1}. But by the
monotonicity of Case (A) we obtain
$$  \bb_0 \le \bb_k(\hat{e}_{i+m}) \le \bb_k(\hat{e}_i) \qquad \forall\; m\in\N. $$
Therefore $(\bb_k(\hat{e}_{i+m}),\hat{e}_{i+m})$ will get into this region for sufficiently large $m\ge 1$,
which contracts to the definition of $\bb_k(\hat{e}_{i+m})$ in (\ref{1.5}).

If Case (B) happens, fix a subindex $i$, by Theorem \ref{T1.7} and the same argument as in Case (A) there
exists an $e_0>0$ sufficiently close to $1$ such that $\ga_{\bb,e}(2\pi)$ is hyperbolic for all $(\bb,e)$
in the region $[\bb_k(\hat{e}_{i}), 9]\times [e_0,1)$. Then by the monotonicity of Case (B) we obtain
$$  \bb_k(\hat{e}_i) \le \bb_k(\hat{e}_{i+m}) \le \bb_0 \qquad \forall\; m\in\N. $$
Therefore $(\bb_k(\hat{e}_{i+m}),\hat{e}_{i+m})$ will get into this region for sufficiently large $m\ge 1$,
which contracts to the definition of $\bb_k(\hat{e}_{i+m})$ in (\ref{1.5}). Thus (vii) holds.

{\bf (viii)} By our study in \S 3.3, we have $(1,0)\,\in\,\Ga_k\bs\Ga_m$. Thus there exists an $\td{e}\in (0,1]$
such that $\bb_k(e)>\bb_m(e)$ for all $e\in [0,\td{e})$. Therefore $\Ga_k$ is different from $\Ga_m$
when $e\in [0,\td{e})$.

{\bf (xi)} Let $e_0\in [0,1)$ and $\bb_m(e_0)<\bb_0\le \bb_k(e_0)$. Then $M\equiv\ga_{\bb_0,e_0}(2\pi)$
is not hyperbolic by Lemma \ref{L9.1} and thus at least one pair of its eigenvalues is on $\U$.
Note also that no eigenvalues of $M$ can be $\pm 1$ by Theorem \ref{T1.1} and Corollary \ref{C9.2}. Write
\be  \sg(M)=\{\lm_1(\bb_0), \lm_1(\bb_0)^{-1}, \lm_2(\bb_0), \lm_2(\bb_0)^{-1}\}.  \lb{9.9}\ee
Thus we can assume $\lm_1\equiv\lm_1(\bb_0)\in \U\bs\R$ and the other pair of eigenvalues satisfy
$\lm_2\equiv\lm_2(\bb_0) \in (\U\cup\R)\bs\{\pm 1, 0\}$.

{\it Claim.} $\lm_2(\bb_0)\in \U\bs\R$.

In fact, if not, we assume $\lm_2(\bb_0)\in \R\bs\{\pm 1, 0\}$.

In this case, $M$ has normal form $R(\th)\dm D(\lm_2) \in \Om^0(M)$ for some $\th\in (0,\pi)\cup (\pi,2\pi)$.
Thus by (\ref{3.6}), Theorem \ref{T1.1} and (iii) of our Theorem \ref{T1.2}, we obtain the following
contradiction:
\bea 0
&=& i_{-1}(\ga_{\bb_0,e_0})   \nn\\
&=& i_{1}(\ga_{\bb_0,e_0}) + S_{M}^+(1) - S_{M}^-(e^{\pm\sqrt{-1}\th}) + S_{M}^+(e^{\pm\sqrt{-1}\th}) - S_{M}^-(-1)  \nn\\
&=& 0 + 0 - S_{R(\th)}^-(e^{\pm\sqrt{-1}\th}) + S_{R(\th)}^+(e^{\pm\sqrt{-1}\th}) - 0  \nn\\
&=& \pm 1,  \nn\eea
where the last equality follows from Lemma 9.1.6 on p.192 and $\<5\>$ of List 9.1.12 on p.198 of \cite{Lon4}.
Therefore the claim is proved.

Now from $\lm_1(\bb_0)$ and $\lm_2(\bb_0)\;\in\;\U\bs\R$, the matrix $M$ has basic normal form
$R(\th_1)\dm R(\th_2) \in \Om^0(M)$ for some $\th_1$ and $\th_2\in (0,\pi)\cup (\pi,2\pi)$. Then by
the study in Section 9.1 of \cite{Lon4}, we obtain
\bea 0
&=& i_{-1}(\ga_{\bb_0,e_0})    \nn\\
&=& i_1(\ga_{\bb_0,e_0}) + S_{M}^+(1) - S_{R(\th_1)}^-(e^{\pm \sqrt{-1}\th_1}) + S_{R(\th_1)}^+(e^{\pm \sqrt{-1}\th_1}) \nn\\
& & \hskip 3 cm   - S_{R(\th_2)}^-(e^{\pm \sqrt{-1}\th_2}) + S_{R(\th_2)}^+(e^{\pm \sqrt{-1}\th_2}) - S_{M}^-(-1)  \nn\\
&=&  - S_{R(\th_1)}^-(e^{\pm \sqrt{-1}\th_1}) + S_{R(\th_1)}^+(e^{\pm \sqrt{-1}\th_1})
         - S_{R(\th_2)}^-(e^{\pm \sqrt{-1}\th_2}) + S_{R(\th_2)}^+(e^{\pm \sqrt{-1}\th_2}).   \lb{9.10}\eea
By Lemma 9.1.6 on p.192 and $\<5\>$ of List 9.1.12 on p.198 of \cite{Lon4} again, the right hand side of
(\ref{9.10}) would be $\pm 2$, if both $\th_1$ and $\th_2$ are located in only one interval of $(0,\pi)$ and
$(\pi,2\pi)$. Thus we must have $\th_1\in (0,\pi)$ and $\th_2\in (\pi,2\pi)$. Let $\om=\exp(\sqrt{-1}\th_1)$.

If $2\pi-\th_2=\th_1$, we then obtain
$$ \sum_{0\le \bb\le \bb_0}\nu_{\om}(\ga_{\bb,e_0})
    \ge \sum_{0\le \bb\le \bb_m(e_0)}\nu_{\om}(\ga_{\bb,e_0}) + \nu_{\om}(\ga_{\bb_0,e_0})
    \ge 1 + 2.  $$
It contradicts to Lemma \ref{L9.1} and proves $2\pi-\th_2 \not= \th_1$.

By the study in Section 9.1 of \cite{Lon4} again, if $2\pi-\th_2>\th_1$, for $\om=\exp(\sqrt{-1}\th_1)$ we obtain
$$ 0 \le i_{\om}(\ga_{\bb_0,e_0})
= i_1(\ga_{\bb_0,e_0}) + S_{M}^+(1) - S_{R(\th_1)}^-(e^{\sqrt{-1}\th_1})
= - S_{R(\th_1)}^-(e^{\sqrt{-1}\th_1}) = -1. $$
This contradiction proves that the only possible case is $2\pi-\th_2<\th_1$.

The proof of Theorem \ref{T1.2} is complete. \hb

Now we give

{\bf The proof of (ii) of Theorem \ref{T1.6}.} By (v) of Theorem \ref{T1.2}, the curve $\Ga_k$ for $e\in [0,1)$
is the boundary curve of the hyperbolic region of $\ga_{\bb,e}(2\pi)$ in the $(\bb,e)$ rectangle
$[0,9]\times [0,1)$. By the definition (\ref{1.5}), the curve $\Ga_k$ is also the envelop curve of the $\om$
degeneracy curves for all $\om\;\in\;\U\bs\{1\}$ from the right hand side of the rectangle $[0,9]\times [0,1)$.
Then (i) of Theorem \ref{T1.6} implies that $\Ga_k$ can be continuously extended into $[0,9]\times (-1,0]$
so that it is symmetric with respect to $[0,9]\times \{0\}$. \hb

The next lemma is useful in the proof of Theorem \ref{T1.8}.

\begin{lemma}\label{L9.3} If $\ga_{\bb_0,e}(2\pi) \approx M_2(-1,c)$ holds or it possesses
the basic normal form $N_1(-1,a)\dm N_1(-1,b)$ for some $(\bb_0,e)\in (0,9)\times [0,1)$ and
$a$, $b\in \R$ and $c\in\R^2$, then $\ga_{\bb,e}(2\pi)$ is hyperbolic for all $\bb \in (\bb_0,9]$.
\end{lemma}

{\bf Proof.} Note that the basic normal form of $M_2(-1,c)$ is either $N_1(-1,\hat{a})\dm N_1(-1,\hat{b})$
or $N_1(-1,\hat{a})\dm D(\lm)$ for some $\hat{a}$, $\hat{b}\in\R$ and $0>\lm\not= -1$. Thus for any
$\om\in\U\bs\{1\}$, by Corollary \ref{C4.5}, (\ref{2.31}), and the study in Section 9.1 of \cite{Lon4},
we obtain
$$  0\le i_{\om}(\ga_{\bb_0,e}) = i_1(\ga_{\bb_0,e}) + S_{M}^+(1) - S_{M}^-(\om) = - S_{M}^-(\om) \le 0, $$
where $M=\ga_{\bb_0,e}(2\pi)$. This proves $i_{\om}(\ga_{\bb_0,e})=0$ for all $\om\in\U$.
Note that $\phi_{\om}(\bar{A}(\bb_0,e))=i_{\om}(\ga_{\bb_0,e})$ and
$\nu_{\om}(\bar{A}(\beta_0,e))=\nu_{\om}(\ga_{\bb_0,e}(2\pi))$ follow from (\ref{2.27}), (\ref{2.30}),
(\ref{4.12}) and (\ref{4.13}).

Now from $\phi_{\om}(\bar{A}(\bb_0,e))=0$ and (ii) of Lemma \ref{L4.4}, we obtain $\bar{A}(\beta,e)>0$ for
all $\bb \in (\bb_0,9]$ on $\ol{D}(\om,2\pi)$ with $\om\in\U$. Therefore
$\nu_{\om}(\ga_{\bb,e}(2\pi)) = \nu_{\om}(\bar{A}(\bb,e)) = 0$ holds for all $\bb\in (\bb_0,9]$ and
$\om\in\U$, and thus the lemma follows. \hb

Now we can give

{\bf Proof of Theorem \ref{T1.8}}. {\bf (i)}  Let $e\in [0,1)$ satisfy $\bb_s(e)<\bb_m(e)$. Then Corollary
\ref{C9.2} implies $\nu_{-1}(\ga_{\bb_s(e),e})=1$. As the limiting case of Cases (ix) and (x) of Theorem
\ref{T1.2}, the matrix $M=\ga_{\bb_s(e),e}(2\pi)$ must have all eigenvalues in $\U$, and possesses its
normal form either $M\approx M_2(-1,c)$ for some $c_2\not=0$, or $M\approx N_1(-1,1)\dm R(\th)$ for some
$\th\in (\pi,2\pi)$, where to get the second case we have used the Figure 2.1.2 on p.50 of \cite{Lon4}
and the fact $N_1(-1,1)\in \Sp(2)^0_{-1,-}$ in that figure.

Note that $M\approx M_2(-1,c)$ for some $c_2\not=0$ can not hold by Lemma \ref{L9.3} and the fact
$\bb_s(e)<\bb_m(e)$. The following is a direct proof of this fact. In this case, its basic
normal form is $N_1(-1,a)\dm D(\lm)$ for some $a\in \{-1, 1\}$ and $0>\lm\not= -1$. Therefore by Theorem
\ref{T1.1}, (iii) of Theorem \ref{T1.2}, and $\<3\>$ and $\<4\>$ of List 9.1.12 on p.198 of \cite{Lon4},
we obtain the following contradiction
$$  1 = i_{-1}(\ga_{\bb_s(e),e}) = i_1(\ga_{\bb_s(e),e}) + S_{M}^+(1) - S_{N_1(-1,a)}^-(-1)
     = - S_{N_1(-1,a)}^-(-1) \le 0. $$

Thus $M\approx N_1(-1,1)\dm R(\th)$ must hold for some $\th\in (\pi,2\pi)$, so $M$ is spectrally
stable and linearly unstable.

{\bf (ii)} Let $e\in [0,1)$ satisfy $\bb_s(e)=\bb_m(e)<\bb_k(e)$. As the limiting case of the Cases (ix)
and (xi) of Theorem \ref{T1.2} and Corollary \ref{C9.2}, the matrix $M=\ga_{\bb_s(e),e}(2\pi)$ must have
basic normal form either $N_1(-1,a)\dm N_1(-1,b)$ for some $a$ and $b\in \{-1, 1\}$, or $-I_2\dm R(\th)$
for some $\th\in (\pi,2\pi)$, where we have used the Figure 2.1.2 on p.50 of \cite{Lon4}. But the first
case is impossible by Lemma \ref{L9.3}. Therefore $M\approx -I_2\dm R(\th)$ holds for some
$\th\in (\pi,2\pi)$, and it is linear stable and not strongly linear stable.

{\bf (iii)} Let $e\in [0,1)$ satisfy $\bb_s(e)<\bb_m(e)<\bb_k(e)$. As the limiting case of Cases
(x) and (xi) of Theorem \ref{T1.2}, the matrix $M=\ga_{\bb_m(e),e}(2\pi)$ must satisfy either
$M\approx N_1(-1,-1)\dm R(\th)$ for some $\th\in (\pi,2\pi)$, or $M\approx M_2(-1,c)$
with $c_2\not=0$, where we have used the Figure 2.1.2 on p.50 of \cite{Lon4} and the fact
$N_1(-1,-1)\in \Sp(2)^0_{-1,+}$ in that figure. Here the second case is also impossible by Lemma \ref{L9.3},
and the conclusion of (iii) follows.

{\bf (iv)}  Let $e\in [0,1)$ satisfy $\bb_s(e)\le\bb_m(e)<\bb_k(e)$. As the limiting case of the Cases
(v) and (xi) of Theorem \ref{T1.2}, the matrix $M\equiv\ga_{\bb_k(e),e}(2\pi)$ must have Krein collision
eigenvalues $\sg(M)=\{\lm_1, \ol{\lm}_1, \lm_2, \ol{\lm}_2\}$ with $\lm_1=\ol{\lm}_2=e^{\sqrt{-1}\th}$
for some $\th\in (0,\pi)\cup (\pi,2\pi)$. Here we have used Theorem \ref{T1.1} and Corollary \ref{C9.2}
to exclude the possibility of eigenvalues $\pm 1$. Therefore for this angle $\th$, the matrix $M$ must
have its normal form $N_2(\om,b)$ for $\om=e^{\sqrt{-1}\th}$ and some $2\times 2$ matrix
$b=\left(\matrix{b_1 & b_2\cr
                 b_3 & b_4\cr}\right)$, which is of the form (25)-(27) by Theorem 1.6.11 on p.34 of
\cite{Lon4}. Because $(I_2\dm (-I_2))^{-1}N_2(e^{\sqrt{-1}\th},b)(I_2\dm (-I_2))=N_2(e^{\sqrt{-1}(2\pi-\th)},\hat{b})$
holds for $\hat{b}=\left(\matrix{b_1 & -b_2\cr
                                 -b_3 & b_4\cr}\right)$, we can always suppose $\th\in (0,\pi)$
without changing the fact $M\approx N_2(\om,b)$.

Note that by (\ref{3.7}), (\ref{3.9}), Corollary \ref{C4.5} and Proposition \ref{P6.1}, we
have $i_{\om}(\ga_{\bb_k(e),e}) = 0$.

Now if $b_2-b_3=0$, by Lemma 1.9.2 on p.43 of \cite{Lon4}, we get $\nu_{\om}(N_2(\om,b))=2$, and then
$N_2(\om,b)$ has basic normal form $R(\th)\dm R(2\pi-\th)$ by the study in Case 4 on p.40 of \cite{Lon4}.
Thus we arrive at the following contradiction
$$  0 \;=\; i_{\om}(\ga_{\bb_k(e),e})
\;=\; i_1(\ga_{\bb_k(e),e}) + S_{M}^+(1) - S_{R(\th)}^-(\om) - S_{R(\th)}^-(\ol{\om})
\;\le\; -1,   $$
by Lemma 9.1.6 on p.192 and $\<5\>$ of List 9.1.12 on p.198 of \cite{Lon4}.

Therefore $b_2-b_3\not= 0$ must hold. Then we obtain
$$  0 \;=\; i_{\om}(\ga_{\bb_k(e),e})
\;=\;  i_1(\ga_{\bb_k(e),e}) + S_{M}^+(1) - S_{N_2(\om,b)}^-(\om)
\;=\; - S_{N_2(\om,b)}^-(\om).    $$
By $\<6\>$ and $\<7\>$ in List 9.1.12 on p.199 of \cite{Lon4}, we obtain that $N_2(\om,b)$
must be trivial as in our discussion in \S 2.1. Then by Theorem 1 of \cite{ZhL1}, the matrix $M$
is spectrally stable and is linearly unstable as claimed.

{\bf (v)} Let $e\in [0,1)$ satisfy $\bb_s(e)<\bb_m(e)=\bb_k(e)$. Note first that $-1$ must be
an eigenvalue of $M=\ga_{\bb_k(e),e}(2\pi)$ with geometric multiplicity $1$ by Corollary \ref{C9.2}.
As the limiting case of Cases (v) and (x) of Theorem \ref{T1.2}, the matrix $M$ must satisfy
either $M\approx M_2(-1,b)$ with $b_1, b_2\in\R$ and $b_2\not=0$, and thus is spectrally stable
and linearly unstable; or $M\approx N_1(-1,a)\dm D(\lm)$ for some $a\in \{-1, 1\}$ and
$-1\not= \lm <0$.

Then in the later case we obtain
$$  0 \;=\; i_{-1}(\ga_{\bb_k(e),e})
\;=\;  i_1(\ga_{\bb_k(e),e}) + S_{M}^+(1) - S_{N_1(-1,a)}^-(-1)
\;=\; - S_{N_1(-1,a)}^-(-1).  $$
Then by $\<3\>$ and $\<4\>$ in List 9.1.12 on p.199 of \cite{Lon4}, we must have $a=1$. This
case can be seen in Figure 2.1.2 on p.50 of \cite{Lon4} with the fact
$N_1(-1,1)\in \Sp(2)^0_{-1,-}$. Thus $M$ is elliptic-hyperbolic (EH) and linearly
unstable.

Note that by the above argument, the matrix $M_2(-1,b)$ has also basic normal form
$N_1(-1,1)\dm D(\lm)$ for some $-1\not= \lm <0$.

{\bf (vi)} Let $e\in [0,1)$ satisfy $\bb_s(e)=\bb_m(e)=\bb_k(e)$. As the limiting case of Cases
(v) and (ix) of Theorem \ref{T1.2}, $-1$ must be the only eigenvalue of $M=\ga_{\bb_k(e),e}(2\pi)$
with $\nu_{-1}(M)=2$ by Corollary \ref{C9.2}. Thus the matrix $M$ must satisfy
$M\approx M_2(-1,c)$ with $c_2=0$ and $\nu_{-1}(M_2(-1,c))=2$ by Subsection 2.1; or
$M\approx N_1(-1,\hat{a})\dm N_1(-1,\hat{b})$ for some $\hat{a}$ and $\hat{b}\in \{-1, 1\}$.
In both cases, $M$ has basic normal form $N_1(-1,a)\dm N_1(-1,b)$ for some
$a$ and $b\in \{-1, 1\}$. Thus we obtain
\bea 0
&=& i_{-1}(\ga_{\bb_k(e),e})    \nn\\
&=& i_1(\ga_{\bb_k(e),e}) + S_{M}^+(1) - S_{N_1(-1,a)}^-(-1) - S_{N_1(-1,b)}^-(-1)   \nn\\
&=& - S_{N_1(-1,a)}^-(-1) - S_{N_1(-1,b)}^-(-1).  \nn\eea
Then by $\<3\>$ and $\<4\>$ in List 9.1.12 on p.199 of \cite{Lon4}, we must have $a=b=1$ similar
to our above study for (v). Therefore it is spectrally stable and linearly unstable as claimed.

The proof is complete. \hb

We give finally the

{\bf Proof of Theorem \ref{T1.9}.} Note first that $\gamma_{\beta,e}$ is analytic in both $\beta$ and
$e$. Since the property of the spectrum for a $4\times 4$ symplectic matrix being complex saddle is an open
condition in $\Sp(4)$, the set $I_e$ in the theorem must be open in $\beta$ for any fixed $e\in [0,1)$.
Thus for any fixed $e\in [0,1)$, it suffices to show $I_e \neq\emptyset$. We argue by contradiction and
suppose $I_e =\emptyset$. Then for any $(\beta,e)\in (0,9]\times [0,1)$, all the eigenvalues of
$\gamma_{\beta,e}$ form two pairs and are located on the union $(\R\bs\{0\})\cup \U$. By our Theorem
\ref{T1.1}, the matrix $\ga_{\bb,e}$ is non-degenerate when $\bb>0$, and then it has no eigenvalue $1$ at
all. But by our Proposition \ref{P1.4}, the matrix $\gamma_{9,e}(2\pi)$ has a pair of double positive eigenvalues
not equal to $1$ for all $e\in [0,1)$. Therefore fix an $e\in [0,1)$. By the continuity of the spectrum in $\bb$,
the matrix $\ga_{\bb,e}(2\pi)$ would have only positive eigenvalues not equal to $1$ for all $\bb\in (0,9]$. This
then contradicts to our Theorem \ref{T1.2}, which yields the existence of the eigenvalue $-1$ of
$\gamma_{\beta,e}(2\pi)$ for some $\beta\in (0,9)$, and completes the proof. \hb

\setcounter{equation}{0}
\section{More observations}\label{sec:10}

Here we describe briefly some more results of \cite{MSS2} of Mart\'{\i}nez, Sam\`{a} and Sim\'{o}
on the Lagrangian triangular homographic solutions in the Newton potential case.

For $e<1$ and close to $1$, the system is HH for any $\beta$ except in a neighborhood of some
critical value which, numerically, appears to be equal to $6$. Some interesting tangencies
are also observed near the corresponding boundaries.

\begin{itemize}
\item[(i)] The tangency at $(\beta,e)=(0,1)$ between the $e$-vertical axis and the
curve which separates the EE and EH domains is of the form $e = 1-C\beta^{\frac{2}{5}}$;

\item[(ii)] The tangency at $(\beta,e)=(0,1)$ between the $e=1$ horizontal line and the
curve which separates the EH and HH domains is of the form $e = 1-C\beta^{4}$;

\item[(iii)] For fixed $\beta\in(0,9)$, the matrix $\gamma_{\beta,e}(2\pi)$
is HH if $1-e>0$ is small enough under a special "non-degenerate" condition, which is defined
in their Lemma 5 on p.663 of \cite{MSS2}, i.e., $d_g\not= 0$ and $e_g\not= 0$ there. It seems
to us that it is not easy to verify this non-degenerate condition, and that the point
$(\beta,e)=(6,1)$ is a possible degenerate point is only checked numerically in \cite{MSS2}.

\item[(iv)] The tangency at $(\beta,e)=(9,0)$ between the $\beta = 9$ vertical line and
the curve which separates the HH and CS domains is of the form $e=C(9-\beta)^{\frac{1}{4}}$.
\end{itemize}

In all the above expressions $C$ denotes suitable constants. Furthermore there is a
point of contact of four different types of domains located approximately at $(1.2091, 0.3145)$.

In the current paper, we have proved the non-degeneracy of the elliptic Lagrangian
triangular solutions. We have also proved the global existences of separation curves
$\Ga_s$ in (i), $\Ga_m$ and $\Ga_k$ in (ii). Our Theorem \ref{T1.7} is related to (iii).

\begin{figure*}[h]
\centering
\includegraphics [width=100mm,height=60mm]{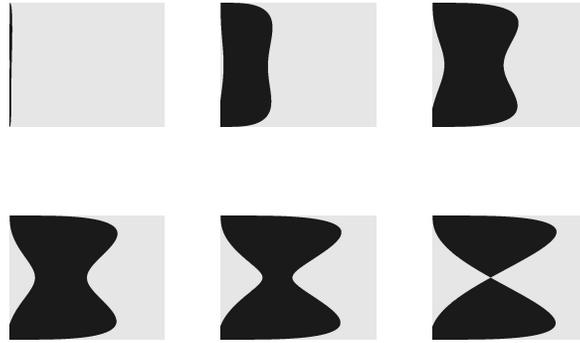}
\caption{$\omega$-degeneracy curves of Lagrangian triangular homographic orbits in the
rectangle ($\beta,e$), $\beta\in(0,2]$ and $e\in(-1,1)$. From the left top to the right bottom,
$\omega$ goes from $1$ to $-1$ along the upper half unit circle of $\U$, more precisely,
$\omega=\exp(\sqrt{-1}\theta)$ with $\theta=\pi/100,\pi/5,2\pi/5, 3\pi/5,4\pi/5,\pi$ respectively.
The last one corresponds to $\omega=-1$. In the black region, the $\omega$-index is $1$.}
\end{figure*}

For $\omega\in\U\bs\{1\}$, we showed that there are two nontrivial degeneracy curves (possibly
tangent to each other at isolated points) $\beta_1(e,\omega)$ and $\beta_2(e,\omega)$ ($\omega=-1$
in Theorem \ref{T1.2} and general $\omega$ in Theorem \ref{T1.5}) which are real analytic in
$e\in [0,1)$. These $\om$ degeneracy curves actually yield a foliation of the non-hyperbolic region
of $\ga_{\bb,e}(2\pi)$ in the $(\bb,e)$ rectangle $[0,9]\times (-1,1)$, when $\om$ runs through $\U$.
We have conducted numerical computations to see this interesting phenomenon according to our above
analysis for $(\bb,e)\in [0,9]\times (-1,1)$. Specially in the Figure 4, we pick up certain figures
from such computations for readers. It is interesting to know how these degeneracy curves behave
under the variation of $\omega$, which we leave for future studies.

In summary, many problems observed numerically already deserve pursuit further.
For example, more precise properties of degeneracy curves including their asymptotic behaviors,
possible intersections and variations with respect to $\omega$, including the above
mentioned interesting properties. In this paper, we have not considered separations between
HH and CS either. We shall study these problems in some forthcoming papers, and we believe that
the ideas and the methods we have developed here can also be used to linear stability problems
for other solutions of the $n$-body problems and systems with periodic coefficients.

\medskip

\noindent {\bf Acknowledgements.} We especially thank Dr. Jie Yang, Yuwei Ou, and Professor Leshun Xu
for helping us to manage the figures and using Matlabs. We would also like to thank valuable
comments and interests of Professors Alain Albouy, Sergey Bolotin, Kuo-Chang Chen, Alain Chenciner
and Susanna Terracini on the earlier manuscript of this paper.

\end{document}